\documentclass[a4paper,11pt]{article}
\usepackage{latexsym,amssymb}
\usepackage{amsmath}
\usepackage[mathscr]{eucal}
\usepackage{xcolor}
\usepackage[abbrev]{amsrefs}

\usepackage{enumerate}

\newtheorem{Theorem}{\bf Theorem}[section]
\newtheorem{Lemma}[Theorem]{\bf Lemma}
\newtheorem{Proposition}[Theorem]{\bf Proposition}
\newtheorem{Corollary}[Theorem]{\bf Corollary}
\newtheorem{Remark}[Theorem]{\bf Remark}
\newtheorem{Example}[Theorem]{\bf Example}
\newtheorem{Definition}[Theorem]{\bf Definition}

\newenvironment{theorem}{\begin{Theorem}$\!\!\!$}{\end{Theorem}}
\newenvironment{lemma}{\begin{Lemma}$\!\!\!$}{\end{Lemma}}
\newenvironment{proposition}{\begin{Proposition}$\!\!\!$}{\end{Proposition}}
\newenvironment{corollary}{\begin{Corollary}$\!\!\!$}{\end{Corollary}}
\newenvironment{remark}{\begin{Remark}$\!\!\!$}{\end{Remark}}

\newenvironment{definition}{\begin{Definition}$\!\!\!$}{\end{Definition}}

\numberwithin{equation}{section}

\numberwithin{equation}{section}

\begin{document}
\title{Characterization of $F$-concavity preserved\\ by the Dirichlet heat flow}
\author{Kazuhiro Ishige, Paolo Salani, and Asuka Takatsu}
\date{}
\maketitle
\vspace{-15pt}
\begin{abstract}
$F$-concavity is a generalization of power concavity and, actually, the largest available generalization of the notion of concavity. 
We characterize the $F$-concavities preserved by the Dirichlet heat flow in convex domains 
on ${\mathbb R}^n$, and 
complete the study of preservation of  concavity properties by the Dirichlet heat flow, started by Brascamp and Lieb in 1976 
and developed in some recent papers. 
More precisely:
\begin{itemize}
  \item[(1)] we discover hot-concavity, which is the strongest $F$-concavity preserved by the Dirichlet heat flow;
  \item[(2)] we show that log-concavity is the weakest $F$-concavity preserved by the Dirichlet heat flow; 
  quasi-concavity is also preserved only for $n=1$;
  \item[(3)] we prove that if $F$-concavity does not coincide with log-concavity and it is not stronger 
  than log-concavity and  $n\ge 2$, 
  then there exists an $F$-concave initial datum such that the corresponding solution to the Dirichlet heat flow is not even quasi-concave, hence losing any reminiscence of concavity.

\end{itemize}
Furthermore, we find a sufficient and necessary condition for $F$-concavity to be preserved by the Dirichlet heat flow. 
We also study the preservation of concavity properties by solutions of the Cauchy--Dirichlet problem 
for linear parabolic equations with variable coefficients
and for nonlinear parabolic equations such as semilinear heat equations, the porous medium equation, and the parabolic $p$-Laplace equation. 
\end{abstract}
\vspace{40pt}

\noindent Addresses:

\smallskip
\noindent 
K. I.: Graduate School of Mathematical Sciences, The University of Tokyo,\\ 
3-8-1 Komaba, Meguro-ku, Tokyo 153-8914, Japan\\
\noindent 
E-mail: {\tt ishige@ms.u-tokyo.ac.jp}\\

\smallskip
\noindent 
P. S.: Dipartimento di Matematica e Informatica``U. Dini", 
Universit\`a di Firenze,\\ viale Morgagni 67/A, 50134 Firenze\\
\noindent 
E-mail: {\tt paolo.salani@unifi.it}\\

\smallskip
\noindent 
A. T.: Department of Mathematical Sciences, Tokyo Metropolitan University,\\
1-1 Minami-osawa, Hachioji-shi, Tokyo 192-0397, Japan\\
\noindent 
E-mail: {\tt asuka@tmu.ac.jp}\\

\vspace{3pt}
\newpage
\tableofcontents
\section{Introduction}
Let $\Omega$ be a (non-empty) convex domain in ${\mathbb R}^n$ with $n\ge 1$ and $\alpha\in[-\infty,\infty]$, throughout this paper.
A nonnegative function $f$ in $\Omega$ is said {\em $\alpha$-concave} in $\Omega$ if 
$$
f((1-\lambda)x+\lambda y)\ge {\mathcal M}_{\alpha}(f(x),f(y);\lambda)
$$
for $x$, $y\in\Omega$ and 
$\lambda\in(0,1)$, where ${\mathcal M}_{\alpha}(a,b;\lambda)$ 
indicates the $\alpha$-mean of two nonnegative numbers $a$ and~$b$ with weight $\lambda\in(0,1)$, defined as follows:
$$
{\mathcal M}_{\alpha}(a,b;\lambda)
:=\left\{\begin{array}{ll}
\max\{a,b\}\quad&\text{if }\alpha=\infty,\vspace{3pt}\\
\left((1-\lambda)a^{\alpha}+\lambda b^{\alpha}\right)^{\frac1\alpha}\quad&\text{if }\alpha\in\mathbb{R}\setminus\{0\},\vspace{5pt}\\
a^{1-\lambda}b^\lambda\quad&\text{if }\alpha=0,\vspace{3pt}\\
\min\{a,b\}\quad&\text{if }\alpha=-\infty,\\
\end{array}\right.
$$
for $a$, $b>0$ and ${\mathcal M}_{\alpha}(a,b;\lambda):=0$ for $a$, $b\ge 0$ with $ab=0$. 
Roughly speaking, we can say that 
a nonnegative function $f$ is $\alpha$-concave
if $\alpha f^\alpha$ is concave  for $\alpha\neq 0$ and if $\log f$ is concave  for $\alpha=0$; 
moreover, it is $\infty$-concave if
it is a positive constant in a convex set  and vanishes elsewhere, 
while it is $(-\infty)$-concave if all the superlevel sets of $f$ are convex.
The case $\alpha=1$ clearly corresponds to the usual concavity 
and $\alpha=0$ corresponds to {\em log-concavity}, 
while the case $\alpha=-\infty$ is usually referred to as {\em quasi-concavity}.
{\em Power concavity} is a generic term for $\alpha$-concavity with $\alpha\in[-\infty,\infty]$.

Due to Jensen's inequality, power concavity has the following nice property: 
\begin{itemize}
  \item
  if $f$ is $\alpha$-concave in $\Omega$ and $\beta\le\alpha$,
  then $f$ is $\beta$-concave in $\Omega$. 
\end{itemize}
This property establishes a hierarchy among power concavities, so that
quasi-concavity is the weakest one while $\infty$-concavity is the strongest one. 

Power concavity is a useful variation of the usual concavity 
and it has been largely studied in the framework of elliptic and parabolic equations (see later for some literature).
Here we are mainly concerned with a classical result by 
Brascamp and Lieb: {\em log-concavity is preserved by the Dirichlet heat flow in~$\mathbb{R}^n$}.
Indeed, in the celebrated paper~\cite{BL},  they proved that
if $f$ is $\alpha$-concave in~${\mathbb R}^{m+n}$ for some $\alpha\in[-1/m,\infty]$, where $m\ge 1$,
then the function
$$
x\mapsto\int_{{\mathbb R}^m}f(x,y)\,dy
$$
is $\gamma$-concave in ${\mathbb R}^n$ with $\gamma=\alpha/(1+m\alpha)$ if $\alpha>-1/m$ and $\gamma=-\infty$ if $\alpha=-1/m$. 
This implies that the solution of
$$
\partial_t u=\Delta u\quad\mbox{in}\quad{\mathbb R}^n\times(0,\infty),
\quad
u(\cdot,0)=\phi\quad\mbox{in}\quad{\mathbb R}^n,
$$
given by
\begin{equation}
\label{eq:1.1}
u(x,t)=(e^{t\Delta_{{\mathbb R}^n}}\phi)(x):=(4\pi t)^{-\frac{n}{2}}\int_{{\mathbb R}^n}e^{-\frac{|x-y|^2}{4t}}\phi(y)\,dy,
\end{equation}
is log-concave in ${\mathbb R}^n$ for all $t>0$ if $\phi$ is log-concave in ${\mathbb R}^n$. 
More generally, 
from the argument of Brascamp and Lieb (see also \cite{Kor} for a different proof),  it can be retrieved that
log-concavity is pushed forward by the Dirichlet heat flow (which is abbreviated as DHF) in any convex domain~$\Omega$ (not only in the whole $\mathbb{R}^n$). 
More precisely,
\begin{itemize}
\item
  if $u$ is a (nonnegative and bounded) solution of 
  \begin{equation}
  \tag{{H}}
  \left\{
  \begin{array}{ll}
  \partial_t u=\Delta u & \quad\mbox{in}\quad\Omega\times(0,\infty),\vspace{3pt}\\
  u=0 & \quad\mbox{on}\quad\partial\Omega\times(0,\infty)\mbox{ if }\partial\Omega\not=\emptyset,\vspace{3pt}\\
  u(\cdot,0)=\phi & \quad\mbox{in}\quad\Omega,
  \end{array}
  \right.
  \end{equation}
  where $\phi\in L^\infty(\Omega)$,
  then $u(\cdot,t)$ is log-concave in $\Omega$ for all $t>0$ provided that 
  $\phi$ is log-concave in $\Omega$.
 \end{itemize}
Throughout the paper we denote by $e^{t\Delta_\Omega}\phi$ the unique (nonnegative and bounded) solution of problem~(H)
(see also the beginning of Section~\ref{section:3}).
\smallskip

The preservation of concavity properties by parabolic equations is an interesting subject of investigation 
with connections and applications to different fields like economy and physics and to other important mathematical questions
in the study of the eigenvalue problems, 
curvature flows, localization of hot spots, and 
functional and geometric inequalities such as Pr\'ekopa--Leindler and Borell--Brascamp--Lieb inequalities. 
Then, in view of the richness of the realm of power concavities, after the result by Brascamp and Lieb, 
it is first of all natural to ask the following question:
\begin{itemize}
\item[(Q)] 
Are there any power concavities preserved by DHF other than log-concavity? 
\end{itemize}
Question (Q) is open in its full generality; there are however partial results, which we recall hereafter.
\begin{proposition}
\label{Theorem:1.1}
Let $\Omega$ be a convex domain in ${\mathbb R}^n$. 
\begin{itemize}
  \item[{\rm (1)}] 
  Let $n=1$. Then $e^{t\Delta_\Omega}\phi$ is quasi-concave in $\Omega$ for all $t>0$ if $\phi$ is quasi-concave in $\Omega$, 
  i.e. quasi-concavity is preserved in dimension $1$.
  \item[{\rm (2)}] Quasi-concavity is in general not preserved in dimension $n\geq 2$: 
  indeed, there exists $\phi\in C_0(\Omega)$ such that   $\phi$ is $\alpha$-concave in~$\Omega$ for some $\alpha\in(-\infty,0)$ and 
  $e^{t\Delta_\Omega}\phi$ is not quasi-concave in~$\Omega$ for some $t>0$.
   \item[{\rm (3)}] Log-concavity is the strongest power concavity which DHF transmits from the initial time to any~$t>0$. 
 \end{itemize}
\end{proposition}
See \cite{An} and \cite{IS01} for assertion~(1); 
see \cite{IS01}*{Theorem~1.1} and \cite{IS02}*{Theorem~4.1} for assertion~(2);  
see \cite{IST03}*{Theorem~1.1} for assertion~(3). 

Proposition~\ref{Theorem:1.1} implies that log-concavity is the strongest power concavity preserved by DHF and that 
quasi-concavity is the weakest power concavity preserved by DHF when $n=1$, while it is not preserved for $n\geq2$.
Thus Proposition~\ref{Theorem:1.1} of course tells us much, but it is far from a  {thorough} answer to question~(Q). 
Furthermore, power concavity is not the only nor the best generalization of concavity. 
Indeed, power concavity is just a particular case of the more general notion of {\em $F$-concavity} (or {\em concavifiability}) and  
$\alpha$-log-concavity (see Example~\ref{Theorem:2.2} for the definition), 
which is a kind of $F$-concavity and not a power concavity, is preserved by DHF in $\Omega$ if $\alpha\in[1/2,1]$ 
(see \cite{IST01}*{Theorem~3.1}).
\begin{definition}
\label{Theorem:1.2}
Let $I=[0,a)$ and $\mbox{{\rm int}}\,I=(0,a)$ with $a\in(0,\infty]$, throughout this paper. 
\begin{itemize}
\item[{\rm (i)}]
A function~$F:I\to[-\infty,\infty)$ is said {\em admissible} on $I$ if $F\in C(\mbox{{\rm int}}\,I)$, 
$F$ is strictly increasing on $I$, and $F(0)=-\infty$. 
Throughout the paper, we denote by $f_F$ the inverse function of $F$ in $J_F:=F({\rm int}\,I)$ 
and set $F(a):=\lim_{r\to a-0}F(r)$.
\item[{\rm (ii)}] 
Let $F$ be admissible on $I$. 
Set 
$$
\mathcal {A}_\Omega(I):=\{f:\Omega\to{\mathbb R}\ |\ f(\Omega) \subset I\}.
$$ 
Given $f\in\mathcal {A}_\Omega(I)$,
we say that $f$ is \emph{$F$-concave} in $\Omega$
if
\[
F(f((1-\lambda)x+\lambda y))
\geq
(1-\lambda)F(f(x))+\lambda F(f(y))
\]
 for all $x,y\in\Omega$ and $\lambda \in (0,1)$. 
 We denote by $\mathcal{C}_\Omega[F]$ the set of $F$-concave functions in $\Omega$.
 \end{itemize}
\end{definition}
In the above definition and throughout this paper, we adhere to the convention that
\begin{align*}
 & -\infty\leq -\infty,\qquad
 \infty\le\infty,\qquad
 -\infty+b
=b-\infty=-\infty, \qquad
\kappa \cdot(\pm\infty)=\pm\infty,\\
 & \,e^{-\infty}=0,\qquad \log 0=-\infty,\qquad -\log 0=\infty,\qquad \log\infty=\infty,
\end{align*}
where $b\in [-\infty,\infty)$ and $\kappa\in(0,\infty)$. 
We retrieve power concavity by considering, for $\alpha\in{\mathbb R}$, the admissible function $\Phi_\alpha$ on $I=[0,\infty)$ defined 
\begin{equation*}
\Phi_{\alpha}(r):=\int_{1}^r s^{\alpha-1}\,ds=
\left\{
\begin{array}{ll}
\displaystyle{\frac{r^{\alpha}-1}{\alpha}} & \mbox{if}\quad \alpha\not=0,\vspace{7pt}\\
\log r& \mbox{if}\quad \alpha=0,\vspace{7pt}\\
\end{array}
\right.
\end{equation*}
for $r\in(0,\infty)$ and $\Phi_{\alpha}(0):=-\infty$. (See also Example~\ref{Theorem:2.2}.)

Clearly, if $F$ is admissible on $[0,a)$, then it is admissible on~$[0,a')$ for every $a'\in(0,a]$. 
Furthermore, for any rigid motion~$T$ on $\mathbb{R}^n$ and $\ell\in(0,\infty)$,  
we have
\begin{equation}
\label{eq:1.2}
f\in {\mathcal C}_{\ell T(\Omega)}[F]\quad\mbox{if and only if}\quad f\circ  (\ell T)\in {\mathcal C}_{\Omega}[F].
\end{equation}
In the universe of $F$-concavities, 
it is possible to introduce a hierarchy, which generalizes the one established among power concavities.
\begin{definition}
\label{Theorem:1.3}
Let $a_1$, $a_2$, $a\in(0,\infty]$ and $a\le\min\{a_1,a_2\}$. 
Set $I_1=[0,a_1)$, $I_2=[0,a_2)$, and $I=[0,a)$. 
Let $F_1$ and $F_2$ be admissible on $I_1$ and $I_2$, respectively. 
We say that {\em $F_1$-concavity is weaker {\rm({\it resp.}\,\,strictly weaker)} than $F_2$-concavity}, 
or equivalently that  {\em $F_2$-concavity is stronger {\rm({\it resp.}\,\,strictly stronger)} than $F_1$-concavity}, in ${\mathcal A}_\Omega(I)$ if 
$$
\mathcal{C}_\Omega[F_2]\cap{\mathcal A}_\Omega(I)\subset\mathcal{C}_\Omega[F_1]
\qquad
\mbox{{\rm (}resp.\,\,$\mathcal{C}_\Omega[F_2]\cap{\mathcal A}_\Omega(I)\subsetneq\mathcal{C}_\Omega[F_1]${\rm)}}.
$$
\end{definition} 
We remark that, in our definition, any $F$-concavity is both stronger and weaker than itself (we use {\em "strictly"} when a strong  comparison applies).
Notice also that, although quasi-concavity does not posses any corresponding admissible function 
(see e.g. \cite{IST03}*{Remark 2.2} and \cite{ConnellRasmusen}), due to the monotonicity of admissible functions, 
an $F$-concave function is always quasi-concave. Similarly, a $\infty$-concave function is $F$-concave for every admissible $F$. 
Then quasi-concavity remains the weakest conceivable concavity property and
$\infty$-concavity is the strongest one. 
We use the expression {\em $\overline{F}$-concavity} when we want to consider all the $F$-concave functions jointly with quasi-concave and $\infty$-concave functions.
\medskip

The main aims of this paper are to strengthen the result by Brascamp and Lieb 
and to investigate its sharpness in the framework of $F$-concavity, 
asking the following questions.
\begin{itemize}
\item[({\bf Q1})] 
What is the strongest $\overline F$-concavity preserved by DHF? 
\item[({\bf Q2})] 
What is the weakest $\overline F$-concavity preserved by DHF? 
\item[({\bf Q3})]
When starting with an $F$-concave initial datum and  $F$-concavity is not preserved by DHF, can we at least hope in maintaining quasi-concavity?
\end{itemize}
For the sake of clarity, let us state explicitly that by saying {\em ``$\overline F$-concavity is preserved by DHF in $\Omega$"}
we mean that, 
\begin{align*}
 & \mbox{if $\phi\in L^\infty(\Omega)$ is $\overline F$-concave in $\Omega$, then
the solution $e^{t\Delta_\Omega}\phi$ of problem~(H)}\\
 & \mbox{is $\overline F$-concave in $\Omega$ for every $t>0$}.
\end{align*}
Due to Proposition~\ref{Theorem:1.1}~(1), 
in the one-dimensional case  {the answers to questions~({\bf Q2}) and ({\bf Q3})
are "Quasi-concavity" and ``Yes", respectively.}
\medskip

We give here complete answers to the above three questions for every $n$. In order to do it, especially for question (Q1), we need to introduce a family of new admissible functions.
\begin{definition}\label{Theorem:1.4}
Let
\begin{equation}
\label{eq:1.3}
h(z):=\left(e^{{\Delta}_{{\mathbb R}}}{\bf 1}_{[0,\infty)}\right)(z)=(4\pi)^{-\frac{1}{2}}\int_0^\infty e^{-\frac{|z-w|^2}{4}}\,dw
\quad\mbox{for $z\in{\mathbb R}$}.
\end{equation}
Then the function $h$ is smooth in ${\mathbb R}$, $\lim_{z\to-\infty}h(z)=0$, $\lim_{z\to\infty}h(z)=1$, and $h'>0$ in ${\mathbb R}$
{\rm ({\it see Lemma}~\ref{Theorem:2.9})}. 
Denote by $H$ the inverse function of $h$. 
For any $a\in(0,\infty]$, we define an admissible function $H_a$ on $[0,a)$ by  
\begin{equation*}
H_a(r):=
\left\{
\begin{array}{ll}
H(r/a) & \mbox{for $r\in(0,a)$ if $a>0$},\vspace{3pt}\\
\log r & \mbox{for $r\in(0,a)$ if $a=\infty$},\vspace{3pt}\\
-\infty & \mbox{for $r=0$ and $a\in(0,\infty]$}.
\end{array}
\right.
\end{equation*}
{\rm ({\it See} Lemma~\ref{Theorem:2.10} {\it for the coherence of this definition})}.
We call $H_a$-concavity {\em hot-concavity}. 
\end{definition}
Hot-concavity and the other already named $F$-concavities can be ordered, and, in particular,  $H_a$-concavity is strictly stronger than log-concavity for $a\in(0,\infty)$ (while they clearly coincide for $a=\infty$).

Now we are ready to give our answers to questions~({\bf Q1})-({\bf Q3}).
\begin{theorem}\label{mainthm}
\label{MainThm}
Let $I=[0,a)$ with $a\in(0,\infty]$ and  $\Omega$ a convex domain in ${\mathbb R}^n$ with $n\ge 1$.  Then the following properties hold.
\begin{itemize}
  \item[{\rm({\bf A1})}] 
  $H_a$-concavity is the strongest $F$-concavity preserved by DHF in~${\mathcal A}_\Omega(I)$.
\item[{\rm ({\bf A2})}] 
  Log-concavity is the weakest $F$-concavity preserved by DHF in~${\mathcal A}_\Omega(I)$ for $n\geq 2$; when $n=1$, the same is true under $C^2$-regularity assumption on $F$ 
  {\rm (}moreover, quasi-concavity is preserved, too{\rm )}.
 \item[{\rm ({\bf A3})}] 
  If $F$-concavity is not stronger than log-concavity in~${\mathcal A}_\Omega(I)$ and  $n\ge 2$, 
  then there exists an $F$-concave initial datum such that the corresponding solution to the Dirichlet heat flow is not even quasi-concave, hence losing any reminiscence of concavity.
%
\end{itemize}
\end{theorem}

With the exception of quasi-concavity for $n=1$, the situation depicted by the above theorem can be nicely summarized by the following picture:
{\small
$$
  \mbox{$\infty$-concavity}\subset\dots\subset
  \overbrace{
  \mbox{$H_a$-concavity}\subset\dots\subset
  \mbox{log-concavity}=\mbox{$H_\infty$-concavity}}^{\mbox{stronger $\leftarrow$\quad\qquad Preserved by DHF\qquad\quad$\rightarrow$ weaker}}\subset\cdots\subset\mbox{quasi-concavity}
$$
}
Let $\Pi$ is the half-space $[0,\infty)\times {\mathbb R}^{n-1}$
of ${\mathbb R}^n$. Since 
\[
h(z)=e^{\Delta_{{\mathbb R}^n}}{\bf 1}_\Pi(ze_1) \quad\mbox{for $z\in{\mathbb R}$}
\]
and the characteristic function 
${\bf 1}_\Pi$ is $\infty$-concave in ${\mathbb R}^n$, 
if $F$ is admissible in $I=[0,1)$ and $F$-concavity is preserved by DHF in ${\mathbb R}^n$, then $h$ must be $F$-concave in ${\mathbb R}$. 
This implies that $H_1$-concavity is stronger than $F$-concavity in ${\mathcal A}_{{\mathbb R}}(I)$ (see Lemma~\ref{Theorem:2.4}). 
Then $H_1$-concavity is the possible strongest $F$-concavity preserved by DHF and hot-concavity naturally appears in the study of the preservation of concavity properties by DHF. 
Fortunately, hot-concavity is preserved by DHF, and thus ({\bf A1}) holds. 
Notice also that answer ({\bf A1}) depends on the interval $I$, precisely on $a$.  The dependance on $a$ can be interpreted as dependance 
on the $L^\infty$ norm of the initial datum. 
In the case of $a=\infty$, it has already been shown in \cite{IST03}*{Theorem~1.1}
that log-concavity (i.e. $H_\infty$-concavity) is the strongest $F$-concavity  preserved by DHF 
in~${\mathcal A}_\Omega(I)$. 
In fact, by coupling ({\bf A1}) and ({\bf A2}),
we get that in the case $a=\infty$ the only $F$-concavity surely preserved by DHF is log-concavity. 
\medskip

We will split Theorem \ref{MainThm} in several steps, see Section 5.
Let us notice here that one of the key ingredients in the proof of one of these steps 
is a result concerning a sufficient condition for the preservation of $F$-concavity by DHF 
(see Proposition~\ref{Theorem:4.4} and Theorem~\ref{Theorem:1.5}). 
Surprisingly, the sufficient condition, under a $C^2$-regularity assumption on admissible functions, is also necessary! 
This takes to another main result of this paper.
\begin{theorem}
\label{Theorem:1.8}
Let $I=[0,a)$ with $a\in(0,\infty]$ and $\Omega$ a convex domain in $\mathbb{R}^n$ with $n\geq 1$. 
Let $F$ be admissible on $I$ such that $F\in C^2(\mbox{{\rm int}}\,I)$. 
Then  $F$-concavity is preserved by DHF in $\Omega$ if and only if
  $$ 
  \lim_{r\to +0}F(r)=-\infty,\quad \mbox{$F'>0$ in ${\rm int}\,I$},\quad\mbox{and}\quad \mbox{$(\log f_F')'$ is concave in $J_F$}.
  $$
\end{theorem}
The latter theorem (together with Proposition \ref{Theorem:1.1}~(1)) gives the one-dimensional part of ({\bf A2}), 
and it also yields that only log-concavity and quasi-concavity are preserved by DHF with $n=1$ among power concavities.

Let us remark that, to our knowledge, Theorem~\ref{Theorem:1.8} is the first result regarding a {\em necessary} and sufficient condition 
for concavity properties of solutions to partial differential equations. 
\medskip

Further than in the case of DHF, the preservation of log-concavity has already been studied in the case of the Cauchy--Dirichlet problem 
for various nonlinear parabolic equations (see e.g. \cites{GK, INS, Keady, Ken, Kor, Kol, L, LV}). 
See also \cites{Kawohl01, Kawohl, Kawohl02} for related topics. 
Moreover, similar investigation about the  preservation or disruption of power concavity along parabolic flows have been pursued in the following cases too: 
the heat equation with a potential \cite{AI}; the one-phase Stefan problem~\cite{CW2}, 
where not even log-concavity is in general preserved; porous medium equation  \cite{CW3} and \cite{IS02}
(with sharp results in some cases); DHF in ring shaped domains \cite{CW1}. 
In this paper, as an application and a generalization of our arguments, 
we study the preservation and the disruption of $F$-concavity by solutions of linear parabolic equations with variable coefficients 
and nonlinear parabolic equations such as semilinear heat equations, the porous medium equation, and the parabolic $p$-Laplace equation, 
and also resolve some related open questions 
(see Sections~\ref{section:4} and~\ref{section:6}).
\vspace{3pt}

The rest of this paper is organized as follows.
In Section~\ref{section:2} we discuss the notion of $F$-concavity, giving examples and collecting some preliminary properties.
In Section~\ref{section:3} 
we first study the disruption of quasi-concavity in ${\mathbb R}^n$ with $n\geq2$, 
proving that starting with a non log-concave initial datum even quasi-concavity can go immediately lost (see Proposition~\ref{Theorem:3.2}). 
Next, we show that 
the disruption of $F$-concavity by DHF in ${\mathbb R}^n$ implies the disruption 
of $F$-concavity by the Dirichlet parabolic flow in any convex domain $\Omega$. 
In Section~\ref{section:4} we investigate the preservation of $F$-concavity 
by generic parabolic flows, finding sufficient (see Proposition~\ref{Theorem:4.2}) and necessary (see Proposition~\ref{Theorem:4.4}) conditions. 
Furthermore, we characterize $F$-concavities preserved by DHF. 
In Section~\ref{section:5} we complete the proof of Theorem \ref{MainThm}, which will be in fact the product of other main theorems stated and proved in this section.
In Section~\ref{section:6} we develop the arguments of Sections~3--5, and we explore the $F$-concavities preserved by solutions of 
linear parabolic equations with variable coefficients and nonlinear parabolic equations. 
\vspace{5pt}

\noindent
{\bf Acknowledgements.} 
K. I. and A. T. were supported in part by JSPS KAKENHI Grant Number 19H05599. 
P. S. was supported in part by  INdAM through a GNAMPA Project.
A. T. was supported in part by JSPS KAKENHI Grant Number 19K03494.
\section{Preliminaries}\label{section:2}
Throughout this paper, for any $x\in{\mathbb R}^n$ and $R>0$, we denote by $B(x,R)$ the open ball in~${\mathbb R}^n$ centered at $x$ of radius $R$. 
Let  $\langle \cdot, \cdot \rangle $ denote the standard inner product on $\mathbb{R}^n$. 
For any set~$E$, let~${\bf 1}_E$ be the characteristic function of $E$. 
We denote by $\mathrm{Sym}\,(n)$ the space of $n\times n$ real symmetric matrices. 
As already said, unless otherwise stated, we denote by $\Omega$ 
a (non-empty) convex domain in ${\mathbb R}^n$.
Let
$$
BC(\Omega):=C(\Omega)\cap L^\infty(\Omega),
\quad
BC_0(\overline{\Omega}):=
\{f\in C(\overline{\Omega})\cap L^\infty(\Omega)\,|\, \mbox{$f=0$ on $\partial\Omega$ if $\partial\Omega\not=\emptyset$}\}. 
$$
For any $T\in(0,\infty]$, we define
\begin{align*}
 & BC(\Omega\times(0,T)):=C(\Omega\times(0,T))\cap L^\infty(\Omega\times(0,T)),\\
 & C^{2;1}(\Omega\times(0,T))\\
 & :=\{f\in C(\Omega\times(0,T))\,|\,\partial_{x_i} f,\partial_{x_i}\partial_{x_j} f, \partial_t f\in C(\Omega\times(0,T))
\,\,\,\mbox{for}\,\,\, i, j=1,\dots,n\},\\
 & BC^{2;1}(\Omega\times(0,T))\\
 & :=\{f\in BC(\Omega\times(0,T))\,|\,\partial_{x_i} f,\partial_{x_i}\partial_{x_j} f, \partial_t f\in BC(\Omega\times(0,T))
\,\,\,\mbox{for}\,\,\, i, j=1,\dots,n\}.
\end{align*}
For any $T\in(0,\infty]$ and $\sigma\in[0,1)$, we set 
\begin{align*}
 & BC^{0,\sigma;0,\sigma/2}(\Omega\times(0,T)):=\{f\in BC(\Omega\times(0,T))\,|\,\|f\|_{C^{0,\sigma;0,\sigma/2}(\Omega\times(0,T))}<\infty\},\\
 & BC^{2,\sigma;1,\sigma/2}(\Omega\times(0,T))\\
 & :=\{f\in BC^{2;1}(\Omega\times(0,T))\,|\,\partial_{x_i}\partial_{x_j}f, \partial_t f\in BC^{0,\sigma;0,\sigma/2}(\Omega\times(0,T)) <\infty\mbox{ for $i$, $j=1,\dots,n$}\},
\end{align*}
where
$$
\|f\|_{C^{0,\sigma;\,0,\sigma/2}(\Omega\times(0,T))}:=
\sup_{(x,t)\in\Omega\times(0,T)}|f(x,t)|
+\sup_{\substack{(x,t), (y,s)\in\Omega\times(0,T)\\ (x,t)\not=(y,s)}}
\frac{|f(x,t)-f(y,s)|}{|x-y|^\sigma+|t-s|^{\sigma/2}}<\infty.
$$ 
Similarly, we define the function spaces $C^{k,\sigma}(\Omega)$, $BC^{k,\sigma}(\Omega)$, $BC^{k,\sigma}(\overline{\Omega})$,
 and $C^{k,\sigma;k,\sigma/2}(\Omega\times(0,T))$, where $k=0,1,2$ and $\sigma\in[0,1)$. 
We use $C$ to denote generic positive constants,  {which}
may take different values within a calculation. 

Let us first give three relevant examples of $F$-concavity.
As we said, power concavity is just a particular case of $F$-concavity, hence this is our first example.
\begin{Example} {\rm (Power concavity)}
\label{Theorem:2.1}
Let $I=[0,\infty)$ and $\alpha\in{\mathbb R}$. 
Define an admissible function $\Phi_\alpha$ on $I$ by  
\begin{equation*}
\Phi_{\alpha}(r):=\int_{1}^r s^{\alpha-1}\,ds=
\left\{
\begin{array}{ll}
\displaystyle{\frac{r^{\alpha}-1}{\alpha}} & \mbox{if}\quad \alpha\not=0,\vspace{7pt}\\
\log r& \mbox{if}\quad \alpha=0,\vspace{7pt}\\
\end{array}
\right.
\end{equation*}
for $r\in(0,\infty)$ and $\Phi_{\alpha}(0):=-\infty$. Then 
$\Phi_\alpha$-concavity corresponds to $\alpha$-concavity, as introduced at the beginning of this paper, and it possesses the following properties. 
\begin{itemize}
  \item[{\rm (1)}] 
  If $\alpha<\beta$, then $\alpha$-concavity is strictly weaker than $\beta$-concavity in ${\mathcal A}_\Omega(I)$.
  \item[{\rm (2)}] 
  Power concavity is closed under positive scalar multiplication, that is, 
  if $f$ is $\alpha$-concave in~$\Omega$, then so is $\kappa f$ for $\kappa\in(0,\infty)$. 
  {\rm ({\it See also Lemma~{\rm\ref{Theorem:2.8}}}.)}
 \end{itemize}
 \end{Example}

The second example is a sort of hybrid between log-concavity and power concavity, introduced in \cite{IST01}.
\begin{Example}
\label{Theorem:2.2} {\rm (Power log-concavity)} 
Let $I=[0,1)$ and $\alpha \in {\mathbb R}$. 
Define an admissible function~$L_\alpha$ on $I$ by
$$
L_\alpha(r):=
-\Phi_\alpha(-\log r)=
\left\{
\begin{array}{ll}
-\displaystyle{\frac{1}{\alpha}}\left[(-\log r)^\alpha-1\right] & \mbox{if}\quad\alpha\not=0,\vspace{7pt}\\
-\log(-\log r) & \mbox{if}\quad\alpha=0,
\end{array}
\right.
$$
for $r\in(0,1)$ and $L_{\alpha}(0):=-\infty$. 
We also refer to $L_\alpha$-concavity as {\rm $\alpha$-log-concavity}
and, generically, as {\rm power log-concavity}. 
The following properties hold {\rm(}see \cites{IST01, IST03}{\rm)}.
\begin{itemize}
 \item[{\rm (1)}] 
 A function $f\in{\mathcal A}_\Omega(I)$ is log-concave in $\Omega$ 
 if and only if $f$ is $1$-log-concave in~$\Omega$.
 \item[{\rm (2)}] 
 If $\alpha<\beta$, then $\beta$-log-concavity is strictly weaker than $\alpha$-log-concavity in ${\mathcal A}_\Omega(I)$.
 \item[{\rm (3)}] 
 If $\alpha\le 1$ and $f$ is $\alpha$-log-concave in $\Omega$, 
 then so is $\kappa f$ for $\kappa\in(0,1]$.
\end{itemize}
Also notice that, for any $\kappa\in(0,1)$, 
the function $\kappa\exp(-|x|^2)$ is $\alpha$-log-concave in ${\mathbb R}^n$ if and only if $\alpha\geq1/2$. 
\end{Example}
Power log-concavity plays an important role in the study of DHF. 
Indeed, 
in \cite{IST01}*{Theorem~3.1}, the authors of this paper proved that $\alpha$-log-concavity is preserved by DHF if $\alpha\in[1/2,1]$ 
and this result is optimal in some suitable sense (see and \cite{IST03}*{Section~4.2} for more details).
Notice that the domain of $\alpha$-log-concavity is different from that of log-concavity and 
$\alpha$-log-concavity is strictly stronger than log-concavity in ${\mathcal A}_\Omega([0,1))$ if $\alpha<1$. 
These suggest that the domains of admissible functions are important in the study for questions~({\bf Q1})--({\bf Q3})
and this makes the study complex and delicate. 

The third example is $H_a$-concavity introduced in Definition~{\rm\ref{Theorem:1.4}}  {and it is crucial to this paper}.
\begin{Example}
\label{Theorem:2.3} {\rm ($H_a$-concavity)} 
Let $I=[0,a)$ with $a\in(0,\infty]$.  
$H_a$-concavity possesses the following properties. 
\begin{itemize}
  \item[{\rm (1)}] 
  For any $b\in(a,\infty]$, 
  $H_b$-concavity is strictly weaker than $H_a$-concavity in ${\mathcal A}_\Omega(I)$.
  \item[{\rm (2)}] 
  $H_1$-concavity is stronger than $\alpha$-log-concavity in ${\mathcal A}_\Omega(I)$ if $\alpha\ge 1/2$.
\end{itemize}
Assertion~{\rm (1)} follows from Theorem~{\rm\ref{Theorem:1.5}}. 
Assertion~{\rm (2)} also follows from Theorem~{\rm\ref{Theorem:1.5}} and {\rm\cite{IST01}*{Theorem~3.1}}.
\end{Example}

In the case $I=[0,1)$, as a result of Theorem~\ref{MainThm} with Example~\ref{Theorem:2.2}
(see also Theorems~\ref{Theorem:1.5} and \ref{Theorem:1.6}), 
we can deduce the following picture of the hierarchy of $F$-concavities preserved by DHF.
{\small
$$
  \cdots\subset
  \overbrace{
  \mbox{$H_1$-concavity}\subset 
  \underbrace{\mbox{$1/2$-log-concavity}\subset\cdots\subset \mbox{$1$-log-concavity}}_
  {\mbox{$\alpha$-log-concavity $(1/2\le\alpha\le 1)$}}
  \subset\mbox{log-concavity}=\mbox{$H_\infty$-concavity}}^{\mbox{stronger $\leftarrow$\quad\qquad Preserved by DHF\qquad\quad$\rightarrow$ weaker}}\subset\cdots
$$
}

Next, we collect and prove some lemmas on $F$-concavity. 
The first lemma concerns the hierarchy of $F$-concavities
(see also \cite{AZ}*{Lemma~3.2}).
\begin{lemma}
\label{Theorem:2.4}
Let $F_1$ and $F_2$ be admissible on $I=[0,a)$ with $a\in(0,\infty]$.
Then $F_1$-concavity is weaker than $F_2$-concavity in ${\mathcal A}_\Omega(I)$
if and only if $F_1\circ f_{F_2} $ is concave in $J_{F_2}$  {\rm ({\it or, equivalently, 
$F_2\circ f_{F_1} $ is convex in $J_{F_1}$})}.
\end{lemma}
{\bf Proof.}
For $i=1,2$, we write $f_i:=f_{F_i}$ and $J_i:=J_{F_i}=F_i({\rm int}\,I)$ for simplicity. 
Assume that $F_1\circ f_2 $ is concave in $J_2$. 
Let $f\in{\mathcal C}_\Omega[F_2]$. 
Since $F_2\circ f$ is concave in $\Omega$ and $F_1 \circ f_2$ is increasing on $J_2$, 
we have 
\begin{equation*}
\begin{split}
\left( F_1 \circ f\right) \left(  (1-\lambda)x+\lambda y \right) 
&=
\left( F_1 \circ f_2 \circ F_2 \circ   f\right) \left(  (1-\lambda)x+\lambda y \right) \\
& \geq 
\left( F_1 \circ f_2\right) \left(  (1-\lambda) ( F_2   \circ f) (x) + \lambda ( F_2   \circ f) (y) \right)\\
& \geq (1-\lambda)
\left( F_1 \circ f_2 \right) \left(\left( F_2   \circ f\right) (x) \right)+ \lambda 
\left( F_1 \circ f_2 \right) \left(\left( F_2   \circ f\right) (y) \right) \\
& = (1-\lambda)
\left( F_1  \circ f\right) (x) + \lambda 
\left( F_1  \circ f\right) (y)
\end{split}
\end{equation*}
for $x$, $y\in\Omega$ and $\lambda\in(0,1)$ if $f(x)f(y)>0$. 
This relation also holds even if  $f(x)f(y)=0$ since $F_1(0)=-\infty$ 
(see Definition~\ref{Theorem:1.2}~(i)).  
These imply that $f$ is $F_1$-concave in $\Omega$,
that is, 
$F_1$-concavity is weaker than $F_2$-concavity in ${\mathcal A}_\Omega(I)$.

Next, we assume that $F_1$-concavity is weaker than $F_2$-concavity in ${\mathcal A}_\Omega(I)$. 
If $F_1\circ f_2$ is not concave in $J_2$, then there exist $z$, $w\in J_2$  and $\lambda\in(0,1)$ such that 
\begin{equation}
\label{eq:2.1}
(F_1\circ f_2) \left( (1-\lambda)z+\lambda w) \right)<(1-\lambda) (F_1\circ f_2)  (z) +\lambda (F_1\circ f_2) (w). 
\end{equation}
Thanks to \eqref{eq:1.2},
we can assume, without loss of generality, that $B(0,R)\subset\Omega$ for some $R>0$. 
Let $\varepsilon\in(0,1)$ be such that $\varepsilon z,  \varepsilon w \in(-R,R)$.
Set 
\[
f(x):=f_2(\varepsilon^{-1}\langle x, e_1\rangle){\bf 1}_{J_2}(\varepsilon^{-1} \langle x, e_1\rangle)
\quad\mbox{for}\quad x\in \Omega,
\]
where $e_1:=(1,0,\dots,0)\in{\mathbb R}^n$. 
We observe from $F_2$-concavity of $f_2$ that $f$ is $F_2$-concave in $\Omega$. 
Since $\varepsilon ze_1$, $\varepsilon we_1\in\Omega$, 
by \eqref{eq:2.1} we have 
\begin{align*}
(F_1 \circ f) \left(  (1-\lambda) \varepsilon z e_1 +\lambda \varepsilon w e_1\right) 
&=
(F_1 \circ f_2) \left(  (1-\lambda) z +\lambda w \right) \\
& <(1-\lambda)(F_1 \circ f_2) (z)+\lambda (F_1 \circ f_2) (w)\\
&=
(1-\lambda)(F_1 \circ f) \left(\varepsilon z e_1 \right)
+
\lambda(F_1 \circ f) \left( \varepsilon w e_1 \right).
\end{align*}
This means that $f$ is not $F_1$-concave in $\Omega$, 
which contradicts ${\mathcal C}_\Omega[F_2]\subset {\mathcal C}_\Omega[F_1]$. 
Thus $F_1 \circ f_2$ is concave in $J_2$. 
The proof is complete.  
$\Box$
\vspace{5pt}
\newline
As a corollary of Lemma~\ref{Theorem:2.4}, 
we have (see also \cite{IST03}*{Theorem~3.2}): 
\begin{lemma}
\label{Theorem:2.5}
Let $F_1$ and $F_2$ be admissible on $I=[0,a)$ with $a\in(0,\infty]$. 
Then the relation ${\mathcal C}_\Omega[F_1]={\mathcal C}_\Omega[F_2]$ holds if and only if 
there exists a pair $(A,B)\in (0,\infty)\times{\mathbb R}$ such that 
$$
F_1(r)=AF_2(r)+B\quad\mbox{for}\quad r\in {\rm int}\,I.
$$ 
\end{lemma}
Lemma~\ref{Theorem:2.4} also implies that 
the hierarchy of $F$-concavities is independent of $\Omega$ and the dimension~$n$. 
\begin{lemma}
\label{Theorem:2.6}
Let $F_1$ and $F_2$ be admissible on $I=[0,a)$ with $a\in (0,\infty]$. 
Then 
\begin{equation*}
\begin{split}
 & {\mathcal C}_\Omega[F_2]\subset{\mathcal C}_\Omega[F_1]
\quad\mbox{if and only if}\quad
{\mathcal C}_{{\mathbb R}}[F_2]\subset{\mathcal C}_{{\mathbb R}}[F_1],\\
 & {\mathcal C}_\Omega[F_2]\subsetneq{\mathcal C}_\Omega[F_1]
\quad\mbox{if and only if}\quad
{\mathcal C}_{{\mathbb R}}[F_2]\subsetneq{\mathcal C}_{{\mathbb R}}[F_1].
\end{split}
\end{equation*}
\end{lemma}
Furthermore, by Lemma~\ref{Theorem:2.4} we have: 
\begin{lemma}
\label{Theorem:2.7}
Let $F_1$ and $F_2$ be admissible on $I=[0,a)$ with $a\in(0,\infty]$. 
If $F_1$-concavity is weaker than $F_2$-concavity in ${\mathcal A}_\Omega(I)$ 
and $\lim_{r\to +0}F_2(r)=-\infty$, 
then $\lim_{r\to +0}F_1(r)=-\infty$. 
\end{lemma}
{\bf Proof.}
It follows from the admissibility of $F_2$ and $\lim_{r\to +0}F_2(r)=-\infty$ that $J_{F_2}=(-\infty,F_2(a))$. 
Lemma~\ref{Theorem:2.4} implies that $F_1 \circ f_{F_2}$ is concave in $J_{F_2}$.
Since $F_1$ and $f_{F_2}$ are strictly increasing, 
$F_1 \circ f_{F_2}$ is strictly increasing in $J_{F_2}$. 
Then we conclude that
\[
\lim_{r\to +0}F_1(r)=
\lim_{z\to -\infty}F_1(f_{F_2}(z))=-\infty.
\]
The proof is complete. 
$\Box$
\vspace{5pt}

Next, we modify the argument in the proof of \cite{Lindberg}*{Theorem~2} in order to 
investigate the closedness under positive scalar multiplication of $F$-concavity.
 (See also \cite{IST03}*{Theorem~3.3}.) 
\begin{lemma}
\label{Theorem:2.8}
Let $F$ be admissible on $I=[0,a)$ with $a\in(0,\infty]$. 
Assume that there exists $\varepsilon_*>0$ such that 
the following holds for every $\varepsilon\in(0,\varepsilon_*]$:
if $f\in{\mathcal C}_\Omega[F]$ and $(1+\varepsilon)f\in{\mathcal A}_\Omega(I)$, then 
$\kappa f\in{\mathcal C}_\Omega[F]$ holds for $\kappa\in[(1+\varepsilon)^{-1},1+\varepsilon]$. 
Then ${\mathcal C}_\Omega[F]={\mathcal C}_\Omega[\Phi_\alpha]\cap{\mathcal A}_\Omega(I)$ for some $\alpha\in{\mathbb R}$. 
\end{lemma}
{\bf Proof.}
For $\kappa>0$, set
$F_\kappa(r):=F(\kappa^{-1} r)$ for $r\in  I_\kappa:=[0, \kappa a)$. 
Let $\varepsilon\in(0,\varepsilon_*]$.
Then the assumption implies that
\[
\mathcal{C}_{\Omega}[F_\kappa] \cap\mathcal {A}_\Omega(I_{\kappa/(1+\varepsilon)})
\subset \mathcal{C}_{\Omega}[F] 
\quad \text{for all\ } \kappa\in [(1+\varepsilon)^{-1}, 1+\varepsilon].
\]
Let $\kappa\in [(1+\varepsilon)^{-1/2}, (1+\varepsilon)^{1/2}]$. 
If  $f\in \mathcal{C}_{\Omega}[F] \cap\mathcal {A}_\Omega(I_{\kappa^2/(1+\varepsilon)})$, 
then the assumption yields $\kappa^{-1} f\in \mathcal{C}_{\Omega}[F]  \cap\mathcal {A}_\Omega(I_{\kappa/(1+\varepsilon)})$, that is, 
$f \in  \mathcal{C}_{\Omega}[F_{\kappa}]  \cap\mathcal {A}_\Omega(I_{\kappa^{2}/(1+\varepsilon)})$. 
Consequently, we have 
\[
\mathcal{C}_{\Omega}[F] \cap\mathcal {A}_\Omega(I_{\kappa^2/(1+\varepsilon)})
\subset 
 \mathcal{C}_{\Omega}[F_{\kappa}] \cap \mathcal {A}_\Omega(I_{\kappa^2/(1+\varepsilon)})
 \subset
  \mathcal{C}_{\Omega}[F] \cap\mathcal {A}_\Omega(I_{\kappa^2/(1+\varepsilon)}).
 \]
By Lemma \ref{Theorem:2.5} we see that $F\circ f_{F_\kappa}$ is affine.
Therefore there exist $c(\kappa), d(\kappa)\in \mathbb{R}$ such that 
\[
F (f_{F_\kappa}(z))=c(\kappa)z+d(\kappa)\quad \text{for\ }z \in F_\kappa({\rm int}\,I_{\kappa^2/(1+\varepsilon)}).
\]
In particular, $c(1)=1$, $c(\kappa)>0$, and $d(1)=0$.
Notice that
$$
J_\varepsilon:=\bigcap_{\kappa \in [(1+\varepsilon)^{-1/2}, (1+\varepsilon)^{1/2}]} F_\kappa({\rm int}\,I_{\kappa^2/(1+\varepsilon)})
=F( (0, (1+\varepsilon)^{-3/2}a)). 
$$
Given $z\in J_{\varepsilon}$, 
$F  \circ f_{F_\kappa}$ is continuous with respect to 
$\kappa\in [(1+\varepsilon)^{-1/2}, (1+\varepsilon)^{1/2}]$, 
hence $c$ and~$d$ are continuous in $[(1+\varepsilon)^{-1/2}, (1+\varepsilon)^{1/2}]$. 
Since $f_F$ is the inverse function of $F$ and $f_{F_\kappa}=\kappa f_F$, 
we observe that 
$$
\kappa f_F (z)=f_{F_\kappa}(z)=f_F( F (f_{F_\kappa}(z)))=f_F(c(\kappa)z+d(\kappa))
$$
for all $\kappa\in[(1+\varepsilon)^{-1/2}, (1+\varepsilon)^{1/2}]$. 

Let $\kappa_1, \kappa_2\in [(1+\varepsilon)^{-1/4}, (1+\varepsilon)^{1/4}]$.  
Then $\kappa_1  \kappa_2\in [(1+\varepsilon)^{-1/2}, (1+\varepsilon)^{1/2}]$ and 
\begin{align*}
\kappa_1 \kappa_2 f_F(z)
&=f_F(c(\kappa_1\kappa_2)z +d(\kappa_1\kappa_2) ) ,\\
\kappa_1 \kappa_2 f_F(z)
&=\kappa_1 f_F \left( c(\kappa_2)z  +d(\kappa_2) \right)
=f_F \left( c(\kappa_1) \left(  c(\kappa_2)z  +d(\kappa_2) \right)+d(\kappa_1) \right).
\end{align*}
These yield 
\begin{equation}
\label{eq:2.2}
 c(\kappa_1 \kappa_2)=c(\kappa_1) c(\kappa_2),
\qquad
d(\kappa_1 \kappa_2)=c(\kappa_1) d(\kappa_2)+d(\kappa_1).
\end{equation}
Setting $C(Z):=\log c(e^Z)$, we have 
\[
C(Z_1+Z_2)=\log |c(e^{Z_1}e^{Z_2})|=\log |c(e^{Z_1}) c(e^{Z_2})|=C(Z_1)+C(Z_2)
\]
for $Z_1, Z_2\in [\log(1+\varepsilon)^{-1/4}, \log(1+\varepsilon)^{1/4}]$.
Thus we find  $\alpha\in \mathbb{R}$ such that 
$C(Z)=\alpha Z $,
that is,
\[
c(\kappa)=\kappa^{\alpha} \quad \text{for}\,\,\,\kappa\in[(1+\varepsilon)^{-1/4}, (1+\varepsilon)^{1/4}].
\]

Consider the case of $\alpha\not=0$. 
Thanks to \eqref{eq:2.2}, for any $\kappa_1$, $\kappa_2\in [(1+\varepsilon)^{-1/4}, (1+\varepsilon)^{1/4}]\setminus\{1\}$, 
we have 
\[
d(\kappa_1 \kappa_2)
=\kappa_1^{\alpha} d(\kappa_2) +d(\kappa_1)
=\kappa_2^{\alpha} d(\kappa_1) + d(\kappa_2), 
\]
hence 
$$
\frac{d(\kappa_1)}{\kappa_1^{\alpha}-1}=\frac{d(\kappa_2)}{\kappa_2^{\alpha}-1}.
$$
Combining this with the continuity of $d$, we find
$\beta\in \mathbb{R}$ such that 
\[
d(\kappa)=\beta \Phi_{\alpha}(\kappa)\quad \text{for\ }\kappa\in[(1+\varepsilon)^{-1/4}, (1+\varepsilon)^{1/4}]. 
\]
Let $z\in J_\varepsilon$, and set
$$
\mathcal{I}_\varepsilon(z):=
\left\{\kappa  f_F(z)\,\Bigr|\,\kappa\in[(1+\varepsilon)^{-1/4}, (1+\varepsilon)^{1/4}]\right\}. 
$$
For any $r=\kappa  f_F(z) \in \mathcal{I}_\varepsilon(z)$, we have 
\begin{align*}
F(r) & =F(\kappa f_F(z))= F(f_{F_{\kappa}}(z))
=c(\kappa)z+d(\kappa)=\kappa^{\alpha} z+\beta \Phi_{\alpha}(\kappa)\\
 & =\left(\frac{r}{f_F(z)}\right)^{\alpha} z+\beta \Phi_{\alpha}\left( \frac{r}{f_F(z)} \right)
=\frac{z}{f_F(z)^\alpha}r^\alpha +\frac{\beta}{\alpha}\left(\frac{1}{f_F(z)^\alpha}r^\alpha-1\right).
\end{align*}
This together with the admissibility of $F$ ensures 
the existence of $(A, B)\in (0,\infty) \times \mathbb{R}$ such that 
$F=A\Phi_\alpha+B$ on $\mathcal{I}_{\varepsilon}(z)$.
This relation implies that $A$ and $B$ are independent of the choice of $z\in J_\varepsilon$, that is,
$$
F=A\Phi_\alpha+B \quad \text{on}\quad
\mathcal{I}_{\varepsilon}:=\bigcup_{z\in J_\varepsilon} \mathcal{I}_{\varepsilon}(z).
$$
Since $J_\varepsilon\to F((0,a))$ as $\varepsilon\to +0$, we see that $\mathcal{I}_{\varepsilon}\to (0,a)$ as $\varepsilon\to +0$. 
Then we deduce that  
$F=A\Phi_\alpha+B$ on $(0,a)$. 
This together with Lemma~\ref{Theorem:2.5} completes the proof of Lemma~\ref{Theorem:2.8} in the case $\alpha\not=0$. 

Consider the case of $\alpha=0$. 
It follows from \eqref{eq:2.2} that  $d(\kappa_1\kappa_2)=d(\kappa_1)+d(\kappa_2)$ 
for $\kappa\in[(1+\varepsilon)^{-1/4}, (1+\varepsilon)^{1/4}]$. 
Then we find $\gamma\in{\mathbb R}$ such that $d(\kappa)=\gamma\log\kappa$ for $\kappa\in[(1+\varepsilon)^{-1/4}, (1+\varepsilon)^{1/4}]$, 
that is,
$$
F(r)=c(\kappa)z+d(\kappa)=z+\gamma\log\kappa
=\gamma\log r+z-\log f_F(z)
$$
for $z\in J_\varepsilon$ and $r=\kappa  f_F(z) \in \mathcal{I}_\varepsilon(z)$. 
By a similar argument in the case $\alpha\not=0$ 
we obtain the desired conclusion in the case $\alpha=0$. 
Thus Lemma~\ref{Theorem:2.8} follows.  
$\Box$\vspace{5pt}

At the end of this section 
we prove two lemmas related with $H_a$-concavity. 
In Lemma~\ref{Theorem:2.10} 
we show that log-concave functions are  approximated by $H_a$-concave functions, 
so justifying the  definition of $H_\infty(r)=\log r$. 
\begin{lemma}
\label{Theorem:2.9}
Let $h$ be as in \eqref{eq:1.3}.
Then 
\begin{equation}
\label{eq:2.3}
\left(e^{t\Delta_{\mathbb{R}}}{\bf 1}_{[0,\infty)}\right)(z)=h(t^{-\frac{1}{2}}z)
\quad\mbox{for}\quad (z,t)\in\mathbb{R}\times(0,\infty).
\end{equation}
Furthermore,
\begin{eqnarray}
\label{eq:2.4}
 & & \lim_{z\to -\infty}h(z)=0,\quad \lim_{z\to \infty}h(z)=1,\\
\label{eq:2.5}
 & & h'(z)=(4\pi)^{-\frac{1}{2}}e^{-\frac{|z|^2}{4}}>0,\quad h''(z)=-\frac{1}{2}zh'(z),\quad \mbox{for}\quad z\in{\mathbb R},\\
\label{eq:2.6}
 & &  h'(z)=-\left(\frac{1}{2}+o(1)\right)zh(z)\quad\mbox{as}\quad z\to -\infty.
\end{eqnarray}
\end{lemma}
{\bf Proof.} 
Let $z\in{\mathbb R}$. 
It follows from \eqref{eq:1.1} that 
\begin{align*}
\left(e^{t\Delta_{\mathbb{R}}}{\bf 1}_{[0,\infty)}\right)(z)
&=(4\pi t)^{-\frac{1}{2}}\int_{\mathbb{R}} e^{-\frac{|z-w|^2}{4t}}{\bf 1}_{[0,\infty)}(w)\,dw
=(4\pi t)^{-\frac12}\int_0^\infty  e^{-\frac{|z-w|^2}{4t}}\,dw\\
&=(4\pi)^{-\frac12}\int_0^\infty e^{-\frac{|t^{-\frac12}z -w|^2}{4}}\,dw
=h(t^{-\frac{1}{2}}z)\quad\mbox{for}\quad t>0.
\end{align*}
Then \eqref{eq:2.3} holds. Furthermore, 
$$
\lim_{t\to +0}h(t^{-\frac{1}{2}}z)=0\quad\mbox{if}\quad z<0,
\qquad
\lim_{t\to +0}h(t^{-\frac{1}{2}}z)=1\quad\mbox{if}\quad z>0,
$$
which yield \eqref{eq:2.4}. 
In addition, we have
\begin{equation*}
\begin{split}
t^{-\frac{1}{2}}h'(t^{-\frac{1}{2}}z)
 & =\partial_zh(t^{-\frac{1}{2}}z)
=(4\pi t)^{-\frac{1}{2}}\int_0^\infty
\partial_z\left(e^{-\frac{|z-w|^2}{4t}}\right)\,dw\\
 & =-(4\pi t)^{-\frac{1}{2}}\int_0^\infty
\partial_w\left(e^{-\frac{|z-w|^2}{4t}}\right)\,dw
=(4\pi t)^{-\frac{1}{2}}e^{-\frac{|z|^2}{4t}}
\end{split}
\end{equation*}
for $t>0$.
This implies that 
$$
h'(z)=(4\pi)^{-\frac{1}{2}}e^{-\frac{|z|^2}{4}}>0,
\qquad
h''(z)=-\frac{1}{2}z(4\pi)^{-\frac{1}{2}}e^{-\frac{|z|^2}{4}}=-\frac{1}{2}zh'(z).
$$
Thus \eqref{eq:2.5} holds.  

It follows from \eqref{eq:1.3} that 
\begin{equation}
\label{eq:2.7}
\begin{split}
h'(z)
&=-\frac{(4\pi)^{-\frac{1}{2}}}{2}\int_0^\infty (z-w)e^{-\frac{|z-w|^2}{4}}\,dw\\
&=-\frac{(4\pi)^{-\frac{1}{2}}}{2}z\int_0^\infty e^{-\frac{|z-w|^2}{4}}\,dw
+\frac{(4\pi)^{-\frac{1}{2}}}{2}
\biggr(\,\int_0^{|z|^{\frac{1}{2}}}+\int_{|z|^{\frac{1}{2}}}^\infty\,\biggr) w e^{-\frac{|z-w|^2}{4}}\,dw\\
&=-\frac{1}{2}z h(z)
+\frac{(4\pi)^{-\frac{1}{2}}}{2} 
\biggr(\,\int_0^{|z|^{\frac{1}{2}}}+\int_{|z|^{\frac{1}{2}}}^\infty\,\biggr) w e^{-\frac{|z-w|^2}{4}}\,dw.
\end{split}
\end{equation}
We have
\begin{equation}
\label{eq:2.8}
\begin{split}
\left| \frac{(4\pi)^{-\frac{1}{2}}}{2} \int_0^{|z|^{\frac12}} w e^{-\frac{|z-w|^2}{4}}\,dw \right|
&\leq \frac{(4\pi)^{-\frac{1}{2}}}{2}  |z|^{\frac12}\int_0^\infty e^{-\frac{|z-w|^2}{4}}\,dw
=\frac{1}{2}|z|^{\frac{1}{2}}h(z)
=o(|z|h(z))
\end{split}
\end{equation}
as $z\to-\infty$.
Furthermore, 
we find $C>1$ such that 
\begin{equation}
\label{eq:2.9}
\begin{split}
\left| \int_{|z|^{\frac{1}{2}}}^\infty w e^{-\frac{|z-w|^2}{4}}\,dw \right|
&= e^{-\frac{z^2}{4}} \int_{|z|^{\frac{1}{2}}}^\infty we^{\frac{zw}{2}-\frac{w^2}{4}}\,dw
\le Ce^{-\frac{z^2}{4}}\int_{|z|^{\frac{1}{2}}}^\infty e^{\frac{zw}{2}}\,dw\\
& =Ce^{-\frac{z^2}{4}} \left[ \frac{2}{z}  e^{\frac{zw}{2}}\right]^{w=\infty}_{w=|z|^{\frac12}}
=2C|z|^{-1}e^{-\frac{z^2+2|z|^{\frac{3}{2}}}{4}}\quad\mbox{if}\quad z<0.
\end{split}
\end{equation}
Since 
\begin{equation*}
\begin{split}
h(z) & \ge(4\pi)^{-\frac{1}{2}}\int_0^1 e^{-\frac{|z-w|^2}{4}}\,dw=(4\pi)^{-\frac{1}{2}}e^{-\frac{z^2}{4}}\int_0^1e^{\frac{zw}{2}-\frac{w^2}{4}}\,dw \\
 & \ge (4\pi)^{-\frac{1}{2}}e^{-\frac{1}{4}}e^{-\frac{z^2}{4}}\int_0^1e^{\frac{zw}{2}}\,dw
=(4\pi)^{-\frac{1}{2}}e^{-\frac{1}{4}}e^{-\frac{z^2}{4}}\left[\frac{2}{z}e^{\frac{zw}{2}}\right]_{w=0}^{w=1}\\
& =(4\pi)^{-\frac{1}{2}}e^{-\frac{1}{4}}e^{-\frac{z^2}{4}}\left(\frac{2}{z}e^{\frac{z}{2}}-\frac{2}{z}\right)
=2(4\pi)^{-\frac{1}{2}}e^{-\frac{1}{4}}|z|^{-1}e^{-\frac{z^2}{4}}(1+o(1))
\end{split}
\end{equation*}
as $z\to-\infty$, 
by \eqref{eq:2.9} we see that
$$
\left| \int_{|z|^{\frac{1}{2}}}^\infty w e^{-\frac{|z-w|^2}{4}}\,dw \right|
=o(|z|h(z))\quad\mbox{as}\quad z\to-\infty.
$$
This together with \eqref{eq:2.7} and \eqref{eq:2.8} implies \eqref{eq:2.6}. Thus Lemma~\ref{Theorem:2.9} follows. 
$\Box$
\begin{lemma}
\label{Theorem:2.10}
For any $f\in {\mathcal C}_\Omega[\Phi_0]\cap L^\infty(\Omega)$, 
there exists a sequence $\{f_a\}_{a>0}$ such that $f_a\in{\mathcal C}_\Omega[H_a]$ and 
$$
\lim_{a\to\infty}f_a(x)=f(x)
$$
uniformly on $\Omega$. 
\end{lemma}
{\bf Proof.}
The proof is a variation of the arguments in \cite{IST03}*{Section~4.2}. 
It suffices to treat large enough $a$. 
By \eqref{eq:2.5} we find $\varepsilon_a\in(0,1)$ such that 
\[
\varepsilon_a h'(-2\varepsilon_a^{-1})=a^{-1}\qquad 
\text{and}
\qquad 
\varepsilon_a\to 0\quad\mbox{as}\quad a\to\infty. 
\]
Set
\[
{h}_a(z):=ah(\varepsilon_a z-2\varepsilon_a^{-1})\quad\mbox{for}\quad z\in{\mathbb R},
\qquad
{h}_a(-\infty):=0.
\]
We observe from  Lemma~\ref{Theorem:2.9} that 
\begin{align*}
{h}_a'(z)=
a\varepsilon_a  h'(\varepsilon_a z-2\varepsilon_a^{-1})  
=a\varepsilon_a \cdot (4\pi)^{-\frac12}e^{-\frac{(\varepsilon_a z-2\varepsilon_a^{-1})^2}{4}} 
=a \varepsilon_a h'(-2\varepsilon_a^{-1})e^{z-\frac{1}{4}\varepsilon_a^2z^2}
=e^{z-\frac{1}{4}\varepsilon_a^2z^2}
\end{align*}
for  $z\in{\mathbb R}$, which together with $\lim_{z\to -\infty}h_a(z)=0$ yields 
\[
h_a(z)=\int_{-\infty}^z e^{w-\frac{1}{4}\varepsilon_a^2w^2}\,dw.
\]
Then we see that 
\begin{equation}
\label{eq:2.10}
\lim_{a\to\infty}h_a(z)=\int_{-\infty}^z e^w\,dw=e^z
\end{equation}
uniformly on $(-\infty,R)$ for any $R\in{\mathbb R}$.  

Let $f\in {\mathcal C}_\Omega[\Phi_0]\cap L^\infty(\Omega)$, and set
$$
f_a(x):=h_a(\Phi_0(f(x)))=h_a(\log f(x))\quad\mbox{for}\quad x\in\Omega. 
$$
Then, by \eqref{eq:2.10} we see that 
$\lim_{a\to\infty}f_a(x)=\lim_{a\to\infty}h_a(\log f(x))=f(x)$
uniformly on $\Omega$. 
Let $\widetilde{H}_a$ be the inverse function of $h_a$ in $[0,a)$. 
Then $\widetilde{H}_a$ is admissble on $[0,a)$
and it follows from $f\in {\mathcal C}_\Omega[\Phi_0]$ that $\widetilde{H}_a(f_a)=\log (f)$ is concave in $\Omega$.
These imply that
$$
f_a\in {\mathcal C}_\Omega[\widetilde{H}_a].
$$ 
On the other hand, since
$h_a(\varepsilon_a^{-1}H(a^{-1}z)+2\varepsilon_a^{-2})
=ah(H(a^{-1}z))=z$ for $z\in(0,a)$, 
we see that 
$$
\widetilde{H}_a(z)=\varepsilon_a^{-1}H(a^{-1}z)+2\varepsilon_a^{-2}
=\varepsilon_a^{-1}H_a(z)+2\varepsilon_a^{-2}\quad\mbox{for}\quad z\in(0,a).
$$
Then Lemma~\ref{Theorem:2.5} implies that 
${\mathcal C}_\Omega[\widetilde{H}_a]={\mathcal C}_\Omega[H_a]$. 
Thus 
$f_a\in {\mathcal C}_\Omega[\widetilde{H}_a]={\mathcal C}_\Omega[H_a]$. 
The proof is complete. 
$\Box$
\section{Disruption of $F$-concavity}\label{section:3}
We consider the Cauchy--Dirichlet problem 
\begin{equation}
\tag{P}
\left\{
\begin{array}{ll}
\partial_t u=Lu & \mbox{in}\quad\Omega\times(0,\infty),\vspace{3pt}\\
u=0 & \mbox{in}\quad\partial\Omega\times(0,\infty)\text{ if }\partial \Omega \neq \emptyset,\vspace{3pt}\\
u(\cdot,0)=\phi\ge 0 & \mbox{in}\quad\Omega,
\end{array}
\right.
\end{equation}
where $\phi\in L^\infty(\Omega)$.
Here $L$ is an elliptic operator of the form 
$$
L:=\sum_{i,j=1}^na^{ij}(x,t)\partial_{x_i}\partial_{x_j}+\sum_{i=1}^n b^i(x,t)\partial_{x_i}-c(x,t),
$$
and the coefficients satisfy the following conditions: 
\begin{itemize}
  \item[(L1)] 
  there exists $\sigma\in(0,1)$ such that 
  $$
  a^{ij},\,\,b^i\in C^{0,\sigma;0,\sigma/2}(\Omega\times[0,\infty)),\qquad
  c\in BC^{0,\sigma;0,\sigma/2}(\Omega\times[0,\infty)),
  $$ 
  where $i$, $j=1,\dots,n$;
  \item[(L2)] 
  $A(x,t):=(a^{ij}(x,t))\in\mathrm{Sym}\,(n)$ for $(x,t)\in\Omega\times[0,\infty)$ 
  and there exists $\Lambda>0$ such that 
  $$
  \Lambda^{-1}|\xi|^2\le
  \langle A(x,t)\xi,\xi\rangle\le\Lambda|\xi|^2
  \quad\mbox{for all $\xi\in{\mathbb R}^n$ and $(x,t)\in\Omega\times[0,\infty)$}.
  $$
\end{itemize}
Then, for any nonnegative initial datum $\phi\in L^\infty(\Omega)$, 
problem~(P) possesses a (unique) minimal nonnegative solution $e^{tL_\Omega}\phi$ 
such that 
\begin{align}
\nonumber
 & e^{tL_\Omega}\phi\in L^\infty((0,\infty):L^2_{{\rm loc}}(\Omega))\cap L^2((0,\infty): H^1_{{\rm loc}}(\Omega)),\\
\label{eq:3.1}
 & \lim_{t\to +0}\|e^{tL_\Omega}\phi-\phi\|_{L^2(\Omega\cap B(0,R))}=0\quad\mbox{for}\quad R>0.
\end{align}
(See e.g. \cite{LSU}*{Chapter III} and \cite{IM}*{Lemma~5.3}.)  
The solution~$e^{tL_\Omega}\phi$ is represented by the minimal Dirichlet heat kernel $G_{L_\Omega}$ associated with the operator $L$ in $\Omega$ 
as follows:
$$
(e^{tL_\Omega}\phi)(x)=\int_\Omega G_{L_\Omega}(x,y,t)\phi(y)\,dy,
\quad (x,t)\in\Omega\times(0,\infty).
$$
Then, under conditions~(L1) and (L2), 
parabolic regularity theorems (see e.g. \cite{LSU}*{Chapter~IV, Theorem~16.3}) imply
that $G_{L_{\Omega}}\in C^{2,\sigma;1,\sigma/2}(\Omega\times\Omega\times(0,\infty))$. 
Furthermore, we observe from the maximum principle and the comparison principle that
\begin{equation}
\label{eq:3.2}
G_{L_\Omega}(x,y,t)>0,\quad 
\int_\Omega G_{L_\Omega}(x,y,t)\,dy\le e^{\int_0^t \|c(s)\|_{L^\infty(\Omega)}\,ds},
\end{equation}
for $x$, $y\in\Omega$ and $t>0$.  
\begin{definition}
\label{Theorem:3.1}
Let $F$ be admissible on $I=[0,a)$ with $a\in(0,\infty]$ and  $\Omega$ a convex domain in $\mathbb{R}^n$. 
Consider problem~{\rm (P)} under conditions~{\rm (L1)} and {\rm (L2)}. 
We say that \emph{$F$-concavity is preserved by the Dirichlet parabolic flow associated with $L$ in $\Omega$} if 
$$
\mbox{$e^{tL_\Omega}\phi \in \mathcal{C}_{\Omega}[F]$ for $t>0$ 
holds for all $\phi\in\mathcal{C}_{\Omega}[F]\cap L^\infty(\Omega)$}.
$$
\end{definition}
In this section we study the disruption of $F$-concavity by the Dirichlet parabolic flow, 
in particular, by DHF. 
\subsection{Disruption of  {quasi-concavity} by DHF in ${\mathbb R}^n$}\label{subsection:3.1}
We first prove that, 
when starting with a non-log-concave initial datum $\phi$ in ${\mathbb R}^n$ with $n\ge 2$, 
then the solution~$e^{t\Delta_{{\mathbb R}^n}}\phi$ may not be quasi-concave for all small enough $t>0$,
hence losing any reminiscence of concavity. 
Proposition~\ref{Theorem:3.2} is one of the main ingredients of this paper. 
\begin{proposition}
\label{Theorem:3.2}
Let $F$ be admissible on $I=[0,a)$ with $a\in(0,\infty]$. 
Let $n\ge 2$ and assume that $F$-concavity
 is not stronger than log-concavity  in ${\mathcal A}_{{\mathbb R}^n}(I)$, 
that is,
\begin{equation}
\label{eq:3.3}
{\mathcal C}_{{\mathbb R}^n}[F]\setminus {\mathcal C}_{{\mathbb R}^n}[\Phi_0]\not=\emptyset.
\end{equation}
Then there exists $\phi\in \mathcal{C}_{{\mathbb R}^n}[F]\cap L^\infty({\mathbb R}^n)$ such that 
$e^{t\Delta_{{\mathbb R}^n}}\phi$ is not quasi-concave in ${\mathbb R}^n$ for all small enough $t>0$. 
\end{proposition}
{\bf Proof.}
The proof heavily depends on the following nice property of the heat flow:
\begin{itemize}
  \item for any $\phi_1$, $\phi_2\in L^\infty({\mathbb R})$, 
  $$
  (e^{t\Delta_{{\mathbb R}^2}}\phi)(z)=(e^{t\Delta_{\mathbb R}}\phi_1)(z_1)(e^{t\Delta_{\mathbb R}}\phi_2)(z_2)
  $$
  for $z=(z_1,z_2)\in{\mathbb R}^2$ and $t>0$, where $\phi(z)=\phi_1(z_1)\phi_2(z_2)$.
\end{itemize}
Setting $\phi_1={\bf 1}_{[0,\infty)}$ and letting $\phi_2$ a suitable modification of the inverse of $F$, 
we see that $\phi=\phi_1\phi_2$ is $F$-concave in ${\mathbb R}^2$. 
Then we prove the non-convexity of a superlevel set of $e^{t\Delta_{{\mathbb R}^2}}\phi$ for all small enough $t>0$. 
The proof is divided into three steps. 
\newline
{\bf Step 1}: 
Assume \eqref{eq:3.3}. It follows from Lemma~\ref{Theorem:2.7} that 
$\lim_{r\to+0} F(r)=-\infty$. 
We write $f:=f_F$ and $J:=J_F$ for simplicity. 
The admissibility of $F$ together with $\lim_{r\to+0} F(r)=-\infty$ implies that $J=(-\infty,F(a))$.
Furthermore, $f$ is positive, continuous, and strictly increasing in $J$. 
By Lemma~\ref{Theorem:2.4} we see that $\log f$ is not concave in $J$, 
that is, 
\begin{equation}
\label{eq:3.4}
f((1-\lambda)\zeta+\lambda \omega)<f(\zeta)^{1-\lambda}f(\omega)^\lambda
\end{equation}
for some $\zeta, \omega \in J$ with $\zeta<\omega$ and $\lambda\in(0,1)$. 
Let $c\in J$ be such that $c>\omega$. Set
$$
\varphi(z):=
\left\{
\begin{array}{ll}
f(z+c) & \mbox{for}\quad z\in(-\infty,0],\vspace{3pt}\\
f(-z+c) & \mbox{for}\quad z\in(0,\infty).
\end{array}
\right.
$$
Then 
\begin{itemize}
  \item $\varphi$ is strictly increasing in $(-\infty,0]$ and $\varphi(z)\to 0$ as $z\to-\infty$;
  \item $\varphi$ is positive and continuous in ${\mathbb R}$ such that $\varphi(z)=\varphi(-z)$ and $\varphi(z)\le\varphi(0)=f(c)<a$ for $z\in{\mathbb R}$;
  \item $\varphi$ is $F$-concave in ${\mathbb R}$ and it is not log-concave in $(-\infty,0]$.
\end{itemize} 

Let 
\begin{equation}
\label{eq:3.5}
v(z,t):=(e^{t\Delta_{{\mathbb R}}}\varphi)(z)\quad\mbox{for}\quad (z,t)\in{\mathbb R}\times(0,\infty).
\end{equation}
Then $v\in C^\infty({\mathbb R}\times(0,\infty))$, $v(z,t)=v(-z,t)$ for $(z,t)\in {\mathbb R}\times(0,\infty)$, and
\begin{equation}
\label{eq:3.6}
\left\{
\begin{array}{ll}
\partial_z v>0\quad & \mbox{in}\quad (-\infty,0)\times(0,\infty),\vspace{3pt}\\
\partial_z v<0\quad & \mbox{in}\quad (0,\infty)\times(0,\infty).
\end{array}
\right. \end{equation}
(See also \cite{An}.) 
Furthermore, it follows from the continuity of $\varphi$ that 
\begin{equation}
\label{eq:3.7}
\lim_{t\to +0}\sup_{z\in K}|v(z,t)-\varphi(z)|=0
\end{equation}
for any bounded interval $K\subset{\mathbb R}$. 
This together with \eqref{eq:3.4} yields 
$$
v((1-\lambda)(\zeta-c)+\lambda (\omega-c),t)<v(\zeta-c,t)^{1-\lambda}v(\omega-c,t)^\lambda
$$
for all small enough $t>0$. 
This means that, 
for any small enough $t>0$, 
$v(\cdot,t)$ is not log-concave in $[\zeta-c,\omega-c]$, 
that is, there exists $z_t\in[\zeta-c,\omega-c]$ such that
\begin{equation}
\label{eq:3.8}
0<(\partial_z^2\log v)(z_t,t)=\frac{v(z_t,t)(\partial_z^2 v)(z_t,t)-((\partial_z v)(z_t,t))^2}
{((\partial_z v)(z_t,t))^2}.
\end{equation}
Since $\varphi>0$ on $[\zeta-c,\omega-c]$, 
by \eqref{eq:3.7} we find $C\in[1,\infty)$ such that 
\begin{equation}
\label{eq:3.9}
C^{-1}\le v(z_t,t)\le C
\end{equation}
for small enough $t>0$. 

Let $t>0$ be small enough such that \eqref{eq:3.8} holds. 
By \eqref{eq:3.6} we find the inverse function $\Gamma_t$ of 
the function $(-\infty,0)\ni z\mapsto v(z,t)$. 
Then $\Gamma_t$ is smooth and $\Gamma_t'>0$ in $(0,v(0,t))$ and 
$$
\Gamma_t(v(z,t))=z\quad\mbox{for}\quad z\in(-\infty,0),
$$
which implies that
$$
\Gamma_t''(v(z,t))((\partial_z v)(z,t))^2+\Gamma_t'(v(z,t))(\partial_z^2 v)(z,t)=0\quad\mbox{for}\quad z\in(-\infty,0).
$$
This together with \eqref{eq:3.8} yields
\begin{equation}
\label{eq:3.10}
\begin{split}
v(z_t,t)\Gamma_t''(v(z_t,t)) & =-v(z_t,t)\frac{\Gamma_t'(v(z_t,t))(\partial_z^2v)(z_t,t)}{((\partial_z v)(z_t,t))^2}\\
 & =-\Gamma_t'(v(z_t,t))\frac{v(z_t,t)(\partial_z^2 v)(z_t,t)-((\partial_z v)(z_t,t))^2}
{((\partial_z v)(z_t,t))^2}-\Gamma_t'(v(z_t,t))\\
 & <-\Gamma_t'(v(z_t,t)).
\end{split}
\end{equation}
{\bf Step 2}: 
Set 
$$
\phi(w,z):={\bf 1}_{[0,\infty)}(w)\varphi(z)
\quad\mbox{for}\quad
(w,z)\in\mathbb{R}^2.
$$
Then $\phi$ is $F$-concave in $\mathbb{R}^2$. 
It follows from \eqref{eq:2.3} and \eqref{eq:3.5} that 
\[
u(w,z,t):=(e^{t\Delta_{\mathbb{R}^2}}\phi)(w,z)=h(t^{-\frac{1}{2}}w)v(z,t)
\quad\mbox{for}\quad(w,z,t)\in\mathbb{R}^2\times(0,\infty).
\]
Here $h$ is as in \eqref{eq:1.3}.  

Let $t>0$ be small enough. 
Since $h'>0$ in ${\mathbb R}$ and 
$\lim_{w\to -\infty}h(w)=0$, 
by \eqref{eq:3.9}, 
for any small enough $\varepsilon>0$, 
we find a unique ${w_\varepsilon}\in(-\infty,0)$ such that
\begin{equation*}
u(w_\varepsilon,z_t,t)=
h(t^{-\frac{1}{2}}{w_\varepsilon})v(z_t,t)=\varepsilon.
\end{equation*}
Then 
\begin{equation}
\label{eq:3.11}
 v(z_t,t)=\frac{\varepsilon}{h(t^{-\frac{1}{2}}{w_\varepsilon})}\quad\mbox{for small enough $\varepsilon>0$},
 \qquad
t^{-\frac{1}{2}}{w_\varepsilon}\to -\infty\quad\mbox{as}\quad \varepsilon\to +0.
\end{equation}
By \eqref{eq:3.6}, applying the implicit function theorem, 
we find a smooth function $g_\varepsilon$ in a neighborhood~${\mathcal N}_\varepsilon$ of $w_\varepsilon$ 
and $\delta_\varepsilon>0$ such that 
\begin{equation}
\begin{split}\label{eq:3.12}
 & \varepsilon=u(w,g_\varepsilon(w),t)
=h(t^{-\frac{1}{2}}w)v(g_\varepsilon(w),t)\quad\mbox{for}\quad w\in {\mathcal N}_\varepsilon,\\
 & g_\varepsilon(w_\varepsilon)=z_t=\Gamma_t(v(z_t,t)),\\
 & u(w,z,t)<\varepsilon\quad\mbox{if $w\in {\mathcal N}_\varepsilon$ 
 and $g_\varepsilon(w)-\delta_\varepsilon<z<g_\varepsilon(w)$},\\
 & u(w,z,t)>\varepsilon\quad\mbox{if $w\in {\mathcal N}_\varepsilon$ 
 and $g_\varepsilon(w)<z<g_\varepsilon(w)+\delta_\varepsilon$}.
\end{split}
\end{equation}
Furthermore, it follows from \eqref{eq:3.12} that 
\[
g_\varepsilon(w)=
\Gamma_t\left(\frac{\varepsilon}{h(t^{-\frac{1}{2}}w)}\right)\quad\mbox{for}\quad w\in {\mathcal N}_\varepsilon.
\]
\noindent
{\bf Step 3}:
Assume that $u(\cdot,t)$ is quasi-concave in ${\mathbb R}^n$ for some small $t>0$.
Then it follows from~\eqref{eq:3.12} that $g_\varepsilon$ is convex in ${\mathcal N}_\varepsilon$, so that
\begin{equation}
\label{eq:3.13}
g''_\varepsilon(w_\varepsilon)\ge 0.
\end{equation}
A direct computation provides
\begin{align*}
g_\varepsilon'(w) & =-\varepsilon t^{-\frac{1}{2}}\Gamma_t'\left(\frac{\varepsilon}{h(t^{-\frac{1}{2}}w)}\right)
\frac{h'(t^{-\frac{1}{2}}w)}{h(t^{-\frac{1}{2}}w)^2},\\
g_\varepsilon''(w) &
=\varepsilon^2 t^{-1}\Gamma_t''\left(\frac{\varepsilon}{h(t^{-\frac{1}{2}}w)}\right)
\frac{h'(t^{-\frac{1}{2}}w)^2}{h(t^{-\frac{1}{2}}w)^4}
+2\varepsilon t^{-1}\Gamma'_t\left(\frac{\varepsilon}{h(t^{-\frac{1}{2}}w)}\right)\frac{h'(t^{-\frac{1}{2}}w)^2}{h(t^{-\frac{1}{2}}w)^3}\\
 & \quad-\varepsilon t^{-1}
\Gamma_t'\left(\frac{\varepsilon}{h(t^{-\frac{1}{2}}w)}\right)\frac{h''(t^{-\frac{1}{2}}w)}{h(t^{-\frac{1}{2}}w)^2}.
\end{align*}
This together with \eqref{eq:2.5} and \eqref{eq:3.11} leads to
\begin{align}
\label{eq:3.14}
\begin{split}
& g_\varepsilon''(w_\varepsilon)\\
&=\varepsilon t^{-1}\frac{\Gamma_t'(v(z_t,t))}{h(t^{-\frac{1}{2}}w_\varepsilon)^3}
\left[
\left(v(z_t,t) \frac{\Gamma_t''(v(z_t,t))}{\Gamma_t'(v(z_t,t))}+2\right)h'(t^{-\frac{1}{2}}w_\varepsilon)^2
-h(t^{-\frac{1}{2}}w_\varepsilon)h''(t^{-\frac{1}{2}}w_\varepsilon)
\right]\\
&=\varepsilon t^{-1}\frac{\Gamma_t'(v(z_t,t))}{h(t^{-\frac{1}{2}}w_\varepsilon)^3}\\
& \quad\times\left[
\left(v(z_t,t) \frac{\Gamma_t''(v(z_t,t))}{\Gamma_t'(v(z_t,t))}+2\right)h'(t^{-\frac{1}{2}}w_\varepsilon)^2
+\frac{1}{2}t^{-\frac{1}{2}}w_\varepsilon h(t^{-\frac{1}{2}}w_\varepsilon) h'(t^{-{\frac{1}{2}}}w_\varepsilon)
\right].
\end{split}
\end{align}
We observe from \eqref{eq:2.6} and \eqref{eq:3.11} that
\[
h'(t^{-\frac{1}{2}}w_\varepsilon)=-\left(\frac{1}{2}+o(1)\right)t^{-\frac12}w_\varepsilon h(t^{-{\frac{1}{2}}}w_\varepsilon)
\quad\mbox{as}\quad \varepsilon\to +0,
\]
hence
\[
\frac{1}{2}t^{-\frac{1}{2}}w_\varepsilon h(t^{-{\frac{1}{2}}}w_\varepsilon)h'(t^{-\frac{1}{2}}w_\varepsilon)
=-(1+o(1))^{-1}h'(t^{-\frac{1}{2}}w_\varepsilon)^2
=-(1+o(1))h'(t^{-\frac{1}{2}}w_\varepsilon)^2
\]
as $\varepsilon\to +0$.
This together with \eqref{eq:3.10} and \eqref{eq:3.14} implies that 
\begin{align*}
g_\varepsilon''(w_\varepsilon)
 & =\varepsilon t^{-1}\frac{\Gamma_t'(v(z_t,t))}{h(t^{-\frac{1}{2}}w_\varepsilon)^3}
h'(t^{-\frac{1}{2}}w_\varepsilon)^2
\left(v(z_t,t) \frac{\Gamma_t''(v(z_t,t))}{\Gamma_t'(v(z_t,t))}+1+o(1)\right)\\
 & =\varepsilon t^{-1}\frac{h'(t^{-\frac{1}{2}}w_\varepsilon)^2}{h(t^{-\frac{1}{2}}w_\varepsilon)^3}
\left(v(z_t,t)\Gamma_t''(v(z_t,t))+(1+o(1))\Gamma_t'(v(z_t,t))\right)<0
\end{align*}
for all small enough $\varepsilon>0$. 
This contradicts \eqref{eq:3.13}. 
We deduce that $u(\cdot,t)$ is not quasi-concave in $\mathbb{R}^2$ for all small enough $t>0$. 
Then Proposition~\ref{Theorem:3.2} follows in the case $n=2$. 
If $n\ge 3$, then 
we set 
\begin{equation}
\label{eq:3.15}
U(w,z,x',t):=u(w,z,t)\quad\mbox{for}\quad
(w,z,x',t)\in\mathbb{R}^2\times\mathbb{R}^{n-2}\times(0,\infty).
\end{equation}
Then $U(\cdot,t)$ is not quasi-concave in $\mathbb{R}^n$ for all small enough $t>0$ 
and $U(\cdot,0)$ is $F$-concave in~${\mathbb R}^n$. 
Thus Proposition~\ref{Theorem:3.2} follows in the case $n\ge 3$, 
and the proof of Proposition~\ref{Theorem:3.2} is complete.~$\Box$
%
\subsection{Disruption of $F$-concavity by the Dirichlet parabolic flow}\label{subsection:3.2}
We show that 
the disruption of $F$-concavity (resp.\,\,quasi-concavity) by DHF in ${\mathbb R}^n$ implies the disruption 
of $F$-concavity (resp.\,\,quasi-concavity) by the Dirichlet parabolic flow in $\Omega$. 
\begin{proposition}
\label{Theorem:3.3}
Let $F$ be admissible on $I=[0,a)$ with $a\in(0,\infty]$ and
$\Omega$ a convex domain in~${\mathbb R}^n$. 
Assume that $F$-concavity is not preserved by DHF in ${\mathbb R}^n$. 
Then, under conditions~{\rm (L1)} and {\rm (L2)}, 
there exists~$\phi\in{\mathcal C}_\Omega[F]\cap BC_0(\overline{\Omega})$ 
 such that $e^{tL_\Omega}\phi$ is not $F$-concave in~$\Omega$ for some $t>0$. 
\end{proposition}
Proposition~\ref{Theorem:3.3} is proved by the following lemma and the similar transformation of DHF. 
\begin{lemma}
\label{Theorem:3.4}
Let $F$ be admissible on $I=[0,a)$ with $a\in(0,\infty]$ and $\phi\in {\mathcal C}_\Omega[F]\cap L^\infty(\Omega)$. 
Assume conditions~{\rm (L1)} and {\rm (L2)}. 
Then there exists a sequence of $\{\phi_j\}\subset{\mathcal C}_\Omega[F]\cap BC_0(\overline{\Omega})$ 
such that 
$$
\lim_{j\to\infty}\left(e^{tL_\Omega}\phi_j\right)(x)=\left(e^{tL_\Omega}\phi\right)(x)
\quad\mbox{for}\quad (x,t)\in\Omega\times(0,\infty).
$$
\end{lemma}
{\bf Proof.} 
For any $\varepsilon>0$ and $\delta\in(0,a)$, 
we set 
$$
\psi_{\varepsilon,\delta}(x):=\left(e^{\varepsilon\Delta_\Omega}e^{F(\min\{\phi,a-\delta\})}\right)(x)\quad\mbox{for}\quad x\in\Omega.
$$
By \eqref{eq:3.2} we have
\begin{equation}
\label{eq:3.16}
\|\psi_{\varepsilon,\delta}\|_{L^\infty(\Omega)}
\le \left\|e^{F(\min\{\phi,a-\delta\})}\right\|_{L^\infty(\Omega)}
\le e^{F(\min\{\|\phi\|_{L^\infty(\Omega)},a-\delta\})}.
\end{equation}
By parabolic regularity theorems (see e.g. \cite{LSU}*{Chapter~III, Theorem~10.1}) 
we see that $\psi_{\varepsilon,\delta}\in BC_0(\overline{\Omega})$. 
Since $e^{F(\min\{\phi,a-\delta\})}$ is log-concave in $\Omega$, 
thanks to the preservation of log-concavity by DHF in $\Omega$, we observe that 
$\psi_{\varepsilon,\delta}$ is log-concave in $\Omega$, which implies that the set 
$$
E_{\varepsilon,\delta}:=\big\{x\in\Omega\,\,|\,\,\psi_{\varepsilon,\delta}(x)>\lim_{r\to +0}e^{F(r)}\big\}
$$ 
is convex. 
Since the function $\psi_{\varepsilon,\delta}{\bf 1}_{E_{\varepsilon,\delta}}$ is log-concave in $\Omega$ 
and 
$$
\log\left(\psi_{\varepsilon,\delta}(x)\right)\in F((0,\min\{\|\phi\|_{L^\infty(\Omega)},a-\delta\}])\quad\mbox{for}\quad x\in E_{\varepsilon,\delta},
$$ 
we define an $F$-concave function $\phi_{\varepsilon,\delta}\in BC_0(\overline{\Omega})$ by 
$$
\phi_{\varepsilon,\delta}(x):=
\left\{
\begin{array}{ll}
f_F\left(\log\left(\psi_{\varepsilon,\delta}(x)\right)\right) & \mbox{if}\quad x\in E_{\varepsilon,\delta},\vspace{3pt}\\
0 & \mbox{otherwise}. 
\end{array}
\right.
$$
Then, by \eqref{eq:3.16} we have
\begin{equation}
\label{eq:3.17}
\|\phi_{\varepsilon,\delta}\|_{L^\infty(\Omega)}\le \min\{\|\phi\|_{L^\infty(\Omega)},a-\delta\}\le \|\phi\|_{L^\infty(\Omega)}.
\end{equation}
Furthermore, by \eqref{eq:3.1}, 
for any $\delta\in(0,a)$,  
we find a sequence $\{\varepsilon_j\}$ with $\lim_{j\to\infty}\varepsilon_j=0$ such that 
\begin{equation*}
\lim_{j\to\infty} \psi_{\varepsilon_j,\delta}(x)=e^{F(\min\{\phi(x),a-\delta\})}
\end{equation*}
for almost all $x\in\Omega$. 
This implies that 
\begin{equation*}
\lim_{j\to\infty} \phi_{\varepsilon_j,\delta}(x)=\min\{\phi(x),a-\delta\}
\end{equation*}
for almost all $x\in\Omega$. 
By \eqref{eq:3.2} and \eqref{eq:3.17} 
we apply Lebesgue's dominated convergence theorem to obtain  
\begin{align*}
\lim_{\delta\to +0}\lim_{j\to\infty}\left(e^{tL_\Omega}\phi_{\varepsilon_j,\delta}\right)(x)
 & =\lim_{\delta\to +0}\lim_{j\to\infty}\int_\Omega G_{L_\Omega}(x,y,t)\phi_{\varepsilon_j,\delta}(y)\,dy\\
 & =\lim_{\delta\to +0}\int_\Omega G_{L_\Omega}(x,y,t)\min\{\phi(y),a-\delta\}\,dy
=\left(e^{tL_\Omega}\phi\right)(x)
\end{align*}
for $(x,t)\in\Omega\times(0,\infty)$. 
Then we obtain the desired conclusion, and the proof is complete. 
$\Box$
\vspace{5pt}
\newline
\noindent{\bf Proof of Proposition~\ref{Theorem:3.3}.}
Assume that $F$-concavity is not preserved by DHF in ${\mathbb R}^n$. 
Then we find $\phi\in {\mathcal C}_{{\mathbb R^n}}[F]\cap L^\infty({\mathbb R}^n)$ such that 
$e^{\tau\Delta_{{\mathbb R}^n}}\phi$ is not $F$-concave in ${\mathbb R}^n$ for some $\tau>0$,
that is, 
there exist $\xi$, $\eta\in{\mathbb R}^n$ and $\lambda\in(0,1)$ such that
\begin{equation}
\label{eq:3.18}
F\left(\left(e^{\tau\Delta_{{\mathbb R}^n}}\phi\right)((1-\lambda)\xi+\lambda \eta)\right)
<(1-\lambda)F\left(\left(e^{\tau\Delta_{{\mathbb R}^n}}\phi\right)(\xi)\right)
+\lambda F\left(\left(e^{\tau\Delta_{{\mathbb R}^n}}\phi\right)(\eta)\right).
\end{equation}
Thanks to \eqref{eq:1.2},
we can assume, without loss of generality, that 
\begin{equation}
\label{eq:3.19}
0\in\Omega,\qquad a^{ij}(0,0)=\delta^{ij}\quad\mbox{for}\quad i,j=1,\dots,n,
\end{equation}
where $\delta^{ij}=1$ if $i=j$ and $\delta^{ij}=0$ if $i\not=j$. 
For any $\ell=1,2,\dots$, set
\begin{equation*}
\begin{split}
 & \phi_\ell(x):=\phi(\ell x),\quad 
u_\ell(x,t):=(e^{tL_\Omega}\phi_\ell)(x),\quad\mbox{for $x\in\Omega$ and $t>0$},\\
 & U_\ell(x,t):=u_\ell(\ell^{-1}x,\ell^{-2}t)\quad\mbox{for $x\in\Omega_\ell:=\ell\Omega$ and $t>0$}. 
\end{split}
\end{equation*}
It follows from \eqref{eq:3.2} that
\begin{equation}
\label{eq:3.20}
\|U_\ell\|_{L^\infty(\Omega_\ell\times(0,\ell^2T))}
=\|u_\ell\|_{L^\infty(\Omega\times(0,T))}\le e^{T\|c\|_{L^\infty(\Omega\times(0,T))}}\|\phi_\ell\|_{L^\infty(\Omega)}<\infty
\end{equation}
for any $T>0$.
Furthermore, $U_\ell$ satisfies 
\begin{equation*}
\left\{
\begin{array}{ll}
\partial_t U_\ell=L_\ell U_\ell & \mbox{in}\quad\Omega_\ell\times(0,\infty),\vspace{3pt}\\
U_\ell=0  & \mbox{in}\quad\partial\Omega_\ell\times(0,\infty)\text{ if }\partial \Omega_\ell \neq \emptyset,\vspace{3pt}\\
U_\ell(x,0)=\phi(x) & \mbox{in}\quad\Omega_\ell,
\end{array}
\right.
\end{equation*}
where
\begin{equation*}
\begin{split}
 & L_\ell:=\sum_{i,j=1}^n a^{ij}_\ell(x,t)\partial_{x_i} \partial_{x_j}
+\sum_{i=1}^n b^i_\ell(x,t)\partial_{x_i}-c_\ell(x,t),\\
 & a^{ij}_\ell(x,t):=a^{ij}(\ell^{-1}x,\ell^{-2}t),
\quad
b^i_\ell(x,t):=\ell^{-1}b^i(\ell^{-1}x,\ell^{-2}t),
\quad
c_\ell(x,t):=\ell^{-2}c(\ell^{-1}x,\ell^{-2}t),
\end{split}
\end{equation*}
for $(x,t)\in\Omega_\ell\times(0,\infty)$ and $i,j=1,\dots,n$. 
Furthermore, by \eqref{eq:3.19} and condition~(L1) 
we see that 
$\Omega_\ell\to{\mathbb R}^n$ as $\ell\to\infty$ and 
\begin{itemize}
  \item 
  the coefficients $a^{ij}_\ell$, $b^i_\ell$, $c_\ell$ are bounded in $C^{0,\sigma;\,0,\sigma/2}(K)$,
  \item 
  $a^{ij}_\ell(x,t)\to \delta^{ij}$, $b^i_\ell(x,t)\to 0$, and $c_\ell(x,t)\to 0$,  
  as~$\ell\to\infty$ uniformly on~$K$,
\end{itemize}
for all compact sets $K\subset{\mathbb R}^n\times[0,\infty)$, where $ i,j=1,\dots,n$. 
Applying parabolic regularity theorems to $\{U_\ell\}$ (see e.g. \cite{LSU}*{Chapters~III, Theorem~10.1 and Chapter~IV, Theorem~10.1}) with \eqref{eq:3.20}, 
we have:
\begin{itemize}
  \item 
  $\{U_\ell\}$ are uniformly bounded and equicontinuous for 
  all compact sets in ${\mathbb R}^n\times[0,\infty)$;
  \item 
  $\displaystyle{\sup_{\ell}}\,\|U_\ell\|_{C^{2,\sigma;\,1,\sigma}(K)}<\infty$ 
  for all compact sets $K\subset{\mathbb R}^n\times(0,\infty)$.
\end{itemize}
Applying the Arzel\`a--Ascoli theorem and the diagonal argument, 
we find a subsequence~$\{U_{\ell_j}\}$ of $\{U_\ell\}$ such that
$$
U_{\ell_j}(x,t)\to (e^{t\Delta_{{\mathbb R}^n}}\phi)(x)\quad\mbox{as}\quad j\to\infty
$$
uniformly for all compact sets in $\mathbb{R}^n\times(0,\infty)$. 
This together with \eqref{eq:3.18} implies that 
\begin{equation*}
F(U_{\ell_j}((1-\lambda)\xi+\lambda \eta,\tau))<(1-\lambda)F(U_{\ell_j}(\xi,\tau))+\lambda F(U_{\ell_j}(\eta,\tau)),
\end{equation*}
that is, 
\begin{equation}
\label{eq:3.21}
\begin{split}
 & F(u_{\ell_j}((1-\lambda)\ell_j^{-1}\xi+\lambda \ell_j^{-1}\eta, \ell_j^{-2}\tau))\\
 & <(1-\lambda)F(u_{\ell_j}(\ell_j^{-1}\xi,\ell_j^{-2}\tau))+\lambda F(u_{\ell_j}(\ell_j^{-1}\eta,\ell_j^{-2}\tau))
\end{split}
\end{equation}
for all large enough $j$. 
Since $\ell^{-1}\xi, \ell^{-1}\eta\in\Omega$ for large enough $\ell$, 
this means that $u_{\ell_j}(\cdot,\ell_j^{-2}\tau)$ is not $F$-concave in $\Omega$ for large enough $j$. 
Taking into account that $\phi_{\ell_j}$ is $F$-concave in $\Omega$, 
thanks to Lemma~\ref{Theorem:3.4},  
we approximate $\phi_{\ell_j}$ by 
$F$-concave functions belonging to $BC_0(\overline{\Omega})$ 
to obtain the desired conclusion. 
Thus Proposition~\ref{Theorem:3.3} follows. 
$\Box$
\begin{proposition}
\label{Theorem:3.5}
Let $F$ be admissible on $I=[0,a)$ with $a\in(0,\infty]$ and
$\Omega$ a convex domain in~${\mathbb R}^n$. 
Assume that there exists $\phi\in{\mathcal C}_{{\mathbb R}^n}[F]\cap L^\infty({\mathbb R}^n)$
such that $e^{\tau\Delta_{{\mathbb R}^n}}\phi$ is not quasi-concave in~${\mathbb R}^n$ for some $\tau>0$. 
Then, under conditions~{\rm (L1)} and {\rm (L2)}, 
there exists $\psi\in{\mathcal C}_\Omega[F]\cap BC_0(\overline{\Omega})$ 
such that $e^{tL_\Omega}\psi$ is not quasi-concave in $\Omega$ for some $t>0$. 
\end{proposition}
{\bf Proof.}
The proof is similar to that of Proposition~\ref{Theorem:3.3}. 
Indeed, 
under the assumption of Proposition~\ref{Theorem:3.5}, 
we find $\xi$, $\eta\in{\mathbb R}^n$, $\lambda\in(0,1)$, and $\tau>0$ such that 
$$
\left(e^{\tau\Delta_{{\mathbb R}^n}}\phi\right)((1-\lambda)\xi+\lambda \eta)
<\min\left\{\left(e^{\tau\Delta_{{\mathbb R}^n}}\phi\right)(\xi),
\left(e^{\tau\Delta_{{\mathbb R}^n}}\phi\right)(\eta)\right\},
$$
instead of \eqref{eq:3.18}. 
The same argument as in the proof of Proposition~\ref{Theorem:3.3} implies that 
$$
u_{\ell_j}((1-\lambda)\ell_j^{-1}\xi+\lambda \ell_j^{-1}\eta,\ell_j^{-2}\tau)
<\min\{u_{\ell_j}(\ell_j^{-1}\xi,\ell_j^{-2}\tau),u_{\ell_j}(\ell_j^{-1}\eta,\ell_j^{-2}\tau)\},
$$
for some large enough $j$, instead of \eqref{eq:3.21}. 
Thus quasi-concavity is not preserved by the Dirichlet parabolic flow associated with $L$ in $\Omega$. 
Consequently, by Lemma~\ref{Theorem:3.4} we obtain the desired conclusion. 
The proof is complete. 
$\Box$
\section{Preservation of $F$-concavity}\label{section:4}
The preservation of log-concavity by solutions of 
the Cauchy--Dirichlet problem for parabolic equations has been studied in several papers 
(see e.g. \cites{GK, INS, IST02, Ken, Kor, LV} and references therein). 
In this section we investigate 
sufficient conditions and necessary conditions for the preservation of $F$-concavity by classical solutions 
of the Cauchy--Dirichlet problem
\begin{equation}
\tag{N}
\left\{\begin{array}{ll}
\partial_tu=\displaystyle{\sum_{i,j=1}^na^{ij}(x,t)\partial_{x_i}\partial_{x_j}u}+{\mathcal G}(x,t,u,\nabla u)
\quad & \mbox{in}\quad\Omega\times(0,T),\vspace{3pt}\\
u(x,t)>0 &\mbox{in}\quad\Omega\times(0,T), \vspace{3pt}\\
u(x,t)=0 & \mbox{on}\quad \partial\Omega\times[0,T)\mbox{ if }\partial\Omega\not=\emptyset,\vspace{3pt}\\
u(x,0)=\phi(x)\ge 0 & \mbox{in}\quad\Omega,
\end{array}
\right.
\end{equation}
where $T\in(0,\infty]$, $\phi\in BC(\overline{\Omega})$, 
$a^{ij}\in BC(\Omega\times[0,T))$ with condition~(L2),
and ${\mathcal G}\in C(\Omega\times[0,T)\times(0,\infty)\times\mathbb{R}^n)$. 
A function 
$$
u\in C((\overline{\Omega}\times(0,T))\cup (\Omega\times[0,T))
$$
is called a {\it classical solution} of problem~(N) 
if $u\in C^{2;1}(\Omega\times(0,T))$ and $u$ satisfies problem~(N) pointwisely.  
See e.g. \cite{LSU}*{Chapter V} for the existence of classical solutions of problem~(N). 
\subsection{Sufficient conditions}\label{subsection:4.1}
In this subsection we develop the arguments of \cites{INS, IST01} to 
obtain sufficient conditions for the preservation of $F$-concavity by classical solutions of problem (N). 
Following the strategy in the proof of \cite{IST01}*{Theorem~3.1}, 
we will reduce the $F$-concavity of the solution~$u$ to the log-concavity of the function $v:=e^{F(u)}$. 

Let 
$$
\lambda\in\Lambda_{n+1}:=\left\{\lambda=(\lambda_1,\dots,\lambda_{n+1})\,\,\biggr|\,\,
\mbox{$0<\lambda_i<1$ for $i=1,\dots,n+1$},\,\,\,\sum_{i=1}^{n+1}\lambda_i=1\right\}. 
$$ 
Let $F$ be admissible on $I$ such that $\lim_{r\to +0}F(r)=-\infty$, 
and $\Omega$ a smooth, bounded, and convex domain in $\mathbb{R}^n$. 
Let $\phi\in BC_0(\overline{\Omega})$. 
Assume that a classical solution~$u$ of problem~(N) satisfies $u(\cdot,t)\in{\mathcal A}_\Omega(I)$ for all $t\in[0,T)$. 
Then we can define  the {\it spatially $F$-concave envelope} $u_F$ of $u$ by 
\begin{equation}
\label{eq:4.1}
u_F(x,t):=\sup_{\lambda\in\Lambda_{n+1}}\,u_{F,\lambda}(x,t)
\quad\mbox{for}\quad (x,t)\in\overline{\Omega}\times [0,T),
\end{equation}
where 
\begin{equation*}
u_{F,\lambda}(x,t):=\sup\left\{
f_F\left(\sum_{i=1}^{n+1}\lambda_i F(u(x_i,t))\right)
\,\,\biggr|\,\,\{x_i\}_{i=1}^{n+1}\subset\overline{\Omega},\,\,x=\sum_{i=1}^{n+1}\lambda_i x_i\right\}.
\end{equation*}
Notice that $u_F\geq u_{F,\lambda}\geq u$ for every $\lambda\in\Lambda_{n+1}$ and that, 
for any $t\in[0,T)$, $u(\cdot,t)$ is $F$-concave in $\Omega$ if and only if $u(\cdot,t)=u_F(\cdot,t)$ in $\Omega$,
since the function $u_F(\cdot,t)$ is the smallest $F$-concave function greater than or equal to $u(\cdot,t)$ in $\Omega$
(see e.g. \cite{sc}*{Theorem~1.1.4}). 
It follows from the convexity of $\Omega$ and $F(0)=-\infty$
that
\begin{equation}
\label{eq:4.2}
u_{F,\lambda}\in C(\overline{\Omega}\times[0,T)),\quad
u_{F,\lambda}>0\quad\mbox{in}\quad\Omega\times(0,T),
\quad
u_{F,\lambda}=0\quad\mbox{on}\quad\partial\Omega\times[0,T).
\end{equation}

We recall the notion of viscosity 
subsolution, supersolution, and solution of problem~(N). 
An upper semicontinuous function $U$ in $\Omega\times(0,T)$ 
is called a \emph{viscosity subsolution} of problem~(N) 
if, for any $(\xi,\tau)\in\Omega\times(0,T)$, 
the inequality 
$$
\partial_t\psi(\xi,\tau)\le \displaystyle{\sum_{i,j=1}^na^{ij}(\xi,\tau)(\partial_{x_i}\partial_{x_j}\psi)}(\xi,\tau)+{\mathcal G}(\xi,\tau,\psi(\xi,\tau),\nabla\psi(\xi,\tau))
$$
holds for every $C^{2;1}(\Omega\times(0,T))$ test function $\psi$ {\em touching $U$ from above at $(\xi,\tau)$}, i.e. satisfying 
$$
\mbox{$\psi(\xi,\tau)=U(\xi,\tau)$ and $\psi\ge U$ in a neighborhood of $(\xi,\tau)$}.
$$
Analogously, a lower semicontinuous function $U$ in $\Omega\times(0,T)$ 
is called a \emph{viscosity supersolution} of problem~(N)  
if, for any $(\xi,\tau)\in\Omega\times(0,T)$, 
the inequality 
$$
\partial_t\psi(\xi,\tau)\ge \displaystyle{\sum_{i,j=1}^na^{ij}(\xi,\tau)(\partial_{x_i}\partial_{x_j}\psi)}(\xi,\tau)+{\mathcal G}(\xi,\tau,\psi(\xi,\tau),\nabla\psi(\xi,\tau))
$$
holds for every $C^{2;1}(\Omega\times(0,T))$ test function $\psi$ {\em touching $U$ from below at $(\xi,\tau)$}, i.e. satisfying 
$$
\mbox{$\psi(\xi,\tau)=U(\xi,\tau)$ and $\psi\le U$ in a neighborhood of $(\xi,\tau)$.}
$$
A continuous function $U$ in $\Omega\times(0,T)$ is called a \emph{viscosity solution} of problem~(N)
if it is a~viscosity subsolution and supersolution of problem~(N) at the same time.
The main issue of this subsection is to prove that, under suitable assumptions, $u_{F,\lambda}$ is a viscosity subsolution of problem~(N).
\begin{proposition}
\label{Theorem:4.1}
Let $F$ be admissible on $I=[0,a)$ with $a\in(0,\infty]$, $\Omega$ a smooth, bounded, and convex domain in $\mathbb{R}^n$, and $T\in(0,\infty]$.
Let $u$ be a classical solution of problem~{\rm (N)} 
such that $u(\cdot,t)\in{\mathcal A}_\Omega(I)$ for $t\in[0,T)$ and $\phi\in BC_0(\overline{\Omega})$. 
Assume that 
$\lim_{r\to +0}F(r)=-\infty$, $F\in C^2({\rm int}\,I)$, $F'>0$ in ${\rm int}\,I$, 
and the following condition holds:
\begin{itemize}
 \item[{\rm (I)}] 
  for any $\theta\in\mathbb{R}^n$ and $t\in(0,T)$, the function ${\mathcal H}_{\theta,t}$ defined by 
  $$
  {\mathcal H}_{\theta,t}(x,z,M)\\
  :=\mbox{trace}\,(A(x,t)M)+
   \frac{{\mathcal G}\left(x,t,f_F(z),f_F'(z)\theta\right)}{f_F'(z)}+\left(\frac{f_F''(z)}{f_F'(z)}-1\right)
   \langle A(x,t)\theta,\theta\rangle 
  $$
  is concave with respect to $(x,z,M)\in\Omega\times J_F\times\mathrm{Sym}\,(n)$, where $A(x,t):=(a^{ij}(x,t))$.
\end{itemize}
Then $u_{F,\lambda}$ defined by \eqref{eq:4.1} is a viscosity subsolution of problem~{\rm (N)} 
for all $\lambda\in\Lambda_{n+1}$.    
\end{proposition}
{\bf Proof.}
It follows from the admissibility of $F$ and $\lim_{r\to +0}F(r)=-\infty$ 
that $J_F=(-\infty,F(a))$. 
Let $v=e^{F(u)}$. It turns out that
\begin{equation*}
\begin{split}
 & \partial_t v-\sum_{i,j=1}^n a^{ij}(x,t)\partial_{x_i}\partial_{x_j} v\\
 & =e^{F(u)}F'(u)\biggr(\partial_t u-\sum_{i,j=1}^n a^{ij}(x,t)\partial_{x_i}\partial_{x_j} u\biggr)
 -(e^{F(u)}F''(u)+e^{F(u)}F'(u)^2)\langle A(x,t)\nabla u,\nabla u\rangle\\
 & =e^{F(u)}F'(u){\mathcal G}(x,t,u,\nabla u)-(e^{F(u)}F''(u)+e^{F(u)}F'(u)^2)\langle A(x,t)\nabla u,\nabla u\rangle\\
 & =e^{F(u)}\frac{1}{f_F'(F(u))}{\mathcal G}(x,t,u,\nabla u)-\left(\frac{F''(u)}{F'(u)^2}+1\right)\frac{1}{v}\langle A(x,t)\nabla v,\nabla v\rangle\\
 & ={\mathcal G}_F(x,t,v,\nabla v)
\end{split}
\end{equation*}
in $\Omega\times(0,T)$, where 
\begin{equation*}
\begin{split}
{\mathcal G}_F(x,t,z,\theta)
:= & \frac{z}{f_F'(\log z)}\,{\mathcal G}\left(x,t,f_F(\log z),\frac{f_F'(\log z)\theta}{z}\right)\\
 & \qquad\qquad\qquad\qquad\qquad\qquad
 +\biggr(\frac{f_F''(\log z)}{f_F'(\log z)}-1\biggr)\frac{1}{z}\langle A(x,t)\theta,\theta\rangle
\end{split}
\end{equation*}
for $(x,t,z,\theta)\in \Omega\times(0,T)\times e^{J_F}\times\mathbb{R}^n$. 

For any $\theta\in\mathbb{R}^n$ and $t\in(0,T)$, 
it follows that
$$
{\mathcal H}_{\theta,t}(x,z,M)=
\mbox{trace}\,(A(x,t)M)+
e^{-z}{\mathcal G}_F(x,t,e^z,e^z\theta)
$$
for $(x,z,M)\in\Omega\times J_F\times\mathrm{Sym}\,(n)$. 
By condition~(I) we apply \cite{INS}*{Theorem~4.3} to see that, 
for any $\lambda\in\Lambda_{n+1}$, 
the function~$v_{\Phi_0,\lambda}$ defined by 
\begin{equation*}
v_{\Phi_0,\lambda}(x,t)
=
\sup\left\{\prod_{i=1}^{n+1} v(x_i,t)^{\lambda_i}
\,\biggr|\,\{x_i\}_{i=1}^{n+1}\subset\overline{\Omega},\,\,x=\sum_{i=1}^{n+1}\lambda_i x_i\right\}
\end{equation*}
is a viscosity subsolution of problem~(N)
with ${\mathcal G}$ replaced by ${\mathcal G}_F$.
This implies that $u_{F,\lambda}$ is a viscosity subsolution of problem~(N), 
and the proof is complete. 
$\Box$\vspace{5pt}

Assume that 
the following comparison principle holds: 
\begin{equation}
\tag{WCP}
\left\{
\begin{array}{l}
\mbox{Let $v$, $w\in C(\overline{\Omega}\times[0,T))$ be a classical solution and a viscosity}\vspace{3pt}\\
\mbox{subsolution of problem~(N), respectively, such that $v\geq w$ on $\overline{\Omega}\times\{0\}$}\vspace{3pt}\\
\mbox{and $v=w=0$ on $\partial\Omega\times[0,T)$. Then $v\geq w $ in $\overline{\Omega}\times[0,T)$.}
\end{array}
\right.
\end{equation}
See e.g.~\cite{UG}*{Section~8} for sufficient conditions of (WCP). 
Let $u$ be a classical solution of problem~(N) such that $\phi\in{\mathcal C}_\Omega[F]\cap BC_0(\overline{\Omega})$. 
Then 
$$
u_{F,\lambda}(\cdot,0)=\phi\quad\mbox{in $\overline{\Omega}$}\,\,\,\,\mbox{for all $\lambda\in\Lambda_{n+1}$}. 
$$
Thanks to \eqref{eq:4.2}, 
by Proposition~\ref{Theorem:4.1} and (WCP) we see that
$$
u\ge u_{F,\lambda} \quad\mbox{in $\overline{\Omega}\times[0,T)$}\,\,\,\,\mbox{for all $\lambda\in\Lambda_{n+1}$}.
$$
This yields $u\ge u_F$ in $\overline{\Omega}\times[0,T)$. 
Since $u\le u_F$ in $\overline{\Omega}\times[0,T)$ (see the definition of $u_F$), 
we find that $u= u_F$ in $\overline{\Omega}\times[0,T)$, 
that is, 
$u(\cdot,t)$ is $F$-concave in $\Omega$ for all $t\in [0,T)$. 
Then we have: 
\begin{proposition}
\label{Theorem:4.2}
Assume the same assumptions as in Proposition~{\rm\ref{Theorem:4.1}} 
and that {\rm (WCP)} holds. 
Let $u$ be a classical solution of problem~{\rm (N)} with the initial datum~$\phi\in {\mathcal C}_\Omega[F]\cap BC_0(\overline{\Omega})$. 
Then 
$$
u(\cdot,t)\in{\mathcal C}_\Omega[F]\quad\mbox{for all $t\in[0,T)$}. 
$$ 
\end{proposition}
%
\subsection{Necessary conditions}\label{subsection:4.2}
First, we give a necessary condition for the preservation of $F$-concavity by DHF in ${\mathbb R}^n$ with~$n\ge 2$. 
The proof of Proposition~\ref{Theorem:4.3} is a modification of the proof of Proposition~\ref{Theorem:3.2}. 
\begin{proposition}
\label{Theorem:4.3}
 {
Let $F$ be admissible on $I=[0,a)$ with $a\in(0,\infty]$ and $n\ge 2$. 
If $F$-concavity is preserved by DHF in $\mathbb{R}^n$, then
$$
\kappa\,\mathcal{C}_{{\mathbb R}^n}[F]\subset \mathcal{C}_{{\mathbb R}^n}[F]\quad\text{for every }\kappa\in(0,1)\,.
$$}
\end{proposition}
{\bf Proof.}
The statement is a consequence of the following slightly stronger property.
\begin{itemize}
  \item[(A)] 
  Assume that there exist $f\in{\mathcal C}_{{\mathbb R}^n}[F]$ and $\kappa\in(0,1)$ 
  such that $\kappa f\not\in{\mathcal C}_{{\mathbb R}^n}[F]$. 
  Then there exists $\phi\in \mathcal{C}_{{\mathbb R}^n}[F]\cap L^\infty({\mathbb R}^n)$ such that 
  $e^{t\Delta_{{\mathbb R}^n}}\phi$ is not $F$-concave in ${\mathbb R}^n$ for all small enough $t>0$. 
\end{itemize}
Let $f\in{\mathcal C}_{{\mathbb R}^n}[F]$ and $\kappa\in(0,1)$ such that 
$\kappa f\not\in{\mathcal C}_{{\mathbb R}^n}[F]$.
We find $\xi$, $\eta\in{\mathbb R}^n$ and $\lambda\in(0,1)$ such that 
\begin{equation}
\label{eq:4.3}
F(\kappa f((1-\lambda)\xi+\lambda \eta))<(1-\lambda)F(\kappa f(\xi))+\lambda F(\kappa f(\eta)).
\end{equation}
Since $F(0)=-\infty$, it follows that $f(\xi)>0$ and $f(\eta)>0$. 
By the concavity of $F(f)$ we see that 
$$
F(f((1-z)\xi+z \eta))\ge(1-z)F(f(\xi))+z F(f(\eta))>-\infty\quad\mbox{for}\quad z\in[0,1]. 
$$
It follows from $f\in{\mathcal C}_{{\mathbb R}^n}[F]$ that 
$f$ is continuous in the set $\{x\in{\mathbb R}^n\,|\,f(x)>0\}$. 
Then we find~$\delta>0$ such that 
$$
f((1-z)\xi+z \eta)>0\quad\mbox{for}\quad z \in(-\delta,1+\delta). 
$$

Set 
\begin{equation*}
\varphi(z):=
\left\{
\begin{array}{ll}
f((1-z )\xi+z  \eta)\quad & \mbox{for}\quad z \in(-\delta,1+\delta),\vspace{5pt}\\
0 & \mbox{for}\quad z\in{\mathbb R}\setminus (-\delta,1+\delta).
\end{array}
\right.
\end{equation*}
Then $\varphi$ is $F$-concave in $\mathbb{R}$. 
Furthermore, $\varphi$ is continuous in $(-\delta,1+\delta)$, and 
\begin{equation}
\label{eq:4.4}
\lim_{t\to +0}(e^{t\Delta_{\mathbb{R}}}\varphi)(z)=\varphi(z)\quad\mbox{for}\quad z \in(-\delta,1+\delta). 
\end{equation}
It follows from \eqref{eq:4.3} that 
\begin{equation}
\label{eq:4.5}
\begin{split}
F(\kappa \varphi(\lambda)) & =F(\kappa f((1-\lambda)\xi+\lambda \eta))\\
 & <(1-\lambda)F(\kappa f(\xi))+\lambda F(\kappa f(\eta))
=(1-\lambda)F(\kappa \varphi(0))+\lambda F(\kappa \varphi(1)).
\end{split}
\end{equation}
Combining \eqref{eq:4.4} and \eqref{eq:4.5}, we have
\begin{equation}
\label{eq:4.6}
F(\kappa   (e^{t\Delta_{\mathbb{R}}}\varphi) (\lambda) )
<(1-\lambda)F(\kappa  (e^{t\Delta_{\mathbb{R}}}\varphi) (0) )+\lambda F(\kappa  (e^{t \Delta_{\mathbb{R}}}\varphi) (1))
\end{equation}
for all small enough $t>0$. 

Set 
$$
\phi(w,z):={\bf 1}_{[0,\infty)}(w)\varphi(z)\quad\mbox{for}\quad(w,z)\in\mathbb{R}^2.
$$
Since $\varphi$ is $F$-concave in ${\mathbb R}$, 
we see that $\phi \in \mathcal{C}_{\mathbb{R}^2}[F]$. 
Furthermore, 
by \eqref{eq:2.3} we have
$$
(e^{t\Delta_{\mathbb{R}^2}} \phi)(w,z)
=(e^{t\Delta_{\mathbb{R}}}{\bf 1}_{[0,\infty)})(w)(e^{t\Delta_{\mathbb{R}}}\varphi)(z)
=h(t^{-\frac{1}{2}}w)(e^{t\Delta_{\mathbb{R}}}\varphi)(z)
$$
for $(w,z,t)\in\mathbb{R}^2\times(0,\infty)$.
In addition, by \eqref{eq:2.4} and \eqref{eq:2.5} 
we find a unique $\omega\in\mathbb{R}$  such that 
$h(\omega)=\kappa$. Then 
\begin{equation}
\label{eq:4.7}
(e^{t\Delta_{\mathbb{R}^2}} \phi)(t^{\frac{1}{2}}\omega,z)
=\kappa (e^{t\Delta_{\mathbb{R}}}\varphi)(z)\quad\mbox{for}\quad z\in{\mathbb R}.
\end{equation}
We deduce from \eqref{eq:4.6} and \eqref{eq:4.7} that
\begin{equation*}
\begin{split}
& F\left( (e^{t\Delta_{\mathbb{R}^2}} \phi)(  (1-\lambda) (t^{\frac{1}{2}}\omega,0)+\lambda (t^{\frac{1}{2}}\omega, 1) )\right)\\
 & =F\left( (e^{t \Delta_{\mathbb{R}^2}} \phi)(t^{\frac{1}{2}}\omega,\lambda )\right)
 =F\left(\kappa(e^{t\Delta_{\mathbb{R}}}\varphi)(\lambda)\right)\\
&<(1-\lambda)F(\kappa  (e^{t\Delta_{\mathbb{R}}}\varphi) (0) )+\lambda F(\kappa  (e^{t\Delta_{\mathbb{R}}}\varphi) (1))\\
&=
(1-\lambda)
F\left( (e^{t\Delta_{\mathbb{R}^2}} \phi)(t^{\frac{1}{2}}\omega,0 )\right)
+\lambda
F\left( (e^{t\Delta_{\mathbb{R}^2}} \phi)(t^{\frac{1}{2}}\omega,1 )\right)
\end{split}
\end{equation*}
for all small enough $t>0$.
This means that $e^{t\Delta_{\mathbb{R}^2}} \phi$ is not $F$-concave in ${\mathbb R}^2$ 
for all small enough~$t>0$. 
Thus property~(A) follows in the case $n=2$.
Similarly to \eqref{eq:3.15}, property~(A) also follows in the case $n\ge 3$. 
The proof is complete.
$\Box$\vspace{5pt}

Now let us consider problem~(N).
In~\cite{Kol}*{Theorem~1.2} Kolesnikov obtained necessary conditions
for log-concavity to be preserved by classical solutions 
of the Cauchy problem for linear parabolic equations under high regularity assumptions on the coefficients
(see also Section~\ref{subsection:6.1}). 
 {Hereafter},
we improve and develop the argument in the proof of \cite{Kol}*{Theorem~1.2}  
to obtain a necessary condition for the preservation of $F$-concavity 
by classical solutions of the Cauchy--Dirichlet problem~(N). 
\begin{proposition}
\label{Theorem:4.4}
Let $F$ be admissible on $I=[0,a)$ with $a\in(0,\infty]$ such that 
$F\in C^{2,\sigma}({\rm int}\,I)$ for some $\sigma\in[0,1)$, $\lim_{r\to +0}F(r)=-\infty$, and $F'>0\mbox{ in } {\rm int}\,I$ with $\liminf_{r\to+0}F'(r)>0$. 
Let $\Omega$ be a convex domain in $\mathbb{R}^n$. 
Assume that the following condition holds. 
\begin{itemize}
  \item[{\rm (II)}] 
  For any initial datum~$\phi\in {\mathcal C}_\Omega[F] \cap C^{2,\sigma}(\Omega)\cap BC^1(\overline{\Omega})$ with $\phi\not\equiv0$ in $\Omega$, 
  there exist $T\in(0,\infty]$ and a classical solution~$u\in C^{2;1}(\Omega\times[0,T))$ of problem~{\rm (N)} such that 
  $u(\cdot,t)\in{\mathcal C}_\Omega[F]$ for all $t\in(0,T)$.
\end{itemize}
Then, for any $\theta\in\mathbb{R}^n$ and $\ell\in{\mathbb R}$, the function $\widetilde{\mathcal{H}}$ defined by 
  $$
\widetilde{\mathcal{H}}(x)\\
  :=\frac{{\mathcal G}\left(x,0,f_F(\langle \theta,x\rangle+\ell),f_F'(\langle \theta,x\rangle+\ell)\theta\right)}{f_F'(\langle \theta,x\rangle+\ell)}+\frac{f_F''(\langle \theta,x\rangle+\ell)}{f_F'(\langle \theta,x\rangle+\ell)}
   \langle A(x,0)\theta,\theta\rangle 
  $$
must be concave in $\{x\in\Omega\,|\,\langle \theta,x\rangle+\ell\in J_F\}$, where $A(x,0):=(a^{ij}(x,0))$.
\end{proposition}
{\bf Proof.}
It follows from the admissibility of $F$ and $\lim_{r\to +0}F(r)=-\infty$ that $J_F=(-\infty,F(a))$. 
For any $b$, $c\in J_F$ with $b<c$, let $\psi\in C^3({\mathbb R})$ be concave in ${\mathbb R}$ such that 
$$
\psi(\xi)=\xi\,\,\,\mbox{for $\xi\in (-\infty,b]$},
\qquad \sup_{\xi\in{\mathbb R}}\psi\le c,
\qquad \sup_{\xi\in{\mathbb R}} |\psi'|<\infty.
$$
Let $B$ be an open ball in $\Omega$, $\theta\in{\mathbb R}^n$, and $\ell\in{\mathbb R}$ such that 
$\langle \theta,x\rangle+\ell\in (-\infty,b]$ for $x\in \overline{B}$. 
Set 
$$
\phi(x):=f_F(\psi(\langle \theta,x\rangle+\ell))\quad\mbox{for}\quad x\in\overline{\Omega}.
$$
Since $F\in C^{2,\sigma}({\rm int}\,I)$ and $F'>0\mbox{ in } {\rm int}\,I$ with $\liminf_{r\to+0}F'(r)>0$, we see that 
$f_F\in C^{2,\sigma}(J_F)$ and $f_F\in BC^1((-\infty,c])$. Then we observe that
$\phi\in {\mathcal C}_\Omega[F]\cap C^{2,\sigma}(\Omega)\cap BC^1(\overline{\Omega})$ and
$$
\phi(x)=f_F(\langle \theta,x\rangle+\ell)\quad\mbox{for}\quad x\in B. 
$$
By condition~(II) 
we find a unique classical solution~$u\in C^{2;1}(\Omega\times[0,T))$ of problem~(N) 
for some $T\in(0,\infty]$ such that $u(\cdot,t)\in{\mathcal C}_\Omega[F]$ for all $t\in(0,T)$. 
Then 
\begin{equation*}
\begin{split}
(\partial_t u)(x,0) & =f_F''(\langle \theta,x\rangle+\ell)\langle A(x,0)\theta,\theta\rangle 
+{\mathcal G}(x,0,f_F(\langle \theta,x\rangle+\ell),f_F'(\langle \theta,x\rangle+\ell)\theta)\\
 & =f_F'(\langle \theta,x\rangle+\ell)\widetilde{\mathcal{H}}(x)
\end{split}
\end{equation*}
for $x\in B$. 
Furthermore, 
$$
\Psi_\lambda(\xi,\eta,t):=F(u((1-\lambda)\xi+\lambda \eta,t))-(1-\lambda)F(u(\xi,t))-\lambda F(u(\eta,t))\ge 0
$$
for $\xi$, $\eta\in\Omega$, $t\in[0,T)$, and $\lambda\in(0,1)$. 
Since $F(\phi(x))=\langle \theta,x\rangle+\ell$ for $x\in B$, we see that 
$$
\Psi_\lambda(\xi,\eta,0)=0\quad\mbox{for $\xi$, $\eta\in B$ and $\lambda\in(0,1)$}. 
$$
These imply that 
\begin{equation}
\label{eq:4.8}
0\le(\partial_t\Psi_\lambda)(x,y,0)
=\widetilde{\mathcal{H}}((1-\lambda)x+\lambda y)
-(1-\lambda) \widetilde{\mathcal{H}}(x)-\lambda \widetilde{\mathcal{H}}(y)
\end{equation}
for $x$, $y\in B$.
Since $b$, $c$, and $B$ are arbitrary, 
we see that $\widetilde{\mathcal{H}}$ is concave in $\{x\in\Omega\,|\,\langle \theta,x\rangle+\ell\in J_F\}$. 
Thus Proposition~\ref{Theorem:4.4} follows.
$\Box$\vspace{5pt}
\newline
In the proof of Proposition~\ref{Theorem:4.4}, 
the assumption that $\liminf_{r\to+0}F'(r)>0$ is used only for proving that $\phi\in BC^1(\overline{\Omega})$. 
Then we can remove it at the price of strengthening slightly condition~(II) to obtain the following result,
whose proof is similar to that of Proposition~\ref{Theorem:4.4}. 
\begin{proposition}
\label{Theorem:4.5}
Let $F$ be admissible on $I=[0,a)$ with $a\in(0,\infty]$ such that 
$F\in C^{2,\sigma}({\rm int}\,I)$ for some $\sigma\in[0,1)$, $\lim_{r\to +0}F(r)=-\infty$, and $F'>0\mbox{ in } {\rm int}\,I$.
Let $\Omega$ be a convex domain in $\mathbb{R}^n$.  
Assume that the following condition holds. 
\begin{itemize}
  \item[{\rm (II')}] 
  For any initial datum~$\phi\in {\mathcal C}_\Omega[F]\cap C^{2,\sigma}(\Omega)\cap BC(\overline{\Omega})$ with $\phi\not\equiv0$ in $\Omega$, 
  there exist $T\in(0,\infty]$ and a classical solution~$u\in C^{2;1}(\Omega\times[0,T))$ of problem~{\rm (N)} such that 
  $u(\cdot,t)\in{\mathcal C}_\Omega[F]$ for all $t\in(0,T)$.
\end{itemize}
Then, for any $\theta\in\mathbb{R}^n$ and $\ell\in{\mathbb R}$, the function $\widetilde{\mathcal{H}}$ given in Proposition~{\rm\ref{Theorem:4.4}} 
is concave in $\{x\in\Omega\,|\,\langle \theta,x\rangle+\ell\in J_F\}$.
\end{proposition}

As a corollary of Proposition~\ref{Theorem:4.4}, 
we obtain a  necessary condition for the preservation of log-concavity 
by classical solutions of problem~(N) with $(a^{ij})=(\delta^{ij})$. 
\begin{corollary}
\label{Theorem:4.6}
Let $\Omega$ be a convex domain in $\mathbb{R}^n$. 
Assume that condition~{\rm (II)} in Proposition~{\rm\ref{Theorem:4.4}} holds for the parabolic equation
$$
\partial_tu=\Delta u+{\mathcal G}(x,u,\nabla u),
$$
where ${\mathcal G}\in C(\Omega\times(0,\infty)\times\mathbb{R}^n)$, 
with $I=[0,\infty)$ and $F=\Phi_0$. 
Then, for any $\theta\in\mathbb{R}^n$ and $\ell\in{\mathbb R}$, 
the function
$$
e^{-(\langle \theta,x\rangle+\ell)}{\mathcal G}\left(x,e^{\langle \theta,x\rangle+\ell},e^{\langle \theta,x\rangle+\ell}\theta\right)
$$ 
must be concave in $\Omega$.
\end{corollary}
As a direct consequence of Corollary~\ref{Theorem:4.6}, 
we obtain negative answers to the following question, 
which is motivated by the arguments in \cite{GK}*{Section~5}.
\begin{itemize}
  \item[({\bf Q4})] 
  Is log-concavity preserved by classical solutions of problem~(N) for the nonlinear parabolic equations listed below?
  \begin{equation*}
  \left\{
  \begin{array}{l}
  \partial_t u=\Delta u+\kappa u^p,\vspace{4pt}\\
  \partial_t u=\Delta u+\kappa e^u,\vspace{4pt}\\
  \partial_t u=\Delta u+\kappa' u^p\log u,\vspace{4pt}\\
  \partial_t u=\Delta u+\mu u^p+\kappa|\nabla u|^q,\vspace{4pt}\\
  \partial_t u=\Delta u+\langle {\bf b},\nabla u^p\rangle.
  \end{array}
  \right.
\end{equation*}
Here $p$, $q\in(1,\infty)$, $\kappa\in(0,\infty)$, $\kappa'\in{\mathbb R}\setminus\{0\}$, 
$\mu\in{\mathbb R}$, and ${\bf b}\in{\mathbb R}^n\setminus\{0\}$.
\end{itemize}
\begin{corollary}
\label{Theorem:4.7}
Let $\Omega$ be a convex domain in $\mathbb{R}^n$. 
Consider one of the nonlinear parabolic equations listed in {\rm ({\bf Q4})}. 
Then, for any $\sigma\in(0,1)$, 
there exist
$$
\phi\in {\mathcal C}_\Omega[\Phi_0]\cap C^{2,\sigma}(\Omega)\cap BC^1(\overline{\Omega})
$$ 
and a corresponding classical solution $u$ of problem~{\rm (N)} for some $T\in(0,\infty]$ 
such that $u(\cdot,t)$ is not log-concave in $\Omega$ for some $t\in(0,T)$. 
\end{corollary}
{\bf Proof.}
Consider one of the nonlinear parabolic equations listed in ({\bf Q4}).
Let $\sigma\in(0,1)$. 
For any $\phi\in C^{2,\sigma}(\Omega)\cap BC^1(\overline{\Omega})$, 
there {exists a} classical solution 
$u\in C^{2;1}(\Omega\times[0,T))$ of problem~(N) for some $T\in(0,\infty]$.
(See e.g. \cites{LSU, QS}.)
Since each of the functions
$$
e^{-z}{\mathcal G}\left(x,e^z,e^z \theta\right)=
\left\{
\begin{array}{l}
\kappa e^{(p-1)z},\vspace{3pt}\\
\kappa e^{-z}e^{e^z},\vspace{3pt}\\
\kappa' z e^{(p-1)z},\vspace{3pt}\\
\mu e^{(p-1)z}+\kappa e^{(q-1)z}|\theta|^q,\vspace{3pt}\\
pe^{(p-1)z}\langle {\bf b},\theta\rangle,
\end{array}
\right.
\quad\mbox{with}\quad z=\langle \theta,x\rangle+\ell,
$$ 
is not concave in $\Omega$ for some $\theta\in{\mathbb R}^n$ and $\ell\in{\mathbb R}$, 
by Corollary~\ref{Theorem:4.6} we see that condition~(II) does not hold. 
Consequently, Corollary~\ref{Theorem:4.7} follows. 
$\Box$
\vspace{5pt}

Next, we modify the argument in the proof of Proposition~\ref{Theorem:4.4} to obtain 
a necessary condition for $F$-concavity to be preserved by DHF in $\Omega$.
\begin{proposition}
\label{Theorem:4.8}
Let $F$ be admissible on $I=[0,a)$ with $a\in(0,\infty]$ such that $F\in C^2({\rm int}\,I)$.
Let~$\Omega$ be a convex domain in $\mathbb{R}^n$. 
Assume that $F$-concavity is preserved by DHF in $\Omega$. 
Then $F'>0$ in ${\rm int}\,I$ and $(\log f_F')'$ is concave in $J_F$.
\end{proposition}
{\bf Proof.}
Since $F$ is strictly increasing in ${\rm int}\,I$, for any $\delta\in(0,a)$, 
there exists $r_*\in(0,\delta)$ such that $F'(r_*)>0$. 
Then we find $r_1$, $r_2\in [0,a]$ with $r_1<r_*<r_2$ such that 
$$
F'(r)>0\quad\mbox{for}\quad r\in(r_1,r_2).
$$
Set $I':=(r_1,r_2)$ and $J':=F(I')$.  
Let $\theta\in{\mathbb R}^n$. 
Let $B$ be an open ball with $\overline{B}\subset\Omega$ and $\kappa\in{\mathbb R}$ 
such that 
$$
\langle \theta,x\rangle+\kappa\in J'\quad\mbox{for}\quad x\in B.
$$
Set 
\begin{equation}
\label{eq:4.9}
\phi(x):=f_F(\langle \theta,x\rangle+\kappa){\bf 1}_B(x)\quad\mbox{for}\quad x\in\Omega. 
\end{equation}
Then $\phi\in {\mathcal C}_\Omega[F]\cap L^\infty(\Omega)$, which together with the preservation of $F$-concavity by DHF implies that 
\begin{equation}
\label{eq:4.10}
e^{t\Delta_\Omega}\phi\in {\mathcal C}_\Omega[F]\quad\mbox{for}\quad t>0.
\end{equation}
For any open ball $B'$ with $\overline{B'}\subset B$, let $\psi\in C_0^\infty(\Omega)$ be such that 
$0\le\psi\le 1$ in $\Omega$, $\psi=1$ in $B'$, and $\mbox{supp}\,\psi\subset B$.  
It follows that
$$
(e^{t\Delta_\Omega}\phi)(x)
=\int_\Omega G_{\Delta_\Omega}(x,y,t)\psi(y)\phi(y)\,dy+\int_\Omega G_{\Delta_\Omega}(x,y,t)(1-\psi(y))\phi(y)\,dy
\quad\mbox{for}\quad x\in\Omega.
$$
Since $\psi\phi\in C^2(\Omega)$ and $(1-\psi)\phi=0$ in $B'$, 
we see that 
\begin{equation}
\label{eq:4.11}
(e^{t\Delta_\Omega}\phi)\in C^{2;1}(B'\times[0,\infty)).
\end{equation}
As in the proof of Proposition~\ref{Theorem:4.4}, we introduce a function
$$
\Psi_\lambda(x,y,t):=F(e^{t\Delta_\Omega}\phi((1-\lambda)x+\lambda y))-(1-\lambda)F(e^{t\Delta_\Omega}\phi(x))-\lambda F(e^{t\Delta_\Omega}\phi(y))
$$
for $x$, $y\in\Omega$, $t\ge 0$, and $\lambda\in(0,1)$. By \eqref{eq:4.9} and \eqref{eq:4.10}, for any $\lambda\in(0,1)$, we have
$$
\Psi_\lambda(x,y,t)\ge 0\quad\mbox{for}\quad (x,y,t)\in\Omega\times\Omega\times(0,\infty),
\quad \Psi_\lambda(x,y,0)=0\quad\mbox{for}\quad (x,y)\in B\times B.
$$
This together with \eqref{eq:4.11} implies that 
$(\partial_t\Psi_\lambda)(x,y,0)\ge 0$ for $x$, $y\in B'$ and $\lambda\in(0,1)$.
Then, by the same argument as in \eqref{eq:4.8} we see that
$$
(\log f_F')'((1-\lambda)z+\lambda w)\ge(1-\lambda)(\log f_F')'(z)+\lambda (\log f_F')'(w)
$$
for $x$, $y\in B'$ and $\lambda\in(0,1)$, where $z:=\langle \theta,x\rangle+\kappa$ and $w:=\langle \theta,y\rangle+\kappa$. 
{Since $\kappa$, $B$, and $B'$ are arbitrary}, we see that 
the function $(\log f_F')'$ is concave in $J'$. 

Assume that $r_2<a$ and $F'(r_2)=0$. It follows from the concavity of $(\log f_F')'$ in $J'$ that 
$$
\limsup_{z\to F(r_2)-0}\, (\log f_F')'(z)<\infty,
$$
which implies that 
$$
\lim_{z\to F(r_2)-0}\log f'_F(z)-\log f'_F(\zeta)=\lim_{z\to F(r_2)-0}\int_\zeta^z (\log f_F'(w))'\,dw<\infty
\quad\mbox{for}\quad \zeta\in J'.
$$
On the other hand, since $F'(f_F(z))f'_F(z)=1$ for $z\in J'$ and $F'(r_2)=0$, 
we see that 
$$
\lim_{z\to F(r_2)-0}f_F'(z)=\lim_{z\to F(r_2)-0} \frac{1}{F'(f_F(z))}=\infty.
$$
This is a contradiction. Then we see that $F'(r)>0$ for $r\in(r_1,a)$.  
Since $\delta$ is arbitrary and $r_1\in(0,\delta)$, 
we deduce that $F'(r)>0$ for $r\in(0,a)$.
Then, setting $I'={\rm int}\,I$ and $J'=J_F$, 
we see that $F'>0$ in ${\rm int}\,I$ and $(\log f_F')'$ is concave in $J_F$. The proof is complete. 
$\Box$
\vspace{5pt}

Further necessary conditions for the preservation of $F$-concavity by solutions of parabolic equations
are discussed in Section~6.
\section{Main theorems}\label{section:5}
In this section we complete the proofs of Theorem \ref{MainThm} and Theorem \ref{Theorem:1.8}.
The proof of the former is better split in some steps which have their own interest and are enucleated in the following two theorems.
\begin{theorem}
\label{Theorem:1.5}
Let $I=[0,a)$ with $a\in(0,\infty]$ and  $\Omega$ a convex domain in ${\mathbb R}^n$ with $n\ge 1$. 
\begin{itemize}
  \item[{\rm (1)}]
  $H_a$-concavity is preserved by DHF in $\Omega$. 
  \item[{\rm (2)}] 
  Let $F$ be admissible on $I$. 
  If $F$-concavity is preserved by DHF in $\Omega$, then 
  $F$-concavity is weaker than $H_a$-concavity in ${\mathcal A}_\Omega(I)$ and $\lim_{r\to +0}F(r)=-\infty$. 
\end{itemize}
\end{theorem}
{\bf Proof.}
Let $a\in(0,\infty]$. 
The proof is divided into three steps.
\vspace{3pt}
\newline
{\bf Step 1}:
Consider the case where $\Omega$ is a smooth, bounded, and convex domain in ${\mathbb R}^n$. 
Let $\phi\in {\mathcal C}_\Omega[H_a]\cap BC_0(\overline{\Omega})$. 
For any $\theta\in{\mathbb R}^n$ and $t>0$, 
let ${\mathcal H}_{\theta,t}$ be as in Proposition~\ref{Theorem:4.1} with ${\mathcal G}=0$ 
and $(a^{ij})=(\delta^{ij})$. 
In the case $a\in(0,\infty)$, 
since the function $ah$ is the inverse function of $H_a$ in ${\mathbb R}$, 
for any $\theta\in{\mathbb R}^n$, 
we observe from Lemma~\ref{Theorem:2.9} that
$$
{\mathcal H}_{\theta,t}(x,z,M) 
=\text{trace}\,(M)+\left(\frac{ah''(z)}{ah'(z)}-1\right)|\theta|^2
=\text{trace}\,(M)+\left(-\frac{1}{2}z-1\right)|\theta|^2
$$
for $(x,z,M)\in\Omega\times{\mathbb R}\times\mathrm{Sym}\,(n)$. 
Similarly, in the case $a=\infty$, 
since $e^z$ is the inverse function of $H_a$, we have 
$$
{\mathcal H}_{\theta,t}(x,z,M)=\text{trace}\,(M)
\quad\mbox{for}\quad (x,z,M)\in\Omega\times{\mathbb R}\times\mathrm{Sym}\,(n).
$$
These imply that ${\mathcal H}_{\theta,t}$ is concave with respect to $(x,z,M)\in\Omega\times{\mathbb R}\times\mathrm{Sym}\,(n)$ 
for any $\theta\in{\mathbb R}^n$ and $t>0$.
Then it follows from Proposition~\ref{Theorem:4.2} that
$e^{t\Delta_\Omega}\phi$ is $F$-concave in $\Omega$ for all $t>0$. 
\vspace{5pt}
\newline
{\bf Step 2}: 
Let $\Omega$ be a convex domain in ${\mathbb R}^n$ and $\phi\in{\mathcal C}_\Omega[H_a]\cap L^\infty(\Omega)$. 
Then there exists a sequence of smooth, bounded, and convex domains $\{\Omega_\ell\}$ 
such that 
$$
\Omega_1\subset\Omega_2\subset\cdots\subset\Omega_\ell\subset\cdots,
\qquad
\bigcup_{\ell=1}^\infty\Omega_\ell=\Omega.
$$
(See e.g. \cite{sc}*{Theorem 2.7.1}.)
By Lemma~\ref{Theorem:3.4} 
we find a sequence $\{\phi_j\}\subset {\mathcal C}_{\Omega_\ell}[H_a]\cap BC_0(\overline{\Omega_\ell})$ 
such that 
$$
\lim_{j\to\infty}\left(e^{t\Delta_{\Omega_\ell}}\phi_j\right)(x)=\left(e^{t\Delta_{\Omega_\ell}}\phi\right)(x)
\quad\mbox{for}\quad (x,t)\in\Omega_\ell\times(0,\infty).
$$
By Step 1 we see that $e^{t\Delta_{\Omega_\ell}}\phi_j$ is $H_a$-concave in $\Omega_\ell$ for all $t>0$, 
consequently so is $e^{t\Delta_{\Omega_\ell}}\phi$.  
On the other hand, we observe that 
$$
\lim_{\ell\to \infty}\left(e^{t\Delta_{\Omega_\ell}}\phi\right)(x)=\left(e^{t\Delta_{\Omega}}\phi\right)(x)
\quad\mbox{for}\quad (x,t)\in\Omega\times(0,\infty).
$$
Then we conclude that $e^{t\Delta_\Omega}\phi$ is $H_a$-concave in $\Omega$ for all $t>0$. 
The proof of assertion~(1) is complete. 
\vspace{5pt}
\newline
{\bf Step 3}: 
We prove assertion~(2). 
Let $F$ be admissible on $I=[0,a)$. 
Assume that $F$-concavity is preserved by DHF in $\Omega$. 
Then Proposition~\ref{Theorem:3.3} implies that 
$F$-concavity is also preserved by DHF in ${\mathbb R}^n$. 
Consider the case of $a<\infty$. 
Let $a'\in(0,a)$ and set 
$$
\phi(x):=a'{\bf 1}_{[0,\infty)}(\langle x, e_1\rangle)\quad\mbox{for}\quad x\in{\mathbb R}^n. 
$$
Then $\phi\in{\mathcal A}_{{\mathbb R}^n}(I)$ and $\phi$ is $F$-concave in ${\mathbb R}^n$. 
This together with the preservation of $F$-concavity by DHF in ${\mathbb R}^n$ implies that 
$e^{\Delta_{{\mathbb R}^n}}\phi$ is $F$-concave in ${\mathbb R}^n$. 
Since 
$$
\left(e^{\Delta_{{\mathbb R}^n}}\phi\right)(x)=a'h(\langle x, e_1\rangle)\in(0,a')
\quad\mbox{for}\quad x\in{\mathbb R}^n,
$$
we see that $F(a'h)$ is concave in ${\mathbb R}$. 
Letting $a'\to a$, we observe that $F(ah)$ is concave in ${\mathbb R}$.
This together with Lemma~\ref{Theorem:2.4} implies that
${\mathcal C}_{{\mathbb R}^n}[H_a]
\subset{\mathcal C}_{{\mathbb R}^n}[F]$. 
Then we deduce from Lemma~\ref{Theorem:2.6} that
${\mathcal C}_\Omega[H_a]\subset{\mathcal C}_\Omega[F]$. 
Furthermore, by Lemma~\ref{Theorem:2.7} we see that $\lim_{r\to +0}F(r)=-\infty$. 
Thus assertion~(2) holds in the case $a<\infty$. 

Assertion~(2) in the case $a=\infty$ follows from \cite{IST03}*{Theorem~1.1}. 
Here we give another proof for the sake of completeness of this paper. 
By assertion~(2) with $a<\infty$ we see that 
\begin{equation}
\label{eq:5.1}
{\mathcal C}_\Omega[H_k]\subset{\mathcal C}_\Omega[F]\cap{\mathcal A}_\Omega([0,k))
\quad\mbox{for all $k>0$}
\end{equation}
and $\lim_{r\to +0}F(r)=-\infty$. 
Let $f$ be log-concave in $\Omega$. 
For any $m>0$, since $\min\{f,m\}$ is log-concave in $\Omega$, 
by Lemma~\ref{Theorem:2.10} we find a sequence $\{f_k\}$ such that $f_k\in {\mathcal C}_\Omega[H_k]$ and 
$$
\lim_{k\to\infty}f_k(x)=\min\{f(x),m\}\quad\mbox{for}\quad x\in\Omega.
$$
This together with \eqref{eq:5.1} yields $f\in {\mathcal C}_\Omega[F]$, 
which implies that that log-concavity is stronger than $F$-concavity in ${\mathcal A}_\Omega([0,\infty))$. 
Thus assertion~(2) holds in the case $a=\infty$, and the proof of Theorem~\ref{Theorem:1.5} is complete. 
$\Box$
\begin{theorem}
\label{Theorem:1.6}
Let $F$ be admissible on $I=[0,a)$ with $a\in(0,\infty]$ and  
$\Omega$ a convex domain in ${\mathbb R}^n$ with $n\ge 2$. 
\begin{itemize}
  \item[{\rm (1)}] 
  Assume that $F$-concavity is preserved by DHF in $\Omega$.
  Then 
    \begin{itemize}
    \item[{\rm (i)}]
    $F$-concavity is stronger than $\log$-concavity in ${\mathcal A}_\Omega(I)$;
    \item[{\rm (ii)}]
    if a function $f$ is $F$-concave in $\Omega$, then so is $\kappa f$ for $\kappa\in(0,1)$. 
    \end{itemize}
  \item[{\rm (2)}] 
  If $F$-concavity is not stronger than log-concavity in ${\mathcal A}_\Omega(I)$, 
  then there exists a bounded continuous function $\phi$ on $\overline{\Omega}$ 
  with the following properties:
  \begin{itemize}
  \item[$\bullet$]
   $\phi$ is $F$-concave in $\Omega$ and $\phi=0$ on $\partial\Omega$ if $\partial\Omega\not=\emptyset$;
  \item[$\bullet$]
   $e^{t\Delta_\Omega}\phi$ is not quasi-concave in $\Omega$ for some $t>0$.
\end{itemize}
\end{itemize}
\end{theorem}
{\bf Proof.}
By Propositions~\ref{Theorem:3.2} and \ref{Theorem:3.5} we obtain 
Theorem~\ref{Theorem:1.6}~(1)-(i) and (2). 
It remains to prove Theorem~\ref{Theorem:1.6}~(1)-(ii). 
Let $\Omega$ be a convex domain in ${\mathbb R}^n$ with $n\ge 2$. 
Assume that 
there exist $f\in{\mathcal C}_\Omega[F]$ and $\kappa\in(0,1)$ satisfying $\kappa f\not\in{\mathcal C}_\Omega[F]$. 
Setting $f=0$ outside $\Omega$, we see that $f\in{\mathcal C}_{{\mathbb R}^n}[F]$ and $\kappa f\not\in{\mathcal C}_{{\mathbb R}^n}[F]$.
Then Proposition~\ref{Theorem:4.3} implies that $F$-concavity is not preserved by DHF in ${\mathbb R}^n$. 
Consequently, by Proposition~\ref{Theorem:3.3} we see that 
$F$-concavity is not preserved by DHF in $\Omega$.
Therefore we observe that if $F$-concavity is preserved by DHF in $\Omega$, 
then $\kappa{\mathcal C}_\Omega[F]\subset{\mathcal C}_\Omega[F]$ for all $\kappa\in(0,1)$. 
Thus Theorem~\ref{Theorem:1.6}~(1)-(ii) follows. The proof of Theorem~\ref{Theorem:1.6} is complete.
$\Box$
\begin{remark} 
Theorem~{\rm\ref{Theorem:1.6}~(1)-(ii)} also gives a necessary condition for $F$-concavity to be preserved by~DHF 
and leads  {to} a nice property of the inverse function $f_F$  {of $F$} {\rm ({\it see Proposition}~\ref{Theorem:6.1}~(3))}, 
which plays an important role in Section~{\rm 6}.
Notice that if $F$-concavity possesses property~{\rm (1)-(ii)} of Theorem~{\rm \ref{Theorem:1.6}} {\em for all $\kappa\in(0,\infty)$},
then $F$-concavity coincides with some power concavity 
{\rm ({\it see Lemma}~\ref{Theorem:2.8})}. 
\end{remark}
%
%
%
%
Now we can proceed to the proofs first of Theorem \ref{Theorem:1.8} and then of Theorem \ref{MainThm}.
\vspace{3pt}
\newline
{\bf Proof of Theorem~\ref{Theorem:1.8}.} 
We prove assertion~(1). 
Let $F$ be admissible on $I=[0,a)$ with $a\in(0,\infty]$. 
Assume that $F$-concavity is preserved by DHF in $\Omega$. 
Then we apply Theorem~\ref{Theorem:1.5}~(2) to obtain $\lim_{r\to +0}F(r)=-\infty$. 
Furthermore, it follows from Proposition~\ref{Theorem:4.8} that $F'>0$ in ${\rm int}\,I$ and the concavity of $(\log f_Z')'$ in $J_F$. 

Conversely, we assume that $\lim_{r\to +0}F(r)=-\infty$, $F'>0$ in ${\rm int}\,I$, and the function $(\log f_F')'$ is concave in~$J_F$. 
By~Proposition~\ref{Theorem:4.2} we see that 
$e^{t\Delta_\Omega}\phi$ is $F$-concave in $\Omega$ for all $t>0$ when $\Omega$ is a smooth, bounded, and convex domain 
and $\phi\in{\mathcal C}_\Omega[F]\cap BC_0(\overline{\Omega})$. 
Then, repeating the same argument as in the proof of Theorem~\ref{Theorem:1.5} (see Step 2), we deduce that 
$F$-concavity is preserved by DHF in convex domains in~${\mathbb R}^n$. 
The proof is complete.~$\Box$
\vspace{3pt}

\begin{remark}
Let us remark again that, to our knowledge, Theorem~{\rm\ref{Theorem:1.8}} is the first result regarding a {\em necessary} and sufficient condition 
for concavity properties of solutions to partial differential equations. 
As an application of Theorem~{\rm\ref{Theorem:1.8}~(1)}, we can also characterize $\alpha$-log-concavity preserved by DHF in~$\Omega$
{\rm ({\it see Corollary}~\ref{Theorem:6.11})}.
\end{remark}

\noindent{\bf Proof of Theorem~\ref{MainThm}.} 
{By Theorem~\ref{Theorem:1.5}
we obtain ({\bf A1}).
Theorem~\ref{Theorem:1.6}~(1)-(i) yields ({\bf A2}) for $n\geq 2$, while Theorem~\ref{Theorem:1.6}~(2) gives ({\bf A3}).}
Theorem \ref{Theorem:1.8} (together with Proposition \ref{Theorem:1.1}~(1)) gives the one-dimensional part of ({\bf A2}).
The proof is complete.~$\Box$


\section{Further study of the preservation of $F$-concavity}\label{section:6}
In this section, based on the arguments of the previous sections, 
we study the preservation of $F$-concavity by solutions of the Cauchy--Dirichlet problem 
for linear parabolic equations with variable coefficients (see Section~\ref{subsection:6.1});
semilinear heat equations (see Section~\ref{subsection:6.2});
the porous medium equation and the parabolic $p$-Laplace equation (see Section~\ref{subsection:6.3}). 
\subsection{Linear parabolic equation with variable coefficients}\label{subsection:6.1}
In this subsection we consider  {a (slightly) simplified version of problem~(P) (treated in Section~\ref{section:3}), 
precisely, the following Cauchy--Dirichlet problem 
for a linear parabolic equation} with variable (but independent of $t$) coefficients: 
\begin{equation}
\tag{P'}
\left\{
\begin{array}{ll}
\partial_t u=Lu & \mbox{in}\quad\Omega\times(0,\infty),\vspace{3pt}\\
u=0 & \mbox{in}\quad\partial\Omega\times(0,\infty)\text{ if }\partial \Omega \neq \emptyset,\vspace{3pt}\\
u(\cdot,0)=\phi\ge 0 & \mbox{in}\quad\Omega,
\end{array}
\right.
\end{equation}
where $\Omega$ is a convex domain in ${\mathbb R}^n$ and $\phi\in L^\infty(\Omega)$. 
Here $L$ is an elliptic operator of the form 
$$
L:=\sum_{i,j=1}^na^{ij}(x)\partial_{x_i}\partial_{x_j}+\sum_{i=1}^n b^i(x)\partial_{x_i}-c(x),
$$
and the coefficients satisfy the following conditions:
\begin{itemize}
  \item[(L1')] 
  there exists $\sigma\in(0,1)$ such that 
  $$
  a^{ij},\,\,b^i\in C^{0,\sigma}(\Omega),\qquad
  c\in BC^{0,\sigma}(\Omega),
  $$ 
  where $i$, $j=1,\dots,n$;
  \item[(L2')] 
  $A(x)=(a^{ij}(x))\in\mathrm{Sym}\,(n)$ for $x\in\Omega$ 
  and there exists $\Lambda>0$ such that 
  $$
  \Lambda^{-1}|\xi|^2\le \langle A(x)\xi,\xi\rangle\le\Lambda|\xi|^2
  \quad\mbox{for all $\xi\in{\mathbb R}^n$ and $x\in\Omega$}.
  $$
\end{itemize}
The preservation of concavity property by solutions of problem~(P') 
has been studied in several papers (see e.g. \cites{AI, Borell, BL, GK, INS, Kol}), 
and the following properties hold.
\begin{itemize}
  \item 
  Let $L=\Delta-c(x)$. If $c$ is nonnegative and convex in $\Omega$, 
  then log-concavity is preserved by the Dirichlet parabolic flow associated with $L$ in $\Omega$ (see e.g. \cites{Borell, BL, INS}).
  \item 
  Let $a^{ij}$, $b^i\in C^{2,\sigma}({\mathbb R}^n)$ for some $\sigma\in(0,1)$ and $c=0$ in ${\mathbb R}^n$. 
  Then log-concavity is preserved by  the Dirichlet parabolic flow associated with $L$ in ${\mathbb R}^n$ 
  if and only if
  $A=(a^{ij})$ is constant in ${\mathbb R}^n$ and $b^i$ is affine in ${\mathbb R}^n$ (see \cite{Kol}). 
\end{itemize}
Below, under conditions~(L1') and (L2'), 
we develop the arguments in the previous sections to characterize the $F$-concavities preserved by the Dirichlet parabolic flow.
\begin{proposition}
\label{Theorem:6.1}
Let $F$ be admissible on $I=[0,a)$ with $a\in(0,\infty]$ such that $F\in C^2({\rm int}\,I)$. 
Let~$\Omega$ be a convex domain in ${\mathbb R}^n$. 
Assume that $F$-concavity is preserved by the Dirichlet parabolic flow associated with $L$ in $\Omega$ 
under conditions~{\rm (L1')} and {\rm (L2')}. 
Then the following conditions hold.
\begin{itemize}
  \item[{\rm (1)}] 
   {$F$-concavity is weaker than $H_a$-concavity and stronger than log-concavity} in ${\mathcal A}_\Omega(I)$.
  \item[{\rm (2)}] 
  $\lim_{r\to +0}F(r)=-\infty$, $F'>0$ in ${\rm int}\,I$, and $(\log f_F')'$ is concave in $J_F$.
  \item[{\rm (3)}]
  $(\log f_F)'$ is $(-1)$-concave in $J_F$. 
\end{itemize}
\end{proposition}
{\bf Proof.}
Assume that $F$-concavity is preserved by the Dirichlet parabolic flow associated with~$L$ in $\Omega$. 
By Proposition~\ref{Theorem:3.3} we see that
$F$-concavity is also preserved by DHF in ${\mathbb R}^n$. 
Then conditions~(1) and (2) follow from Theorem~\ref{Theorem:1.5}~(2), Theorem~\ref{Theorem:1.6}~(1)-(i), and Theorem~\ref{Theorem:1.8}.
Furthermore, 
we deduce from Theorem~\ref{Theorem:1.6}~(1)-(ii) and Theorem~\ref{Theorem:1.8}~(2) that 
$$
e^{-t}e^{t\Delta_{{\mathbb R}^n}}\phi\in {\mathcal C}_{{\mathbb R}^n}[F]\quad\mbox{for $t>0$ if $\phi\in{\mathcal C}_{{\mathbb R}^n}[F]$}.
$$
This implies that 
$F$-concavity is preserved by the Dirichlet parabolic flow associated with $\Delta -1$ in ${\mathbb R}^n$.
Then, by Proposition~\ref{Theorem:4.4} we see that, 
for any $\theta\in{\mathbb R}^n$ and $\ell\in{\mathbb R}$, 
the function 
$$
-\frac{f_F(\langle\theta,x\rangle+\ell)}{f'_F(\langle\theta,x\rangle+\ell)}+\frac{f''_F(\langle\theta,x\rangle+\ell)}{f'_F(\langle\theta,x\rangle+\ell)}|\theta|^2
$$
is concave in $\{x\in\Omega\,|\,\langle\theta,x\rangle+\ell\in J_F\}$, that is,
the function
$$
-\frac{f_F(z)}{f'_F(z)}+\frac{f''_F(z)}{f'_F(z)}|\theta|^2
$$
is concave with respect to $z\in J_F$
for any $\theta \in \mathbb{R}^n\setminus\{0\}$, which implies that $(\log f_F)'$ is $(-1)$-concave in $J_F$. 
Thus condition~(3) follows. 
The proof is complete.
$\Box$\vspace{5pt}

We focus on the case that $c=0$ in $\Omega$, 
and obtain a necessary and sufficient condition for $F$-concavity to be preserved by the Dirichlet parabolic flow.
\begin{theorem}
\label{Theorem:6.2}
Let $F$ be admissible on $I=[0,a)$ with $a\in(0,\infty]$ such that $F\in C^2({\rm int}\,I)$. 
Let $\Omega$ be a convex domain in ${\mathbb R}^n$ and the elliptic operator $L$ satisfy conditions~{\rm (L1')} and {\rm (L2') 
with $c=0$ in $\Omega$.} 
Then $F$-concavity is preserved by the Dirichlet parabolic flow associated with $L$ in~$\Omega$ 
if and only if the following conditions hold:
\begin{itemize}
  \item[{\rm (1)}] $\,\displaystyle{(\log f_F')'(z)\langle A(x)\theta,\theta\rangle}$ are concave with respect to $(x,z)\in\Omega\times J_F$ 
  for any $\theta \in \mathbb{R}^n$;
  \item[{\rm (2)}] $b^i$ is affine in $\Omega$ for $i=1,\dots,n$;
  \item[{\rm (3)}]
  $\lim_{r\to +0}F(r)=-\infty$, $F'>0$ in ${\rm int}\,I$, and $(\log f_F')'$ is concave in $J_F$. 
\end{itemize}
\end{theorem}
{\bf Proof.}
Assume that $F$-concavity is preserved by the Dirichlet parabolic flow associated with~$L$ in $\Omega$ under conditions~(L1') and (L2'). 
By Proposition~\ref{Theorem:6.1}~(2) we see that condition~(3) holds. 
Furthermore, thanks to parabolic regularity theorems (see e.g. \cite{LSU}*{Chapter~IV, Theorem~10.1}), 
we see that assumption~(II') in Proposition~\ref{Theorem:4.5} holds. 
Then Proposition~\ref{Theorem:4.5} implies that, 
for any $\theta=(\theta_1,\dots,\theta_n)\in\mathbb{R}^n$ and $\ell\in{\mathbb R}$, 
the function
\begin{equation}
\label{eq:6.1}
\sum_{i=1}^n b^i(x) {\theta_i}+(\log f_F')'(\langle\theta,x\rangle+\ell)
\langle A(x)\theta,\theta\rangle
\end{equation}
is concave in $\{x\in\Omega\,|\,\langle\theta,x\rangle+\ell\in J_F\}$. 
Then we observe that the function $\sum_{i=1}^n b^i(x) {\theta_i}$ is concave with respect to $x\in\Omega$ 
for any $\theta=(\theta_1,\dots,\theta_n)\in\mathbb{R}^n$, 
which implies that $b^i$ is concave and convex in $\Omega$ for $i=1,\dots,n$. 
Thus condition~(2) holds. 
Furthermore, by \eqref{eq:6.1} we see that,
for any $\theta\in{\mathbb R}^n$, 
\begin{equation}
\label{eq:6.2}
\mbox{$\,\displaystyle{(\log f_F')'(\langle\theta,x\rangle+\ell)\langle A(x)\theta,\theta\rangle}$ are concave 
in $\{x\in\Omega\,|\,\langle\theta,x\rangle+\ell\in J_F\}$}.
\end{equation}
Then we obtain condition~(1). 
Indeed, if not, we find $\lambda\in(0,1)$, $x$, $y\in\Omega$, and $z$, $w\in J_F$ such that 
\begin{equation}
\label{eq:6.4}
\begin{split}
 & (\log f_F')'((1-\lambda)z+\lambda w)\langle A((1-\lambda)x+\lambda y)\theta,\theta\rangle\\
 & <(1-\lambda)(\log f_F')'(z)\langle A(x)\theta,\theta\rangle+\lambda (\log f_F')'(w)\langle A(y)\theta,\theta\rangle.
\end{split}\end{equation}
Thanks to the continuity of $A$, we can assume that $x\not=y$.
Then we find $\theta\in{\mathbb R}$ such that 
$z-w=\langle x-y,\theta\rangle$. 
Setting $\ell=z-\langle \theta,x\rangle$, we obtain
$$
z=\langle \theta,x\rangle+\ell,\quad w=\langle \theta,y\rangle+\ell.
$$
This together with \eqref{eq:6.4} implies that \eqref{eq:6.2} does not hold, which is a contradiction.
Thus condition~(1) holds.

Conversely, under conditions~(1)--(3), 
by Proposition~\ref{Theorem:4.2} we see that, 
if $\Omega$ is a smooth, bounded, and convex domain in ${\mathbb R}^n$,
and $\phi\in {\mathcal C}_\Omega[F]\cap BC_0(\overline{\Omega})$, 
then $e^{tL_\Omega}\phi\in {\mathcal C}_\Omega[F]$ for $t>0$. 
For any convex domain $\Omega$, 
we apply the same argument as in Step 2 of the proof of Theorem~\ref{Theorem:1.5} 
to see that $F$-concavity is preserved by the Dirichlet parabolic flow associated with $L$ in $\Omega$.
Thus Theorem~\ref{Theorem:6.2} follows.
$\Box$
\vspace{3pt}

Similarly, 
we obtain a sufficient and necessary condition for log-concavity to be preserved by the Dirichlet parabolic flow 
in every convex domain in ${\mathbb R}^n$ under conditions~{\rm (L1')} and {\rm (L2')}.
\begin{theorem}
\label{Theorem:AA} 
Let the elliptic operator $L$ satisfy conditions~{\rm (L1')} and {\rm (L2')} with $\Omega$ replaced by~${\mathbb R}^n$. 
Assume that 
log-concavity is preserved by the Dirichlet parabolic flow associated with $L$ in ${\mathbb R}^n$ under conditions~{\rm (L1')} and {\rm (L2')}. 
Then the following conditions hold:
\begin{itemize}
  \item[{\rm (1)}] the matrix $A$ is constant in ${\mathbb R}^n$;
  \item[{\rm (2)}] $b^i$ is affine in ${\mathbb R}^n$ for $i=1,\dots,n$;
  \item[{\rm (3)}] $c$ is convex in ${\mathbb R}^n$. 
\end{itemize}
Conversely, under conditions~{\rm (1)}, {\rm (2)}, and {\rm (3)}, 
log-concavity is preserved by the Dirichlet parabolic flow associated with $L$ in every convex domain in ${\mathbb R}^n$. 
\end{theorem}
{\bf Proof.}
Assume that log-concavity is preserved by the Dirichlet parabolic flow associated with $L$ in ${\mathbb R}^n$. 
Proposition~\ref{Theorem:4.5} implies that, 
for any $\theta=(\theta_1,\dots,\theta_n)\in\mathbb{R}^n$, 
the function
\[
\sum_{i=1}^n b^i(x) {\theta_i}-c(x)+\langle A(x)\theta,\theta\rangle
\]
is concave in ${\mathbb R}^n$. 
Then, for any $\theta\in{\mathbb R}^n$, 
$\langle A(x)\theta,\theta\rangle$ is concave and nonnegative in ${\mathbb R}^n$. 
This implies condition~(1). 
Furthermore, we see that
for any $\theta\in{\mathbb R}^n$, $\sum_{i=1}^n b^i(x) {\theta_i}$ is concave in ${\mathbb R}^n$, 
and we obtain condition~(2). Then condition~(3) also holds.

Conversely, under conditions~(1), (2), (3), 
applying the same argument as in Step 2 of the proof of Theorem~\ref{Theorem:1.5}, 
we see that log-concavity is preserved by the Dirichlet parabolic flow associated with $L$ in every convex domain in ${\mathbb R}^n$. 
Thus Theorem~\ref{Theorem:AA} follows.
$\Box$
\subsection{Semilinear heat equation}\label{subsection:6.2}
In this subsection we consider the Cauchy--Dirichlet problem 
for a semilinear heat equation
\begin{equation}
\tag{SH}
\left\{\begin{array}{ll}
\partial_tu=\Delta u+\kappa|u|^{p-1}u
\quad & \mbox{in}\quad\Omega\times(0,T),\vspace{3pt}\\
u(x,t)=0 & \mbox{on}\quad \partial\Omega\times[0,T)\text{ if }\partial \Omega \neq \emptyset\,,\vspace{3pt}\\
u(x,0)=\phi(x)\ge 0& \mbox{in}\quad\Omega,
\end{array}
\right.
\end{equation}
where $T\in(0,\infty]$, $\phi\in L^\infty(\Omega)$, $\kappa\in{\mathbb R}$, and $p>1$. 
Problem~(SH) possesses a unique classical $L^\infty(\Omega)$-solution~$S_\Omega(\cdot)\phi$ for some $T\in(0,\infty]$ 
(see \cite{QS}*{Section~15} 
for the existence and the uniqueness of classical $L^\infty(\Omega)$-solutions of problem~(SH)).  
Then 
$$
S_\Omega(\cdot)\phi\in BC^{2;1}(\Omega\times(0,T))\cap C(\overline{\Omega}\times(0,T)),
\quad
\lim_{t\to+0}\|S_\Omega(t)\phi-e^{t\Delta_\Omega}\phi\|_{L^\infty(\Omega)}=0,
$$ 
and  
\begin{align*}
(S_\Omega(t)\phi)(x) & =\int_\Omega G_{\Delta_\Omega}(x,y,t)\phi(y)\,dy\\
 & +\kappa\int_0^t\int_\Omega G_{\Delta_\Omega}(x,y,t-s) |(S_\Omega(s)\phi)(y)|^{p-1}
(S_\Omega(s)\phi)(y)\,dy\,ds
\end{align*}
for all $(x,t)\in\Omega\times(0,T)$. 
Let $z$ be a solution of the ODE $z'=\kappa |z|^{p-1}z$ with $z(0)=\|\phi\|_{L^\infty(\Omega)}$, 
that is, 
$$
z(t)=\|\phi\|_{L^\infty(\Omega)}\left(1-\kappa(p-1)t\|\phi\|_{L^\infty(\Omega)}^{p-1}\right)^{-\frac{1}{p-1}}. 
$$
Let $T_M$ (resp.\,\,$T_M^*$) be the maximal existence time of the solution~$S_\Omega(\cdot)\phi$ (resp.\,\,$z$).
Then it follows from the comparison principle that 
\begin{eqnarray}
\notag
 & & 0\le (S_\Omega(t)\phi)(x)\le z(t)\quad\mbox{for}\quad(x,t)\in\Omega\times (0,T_M^*),\\
\label{eq:6.5}
 &  & T_M\ge T_M^*=\frac{1}{\kappa(p-1)}\|\phi\|_{L^\infty(\Omega)}^{-(p-1)}\quad\mbox{if}\quad \kappa>0,
 \quad
 T_M=T_M^*=\infty\quad\mbox{if}\quad\kappa\le 0.
\end{eqnarray}
Notice that, if $\kappa>0$, 
then it does not necessarily hold that $T_M=\infty$ (see e.g. \cites{F, W}). 
We also observe from  the strong maximum principle that 
$(S_\Omega(t)\phi)(x)>0$ for $(x,t)\in\Omega\times(0,T_M)$ if $\phi\not=0$ in $L^\infty(\Omega)$.

Problem~(SH) has been studied from various points of view 
(see e.g. the monograph~\cite{QS} and references therein).  
Behavior of solutions of problem~(SH) depends on the exponent $p$, the sign of $\kappa$, the behavior of the initial datum $\phi$, and 
the shape of the domain~$\Omega$.
In the case where $\kappa>0$ and $p>1+2/n$, 
the large time behavior of solutions with $\Omega={\mathbb R}^n$ depends 
on the size of the initial datum (see e.g. \cites{Kawanago, Souplet}). 
On the other hand, in the case $\kappa<0$, 
the large time behavior of solutions with $\Omega={\mathbb R}^n$ varies widely with the behavior of the initial datum $\phi$ at the space infinity and the sign of $p-(1+2/n)$
(see e.g. \cites{GV, KP, HL}). 

Concavity properties of solutions of problem~(SH) have been studied only in the case $\kappa<0$, 
where the preservation of log-concavity by solutions of problem~(SH) has been established  {in} \cites{GK, INS}. 
In the case $\kappa>0$, 
Corollary~\ref{Theorem:4.7} implies that 
log-concavity is not preserved by classical solutions of problem~(SH). 
The aim of this section is to characterize  {$F$-concavities, in particular, the strongest and the weakest ones}, 
preserved by solutions of problem~(SH). 
\begin{definition}
\label{Theorem:6.5}
Let $F$ be admissible on $I=[0,a)$ with $a\in(0,\infty]$ and  $\Omega$
a convex domain in~${\mathbb R}^n$. 
We say that \emph{$F$-concavity is preserved by solutions of problem~(SH)} if 
$$
\mbox{$S_\Omega(t)\phi \in \mathcal{C}_{\Omega}[F]$ holds 
for all $\phi\in\mathcal{C}_{\Omega}[F]\cap L^\infty(\Omega)$ and $t\in(0,T_M)$}.
$$
\end{definition}
%

We start with proving that the disruption of $F$-concavity by DHF 
implies the same by solutions of problem~(SH).
\begin{proposition}
\label{Theorem:6.6}
Let $F$ be admissible on $I=[0,a)$ with $a\in(0,\infty]$ and 
$\Omega$ a convex domain in~${\mathbb R}^n$.  
If $F$-concavity is not preserved by DHF in ${\mathbb R}^n$, 
then there exists $\phi\in{\mathcal C}_\Omega[F]\cap BC_0(\overline{\Omega})$  
such that $S_\Omega(t)\phi$ is not $F$-concave in $\Omega$ for some $t\in(0,T_M)$. 
\end{proposition}
{\bf Proof.}
Thanks to \eqref{eq:1.2}, 
we can assume, without loss of generality, that $0\in\Omega$. Set $\Omega_\ell:=\ell\Omega$ for $\ell=1,2,\dots$. 
Assume that there exists $\phi\in {\mathcal C}_{{\mathbb R}^n}[F]\cap L^\infty({\mathbb R}^n)$ 
such that $e^{\tau\Delta_{{\mathbb R}^n}}\phi$ is not $F$-concave in ${\mathbb R}^n$ for some $\tau>0$, 
that is, 
\begin{equation}
\label{eq:6.6}
F((e^{\tau\Delta_{{\mathbb R}^n}}\phi)((1-\lambda)\xi+\lambda\eta))
<(1-\lambda)F((e^{\tau\Delta_{{\mathbb R}^n}}\phi)(\xi))+\lambda F((e^{\tau\Delta_{{\mathbb R}^n}}\phi)(\eta))
\end{equation}
for some $\xi$, $\eta\in{\mathbb R}^n$ and $\lambda\in(0,1)$. 
By the same argument as in the proof of Proposition~\ref{Theorem:3.3} 
we find a sequence~$\{\ell_j\}\subset\{\ell\}$ with $\lim_{j\to\infty}\ell_j=\infty$ such that 
there exists $\phi_{\ell_j}\in{\mathcal C}_{\Omega_{\ell_j}}[F]\cap BC_0(\overline{\Omega_{\ell_j}})$ 
satisfying 
\begin{align}
\label{eq:6.7}
 & \sup_j \|\phi_{\ell_j}\|_{L^\infty(\Omega_{\ell_j})}\le\|\phi\|_{L^\infty({\mathbb R}^n)},\\
\notag
 & \lim_{j\to\infty}\left(e^{t\Delta_{\Omega_{\ell_j}}}\phi_{\ell_j}\right)(x)=(e^{t\Delta_{{\mathbb R}^n}}\phi)(x)
\quad\mbox{uniformly on compact sets in ${\mathbb R}^n\times(0,\infty)$}.
\end{align}
We write $\Omega_j:=\Omega_{\ell_j}$ and $\phi_j:=\phi_{\ell_j}$ for simplicity. 

Since $\ell_j\to \infty$ as $j\to\infty$, 
taking large enough $j$, by \eqref{eq:6.5} and \eqref{eq:6.7} we find a unique $L^\infty(\Omega_j)$-classical solution $v_j$ of the problem
\begin{equation}
\label{eq:6.8}
\left\{\begin{array}{ll}
\partial_tv=\Delta v+\kappa\ell_j^{-2} |v|^{p-1}v
\quad & \mbox{in}\quad\Omega_j\times(0,1),\vspace{3pt}\\
v=0\quad & \mbox{on}\quad\partial\Omega_j\times(0,1)\text{ if }\partial \Omega_j \neq \emptyset,\vspace{3pt}\\
v(x,0)=\phi_j(x)& \mbox{in}\quad\Omega_j\,.
\end{array}
\right.
\end{equation}
Set
$$
\gamma_j^\pm(t):=\left(1\mp|\kappa|\ell_j^{-2}(p-1)\int_0^t \|e^{s\Delta_{\Omega_j}}\phi_j\|_{L^\infty(\Omega_j)}^{p-1}\,ds\right)^{-\frac{1}{p-1}},
\quad
V_j^\pm(x,t):=\gamma_j^\pm(t)\left(e^{t\Delta_{\Omega_j}}\phi_j\right)(x),
$$
for $(x,t)\in\Omega_j\times(0,1)$. 
Since  
$$
\partial_t V_j^\pm-\Delta V_j^\pm=\pm|\kappa|\ell_j^{-2}\gamma^\pm_j(t)^p\|e^{t\Delta_{\Omega_j}}\phi_j\|_{L^\infty(\Omega_j)}^{p-1}e^{t\Delta_{\Omega_j}}\phi_j
\quad\mbox{in}\quad\Omega_j\times(0,1),
$$
we have
$$
\partial_t V_j^+-\Delta V_j^+\ge|\kappa|\ell_j^{-2}|V_j^+|^{p-1}V_j^+,
\qquad
\partial_t V_j^--\Delta V_j^+\le -|\kappa|\ell_j^{-2}|V_j^-|^{p-1}V_j^-,
$$
in $\Omega\times(0,1)$, 
that is, $V_j^+$ (resp.\,\,$V_j^-$) is a supersolution (resp.\,\,subsolution) of problem~\eqref{eq:6.8}. 
Furthermore, 
$$
v_j=V_j^\pm
\quad\mbox{on $\partial\Omega_j\times(0,1)$ if $\partial\Omega_j\not=\emptyset$ and on $\Omega_j\times\{0\}$}. 
$$
Then the comparison principle for problem~\eqref{eq:6.8} implies that 
\begin{equation}
\label{eq:6.9}
V_j^-(x,t)\le v_j(x,t)\le V_j^+(x,t)\quad\mbox{in}\quad\Omega_j\times(0,1).
\end{equation}
On the other hand, by \eqref{eq:6.7} we have 
\begin{equation}
\label{eq:6.10}
|\gamma_j^\pm(\tau)-1|\le C\ell_j^{-2}\int_0^\tau \|\phi_j\|_{L^\infty(\Omega_j)}^{p-1}\,ds\to 0\quad\mbox{as}\quad j\to\infty,
\end{equation}
where $\tau$ is as in \eqref{eq:6.6}. 
Combining \eqref{eq:6.6}, \eqref{eq:6.9}, and \eqref{eq:6.10}, 
we obtain 
\begin{equation}
\label{eq:6.11}
F(v_j((1-\lambda)\xi+\lambda\eta,\tau))<(1-\lambda)F(v_j(\xi,\tau))+\lambda F(v_j(\eta,\tau))
\end{equation}
for large enough $j$.

Set $u_j(x,t):=v_j(\ell_j x,\ell_j^2t)$ for $(x,t)\in\Omega\times[0,\ell_j^{-2})$. 
Then $u_j$ is a solution of problem~(SH)  
such that $u_j(\cdot,0)\in {\mathcal C}_\Omega[F]\cap BC_0(\overline{\Omega})$.
Furthermore, it follows from \eqref{eq:6.11} that 
\begin{align*}
  F\left(u_j((1-\lambda)\xi_j+\lambda\eta_j,\tau_j)\right)
    <(1-\lambda)F(u_j(\xi_j,\tau_j))+\lambda F(u_j(\eta_j,\tau_j))
\end{align*}
for large enough $j$, where $\xi_j:=\ell_j^{-1}\xi$, $\eta_j:=\ell_j^{-1}\eta$, and $\tau_j:=\ell_j^{-2}\tau$. 
These mean that $F$-concavity is not preserved by solutions of problem~(SH). 
Thus Proposition~\ref{Theorem:6.6} follows. 
$\Box$\vspace{5pt}
\newline
Similarly to the proof of Proposition~\ref{Theorem:6.6}, we see
that the disruption of quasi-concavity by DHF implies the same by solutions of problem~(SH).
\begin{proposition}
\label{Theorem:6.7}
Let $F$ be admissible on $I=[0,a)$ with $a\in(0,\infty]$ and 
$\Omega$ a convex domain in~${\mathbb R}^n$.  
Assume that there exists $\phi\in{\mathcal C}_{{\mathbb R}^n}[F]\cap L^\infty({\mathbb R}^n)$ 
such that $e^{\tau\Delta_{{\mathbb R}^n}}\phi$ is not quasi-concave in~${\mathbb R}^n$ for some $\tau>0$.
Then there exists $\psi\in{\mathcal C}_\Omega[F]\cap BC_0(\overline{\Omega})$  
such that $S_\Omega(t)\psi$ is not quasi-concave in $\Omega$ for some $t\in(0,T_M)$. 
\end{proposition}
Combining Propositions~\ref{Theorem:6.1}, \ref{Theorem:6.6}, and \ref{Theorem:6.7}
with Theorems~\ref{Theorem:1.5} and \ref{Theorem:1.6}, we have:  
\begin{theorem}
\label{Theorem:6.8}
Let $F$ be admissible on $I=[0,a)$ with $a\in(0,\infty]$ 
and $\Omega$ a convex domain in~${\mathbb R}^n$. 
\begin{itemize}
  \item[{\rm (1)}] 
  Assume that $F$-concavity is preserved by solutions of problem~{\rm (SH)}.
  Then the following properties hold.
  \begin{itemize}
  \item[{\rm (i)}] 
  $F$-concavity is weaker than $H_a$-concavity in ${\mathcal A}_\Omega(I)$ and $\lim_{r\to+0}F(r)=-\infty$.
   \item[{\rm (ii)}] 
   Let $n\ge 2$. Then $F$-concavity is stronger than log-concavity in ${\mathcal A}_\Omega(I)$. 
   Furthermore, 
   if a function $f$ is $F$-concave in $\Omega$, then so is $\kappa f$ for $\kappa\in(0,1)$. 
   \item[{\rm (iii)}] 
   If $F\in C^2({\rm int}\,I)$, then $F$ satisfies conditions~{\rm (1)}--{\rm (3)} of Proposition~{\rm\ref{Theorem:6.1}}.
  \end{itemize}
  \item[{\rm (2)}] 
  If $F$-concavity is not stronger than log-concavity in ${\mathcal A}_\Omega(I)$ and $n\ge 2$, 
  then there exists $\phi\in {\mathcal C}_\Omega[F]\cap BC_0(\overline{\Omega})$ 
  such that $S_\Omega(t)\phi$ is not quasi-concave in $\Omega$ for some $t\in(0,T_M)$.
\end{itemize}
\end{theorem}
{\bf Proof.}
Assume that $F$-concavity is preserved by solutions of problem~(SH). 
Proposition~\ref{Theorem:6.6} implies that $F$-concavity is preserved by DHF in ${\mathbb R}^n$. 
Then, combing Theorem~\ref{Theorem:1.5}, Theorem~\ref{Theorem:1.6}, and Proposition~\ref{Theorem:6.1} with Lemma~\ref{Theorem:2.6},
we obtain assertion~(1). 
Assertion~(2) follows from Theorem~\ref{Theorem:1.6} and Proposition~\ref{Theorem:6.7}. 
The proof is complete.  
$\Box$
\vspace{5pt}

Now we are in position to state the main results of this subsection: 
Theorem~\ref{Theorem:6.9}  and Theorem~\ref{Theorem:6.10}. The former concerns with the case of $\kappa>0$ and the latter with $\kappa\leq 0$.
\begin{theorem}
\label{Theorem:6.9}
Let $F$ be admissible on $I=[0,a)$ with $a\in(0,\infty]$ such that $F\in C^2({\rm int}\,I)$. 
Let~$\Omega$ be a convex domain in~${\mathbb R}^n$ and $\kappa>0$. 
Then $F$-concavity is not preserved by solutions of problem~{\rm (SH)}.
\end{theorem}
{\bf Proof.}
Assume that $F$-concavity is preserved by solutions of problem~(SH). 
It follows from Theorem~\ref{Theorem:1.8} and Proposition~\ref{Theorem:6.6} that 
$F$-concavity is preserved by DHF in $\Omega$. 
This together with Theorem~\ref{Theorem:1.8}~(1) yields
$\lim_{r\to +0}F(r)=-\infty$ and $F'>0\mbox{ in } {\rm int}\,I$. 
By parabolic regularity theorems (see e.g. \cite{LSU}*{Chapter~IV}) 
we see that assumption~(II') in Proposition~\ref{Theorem:4.5} holds for problem (SH). 
Then, by Proposition~\ref{Theorem:4.5}, 
for any $\theta\in{\mathbb R}^n$ and $\ell\in{\mathbb R}$, 
we see that 
$$
\kappa\frac{f_F(\langle \theta,x\rangle+\ell)^p}{f_F'(\langle \theta,x\rangle+\ell)}+\frac{f''_F(\langle \theta,x\rangle+\ell)}{f'_F(\langle \theta,x\rangle+\ell)}|\theta|^2
$$
is concave in $\{x\in \Omega\,|\,\langle \theta,x\rangle+\ell\in J_F\}$. 
This implies that $f_F(z)^p/f_F'(z)$
is concave with respect to $z\in J_F$.
Then we see that 
\begin{equation}
\label{eq:6.12}
\left(\frac{f_F(z)^p}{f_F'(z)}\right)'=f_F(z)^{p-1}\left(p-\frac{f_F(z)f_F''(z)}{f_F'(z)^2}\right)
\quad\mbox{is non-increasing in $J_F$}.
\end{equation}
On the other hand, it follows from Theorem~\ref{Theorem:6.8}~(1)-(iii) (see also Proposition~\ref{Theorem:6.1}~(3)) that 
$f_F/f_F'$ is convex in $J_F$, which implies that 
\begin{equation}
\label{eq:6.13}
\left(\frac{f_F(z)}{f_F'(z)}\right)'=1-\frac{f_F(z)f''_F(z)}{f_F'(z)^2}\quad\mbox{is non-decreasing in $J_F$}.
\end{equation}
Since $f_F$ is strictly increasing in $J_F$, 
by \eqref{eq:6.12} and \eqref{eq:6.13} we see that 
$$
\left(\frac{f_F}{f_F'}\right)'=-(p-1)\quad\mbox{in}\quad J_F,
$$
that is, there exists $C_1\in{\mathbb R}$ such that 
$$
\frac{f_F(z)}{f_F'(z)}=-(p-1)z+C_1\quad\mbox{for}\quad z\in J_F.
$$
Since $f_F/f_F'>0$ in $J_F$, we see that $-(p-1)F(a)+C_1\ge 0$. 
Then there exists $C_2>0$ such that 
$$
f(z)=C_2(-(p-1)z+C_1)^{-\frac{1}{p-1}}\quad\mbox{for}\quad z\in J_F,
$$
that is, 
$$
F(r)=-\frac{1}{p-1}\left(\frac{r}{C_2}\right)^{-(p-1)}+\frac{C_1}{p-1}\quad\mbox{for}\quad r\in I.
$$
Lemma~\ref{Theorem:2.5} implies that 
$F$-concavity coincides with $-(p-1)$-concavity in ${\mathcal A}_\Omega(I)$. 
On the other hand, Theorem~\ref{Theorem:6.8}~(1)-(iii) implies that $F$-concavity is stronger than log-concavity in ${\mathcal A}_\Omega(I)$. 
This is a contradiction (see Example~\ref{Theorem:2.1}~(1)). 
Thus $F$-concavity is not preserved by solutions of problem~(SH), 
and Theorem~\ref{Theorem:6.9} follows.
$\Box$\vspace{5pt}

In the case $\kappa\le 0$, we have:
\begin{theorem}
\label{Theorem:6.10}
Let $F$ be admissible on $I=[0,a)$ with $a\in(0,\infty]$ and $F\in C^2({\rm int}\,I)$.
Let $\Omega$ be a convex domain in~${\mathbb R}^n$ and $\kappa\le 0$. 
Then $F$-concavity is preserved by solutions of problem~{\rm (SH)} if and only if $\lim_{r\to +0}F(r)=-\infty$, 
$F'>0$ in ${\rm int}\,I$, and 
\begin{equation}
\label{eq:6.14}
\mbox{the functions $\mathcal{F}_1:=(\log f_F')'$ and 
$\mathcal{F}_2:=\kappa\displaystyle{\frac{f_F^p}{f_F'}}$ 
are concave in $J_F$}.
\end{equation}
In particular, 
\begin{itemize}
  \item[{\rm (1)}] 
  log-concavity is preserved by solutions of problem~{\rm (SH)}:
  \item[{\rm (2)}] 
  $\alpha$-log-concavity is preserved by solutions of problem~{\rm (SH)} if and only if $\alpha\in[1/2,1]$:
  \item[{\rm (3)}] 
  for any $a\in(0,\infty]$, $H_a$-concavity is preserved by solutions of problem~{\rm (SH)}.
\end{itemize}
\end{theorem}
{\bf Proof.}
Assume that $\lim_{r\to +0}F(r)=-\infty$, $F'>0$ in ${\rm int}\,I$, and \eqref{eq:6.14} holds. 
It follows from \eqref{eq:6.14} that 
the function
$$
\kappa\frac{f_F(z)^p}{f_F'(z)}+\frac{f_F''(z)}{f_F'(z)}|\theta|^2
=\mathcal{F}_2(z)+\mathcal{F}_1(z)|\theta|^2
$$
is concave with respect to $z\in J_F$ for any $\theta\in{\mathbb R}^n$. 
Then Proposition~\ref{Theorem:4.2} implies that 
$F$-concavity is preserved by solutions of problem~(SH) 
in the case $\Omega$ is a smooth bounded convex domain. 
Then, by the same argument as in Step~2 of the proof of Theorem~\ref{Theorem:1.5} 
we obtain the preservation of $F$-concavity by solutions of problem~(SH) for a generic convex domain.

Conversely, assume that $F$-concavity is preserved by solutions of problem~(SH). 
Then, similarly to the proof of Theorem~\ref{Theorem:6.9}, we see that, 
for any $\theta\in{\mathbb R}^n$ and $\ell\in{\mathbb R}$,  
the function
$$
\kappa\frac{f_F(\langle \theta,x\rangle+\ell)^p}{f_F'(\langle \theta,x\rangle+\ell)}+\frac{f_F''(\langle \theta,x\rangle+\ell)}{f_F'(\langle \theta,x\rangle+\ell)}|\theta|^2
=\mathcal{F}_2(\langle \theta,x\rangle+\ell)+\mathcal{F}_1(\langle \theta,x\rangle+\ell)|\theta|^2
$$
is concave in $\{x\in\Omega\,|\,\langle \theta,x\rangle+\ell\in J_F\}$.
Since $\theta$ and $\ell$ are arbitrary, we observe that each of the functions $\mathcal{F}_1$ and $\mathcal{F}_2$ is concave in~$J_F$, 
that is, \eqref{eq:6.14} holds. 
Furthermore, it follows from Theorem~\ref{Theorem:6.8}~(1)-(iii) that $\lim_{r\to +0}F(r)=-\infty$ and $F'>0$ in ${\rm int}\,I$. 

It remains to prove assertions~(1)--(3). 
We prove assertion~(1). 
It suffices to prove that \eqref{eq:6.14} holds 
in the case $F=\Phi_0$:
since $f_F(z)=e^z$ for $z\in J_F={\mathbb R}$, 
we have
$$
{\mathcal F}_1(z)=1,\quad {\mathcal F}_2(z)=\kappa e^{(p-1)z},\quad\mbox{for}\quad z\in{\mathbb R},
$$
which are concave in ${\mathbb R}$ since $\kappa\leq 0$. 
Thus \eqref{eq:6.14} holds, and assertion~(1) follows.  

We prove assertion~(2). 
We write $f_\alpha:=f_{L_\alpha}$ and $J_\alpha:=J_{L_\alpha}$ for simplicity.  
By the above arguments we see that 
$\alpha$-log-concavity is preserved by solutions of problem~(SH) if and only if 
$\lim_{r\to +0}L_\alpha(r)=-\infty$, $L_\alpha'>0$ in ${\rm int}\,I$, and the functions 
$(\log f_\alpha')'$ and $\kappa f_\alpha^p/f_\alpha'$ are concave in $J_\alpha$.

If $\alpha<0$, then $L_\alpha(r)\to 1/\alpha$ as $r\to +0$. 
If $\alpha=0$, since
$$
f_\alpha(z)=e^{-e^{-z}},
\quad f_\alpha'(z)=e^{-e^{-z}}e^{-z},
\quad f_\alpha''(z)=e^{-e^{-z}}e^{-2z}-e^{-e^{-z}}e^{-z},
\quad
\frac{f_\alpha''(z)}{f_\alpha'(z)}=e^{-z}-1,
$$
for $z\in J_\alpha={\mathbb R}$, 
we see that $(\log f_\alpha')'$ is not concave in $J_\alpha$. 
Therefore, by \eqref{eq:6.14} we see that $\alpha$-log-concavity is not preserved by solutions of problem~(SH) if $\alpha\le 0$. 

Let $\alpha>0$. Then $J_\alpha=(-\infty,1/\alpha)$ and 
$$
f_\alpha(z)=\exp\left(-(1-\alpha z)^{\frac{1}{\alpha}}\right)
\quad\mbox{for}\quad z\in J_\alpha.
$$
It follows that 
\begin{align*}
 & f'_\alpha(z)=(1-\alpha z)^{-1+\frac{1}{\alpha}}f_\alpha(z),\\
 & f''_\alpha(z)=\left((\alpha-1)(1-\alpha z)^{-2+\frac{1}{\alpha}}+(1-\alpha z)^{-2+\frac{2}{\alpha}}\right)f_\alpha(z),\\
 & {\mathcal F}_1(z)=(\alpha-1)(1-\alpha z)^{-1}+(1-\alpha z)^{-1+\frac{1}{\alpha}},\\
 & {\mathcal F}_2(z)=\kappa(1-\alpha z)^{1-\frac{1}{\alpha}}f_\alpha(z)^{p-1},
\end{align*}
for $z\in J_\alpha$.
Then
\begin{align*}
 & {\mathcal F}_1'(z)=\alpha(\alpha-1)(1-\alpha z)^{-2}+(\alpha-1)(1-\alpha z)^{-2+\frac{1}{\alpha}},\\
 & {\mathcal F}_1''(z)=2\alpha^2(\alpha-1)(1-\alpha z)^{-3}+(\alpha-1)(2\alpha-1)(1-\alpha z)^{-3+\frac{1}{\alpha}},\\
 & {\mathcal F}_2'(z)=-\kappa(\alpha-1)(1-\alpha z)^{-\frac{1}{\alpha}}f_\alpha(z)^{p-1}+\kappa(p-1)f_\alpha(z)^{p-1},\\
 & {\mathcal F}_2''(z)=-\kappa(\alpha-1)(1-\alpha z)^{-1-\frac{1}{\alpha}}f_\alpha(z)^{p-1}\\
  & \qquad\qquad
  -\kappa(\alpha-1)(p-1)(1-\alpha z)^{-1}f_\alpha(z)^{p-1}+\kappa(p-1)^2(1-\alpha z)^{-1+\frac{1}{\alpha}}f_\alpha(z)^{p-1},
\end{align*}
for $z\in J_\alpha$. 
We observe that ${\mathcal F}_1''\le 0$ in $J_\alpha$ if and only if $\alpha\in[1/2,1]$. 
Furthermore, if $\alpha\in[1/2,1]$, then ${\mathcal F}_2''\le 0$ in $J_\alpha$. 
Therefore we deduce that 
$\alpha$-log-concavity is preserved by solutions of problem~(SH) if and only if $\alpha\in[1/2,1]$. 
Thus assertion~(2) follows. 

We prove assertion~(3). 
It suffices to consider the case of $a<\infty$. 
Indeed, assertion~(3) follows assertion~(1) if $a=\infty$. 
Since $ah$ is the inverse function of $H_a$ in ${\mathbb R}$, 
by \eqref{eq:2.5} we have 
$$
{\mathcal F}_1(z)=\frac{h''(z)}{h'(z)}=-\frac{1}{2}z,\quad {\mathcal F}_2(z)=\kappa a^{p-1}\frac{h(z)^p}{h'(z)},
\quad\mbox{for}\quad z\in J_{H_a}={\mathbb R}. 
$$
Then ${\mathcal F}_1$ is concave in ${\mathbb R}$. 
Furthermore, it follows from \eqref{eq:2.5} that
\begin{align*}
{\mathcal F}_2'(z) & =\kappa a^{p-1}ph(z)^{p-1}-\kappa a^{p-1}h(z)^p\frac{h''(z)}{h'(z)^2}=\kappa a^{p-1}ph(z)^{p-1}+\kappa\frac{a^{p-1}}{2}z\frac{h(z)^p}{h'(z)},\\
{\mathcal F}_2''(z) & =\kappa a^{p-1}p(p-1)h(z)^{p-2}h'(z)
+\kappa\frac{a^{p-1}}{2}\frac{h(z)^p}{h'(z)}
 +\kappa\frac{a^{p-1}}{2}pzh(z)^{p-1}+\kappa\frac{a^{p-1}}{4}z^2\frac{h(z)^p}{h'(z)},
 \end{align*}
for $z\in{\mathbb R}$, and we see that ${\mathcal F}_2''\le 0$ in $[0,\infty)$. 
On the other hand, it follows from \eqref{eq:2.7} that 
$$
h'(z)\ge -\frac{1}{2}zh(z)\quad\mbox{for}\quad z\in{\mathbb R}. 
$$
Then we have
\begin{equation}
\label{eq:6.15}
{\mathcal F}_2''(z)
\le -\kappa\frac{a^{p-1}p}{2}(p-2)zh(z)^{p-1}+\kappa\frac{a^{p-1}}{2}\frac{h(z)^p}{h'(z)}+\kappa\frac{a^{p-1}}{4}z^2\frac{h(z)^p}{h'(z)}
\end{equation}
for $z\in(-\infty,0)$, 
and we see that ${\mathcal F}_2''\le 0$ in $(-\infty,0)$ if $p\ge 2$. 
We deduce that 
${\mathcal F}_2''(z)\le 0$ in ${\mathbb R}$ if $p\ge 2$.  

We consider the case of $1<p<2$. 
Since $H_1$-concavity is preserved by DHF in $\Omega$, 
by Proposition~\ref{Theorem:6.1}~(3) we see that $(\log h)'$ is $(-1)$-concave in ${\mathbb R}$, 
which together with \eqref{eq:2.5} implies that 
\begin{equation}
\label{eq:6.16}
\begin{split}
0 & \le\left(\frac{h(z)}{h'(z)}\right)''
=\left(1-\frac{h(z)h''(z)}{h'(z)^2}\right)'
=\left(\frac{zh(z)}{2h'(z)}\right)'
=\frac{h(z)}{2h'(z)}+\frac{z}{2}-\frac{zh(z)h''(z)}{2h'(z)^2}\\
 & =\frac{h(z)}{2h'(z)}+\frac{z}{2}+\frac{z^2h(z)}{4h'(z)}
\quad\mbox{for}\quad z\in{\mathbb R}.
\end{split}\end{equation}
By \eqref{eq:6.15} and \eqref{eq:6.16} we obtain
\begin{equation*}
\begin{split}
{\mathcal F}_2''(z)
 & \le -\kappa\frac{a^{p-1}}{2}|z|h(z)^{p-1}+\kappa\frac{a^{p-1}}{2}\frac{h(z)^p}{h'(z)}+\kappa\frac{a^{p-1}}{4}z^2\frac{h(z)^p}{h'(z)}\\
 & =\kappa a^{p-1}h(z)^{p-1}
\left[\frac{z}{2}+\frac{h(z)}{2h'(z)}+\frac{z^2h(z)}{4h'(z)}\right]\le 0
\end{split}
\end{equation*}
for $z\in(-\infty,0)$. 
Thus ${\mathcal F}_2''\le 0$ in ${\mathbb R}$ if $1<p<2$. 
Therefore, each of ${\mathcal F}_1$ and ${\mathcal F}_2$ is concave in~${\mathbb R}$, 
and $H_a$-concavity is preserved by solutions of problem~(SH). 
Thus assertion~(3) follows, and the proof of Theorem~\ref{Theorem:6.10} is complete.
$\Box$\vspace{5pt}
\newline
By a direct consequence of Theorem~\ref{Theorem:6.10} with $\kappa=0$ we have:
\begin{corollary}
\label{Theorem:6.11}
Let $\Omega$ be a convex domain in~${\mathbb R}^n$ and $\alpha\in{\mathbb R}$.
Then $\alpha$-log-concavity is preserved by DHF in~$\Omega$ if and only if $\alpha\in[1/2,1]$.
\end{corollary}
%
\subsection{Porous medium equation and parabolic $p$-Laplace equation}\label{subsection:6.3}
In this subsection
we discuss the preservation of concavity property by solutions of the Cauchy--Dirichlet problem
for the porous medium equation and the parabolic $p$-Laplace equation. 

Consider the Cauchy--Dirichlet problem
for the porous medium equation 
\begin{equation}
\tag{PM}
\left\{
\begin{array}{ll}
\partial_t u=\Delta (u^m) & \quad\mbox{in}\quad\Omega\times(0,\infty),\vspace{3pt}\\
u=0 & \quad\mbox{on}\quad\partial\Omega\times(0,\infty)\mbox{ if }\partial\Omega\not=\emptyset,\vspace{3pt}\\
u(\cdot,0)=\phi\ge 0 & \quad\mbox{in}\quad\Omega,
\end{array}
\right.
\end{equation}
where $m>1$ and $\phi\in L^\infty(\Omega)$. 
The porous medium equation provides a simple model in many physical situations  {and the
preservation of power concavity} by solutions of problem~(PM) has been studied in many papers,
see e.g. \cites{BV, CVW, CW3, INS, IS02, DH1, DH2, DHL, LV} and references therein. 
(See also the monograph~\cite{V} for the the existence and the uniqueness of solutions 
of problem~(PM) and for related topics.)
Among others, 
B\'enilan and V\'azquez~\cite{BV}*{Theorem~1} proved the preservation of $(m-1)$\,-concavity by solutions of (PM) in the case $\Omega={\mathbb R}$. 
Furthermore, Daskalopoulos, Hamilton, and Lee \cite{DHL}*{Theorem~1.1}
proved that $(m-1)/2$\,-concavity is preserved by solutions of problem~(PM) in the case $\Omega={\mathbb R}^n$ with $n\ge 1$ 
under the following non-degeneracy condition on  $\phi_m:=\phi^{m-1}$:
\begin{equation}
\label{eq:6.17}
\left\{
\begin{array}{l}
  \mbox{$\phi_m$ has a compact support $D$};\vspace{3pt}\\
  \mbox{$\phi_m$ is smooth on $D$ and $\displaystyle{\min_{D}\,(\phi_m+|\nabla\phi_m|^2)>0}$}. 
\end{array}
\right.
\end{equation}
On the other hand, 
it was also proved in \cite{INS}*{Theorem~1.3} that $(m-1)/2$\,-concavity is preserved by solutions 
of problem~(PM) in the case where $\Omega$ is a smooth, bounded, and convex domain in~${\mathbb R}^n$ 
and $\phi\in BC_0(\overline{\Omega})$ 
with $\phi>0$ in $\Omega$. 
Then, similarly to DHF, the following question naturally arises:
\begin{itemize}
  \item[({\bf Q5})]  
  if $\alpha\not=(m-1)/2$, is $\alpha$-concavity preserved by solutions of problem~(PM)? 
\end{itemize}
(See \cite{V}*{Page 520} for related questions.) 
A partial negative answer to question~({\bf Q5}) was given in \cite{IS02}*{Theorem~1.1} for the case of $n\ge 2$. 
Recently, question~({\bf Q5}) was negatively and perfectly resolved by  {\cite{CW3}*{Theorem~1.1}
in the case where $\Omega={\mathbb R}^2$, but it remains open in the other cases.}

Although problem~(PM) is outside the framework of the previous sections, 
as a direct application of the argument in the proof of Proposition~\ref{Theorem:4.1}, 
we obtain a negative answer to question~({\bf Q5}) for all dimension.
\begin{theorem}
\label{Theorem:6.12}
Let $n\ge 1$, $m>1$, and $-\infty<\alpha<(m-1)/2$.  
\begin{itemize}
  \item[{\rm (1)}] 
  Let $\Omega$ be a bounded convex domain in ${\mathbb R}^n$. 
  Then there exists $\phi\in {\mathcal C}_\Omega[\Phi_\alpha]\cap BC_0(\overline{\Omega})$ with $\phi>0$ in $\Omega$ such that 
  the corresponding solution of problem~{\rm (PM)} is not $\alpha$-concave in $\Omega$ for some $t>0$.  
  \item[{\rm (2)}]  
  Let $D$ be compact and convex in ${\mathbb R}^n$. 
  There exists $\phi\in {\mathcal C}_{{\mathbb R}^n}[\Phi_\alpha]$ 
  such that $\phi_m=\phi^{m-1}$ satisfies condition~\eqref{eq:6.17} 
  and the corresponding solution of problem~{\rm (PM)} with $\Omega={\mathbb R}^n$ is not $\alpha$-concave in ${\mathbb R}^n$ for some~$t>0$.
\end{itemize}
\end{theorem}
{\bf Proof.}
We prove assertion~(1). 
Thanks to the boundedness of $\Omega$ and \eqref{eq:1.2}, 
we can assume, without loss of generality, that 
$$
\langle e_1,x\rangle>0\quad\mbox{for}\quad x\in\Omega. 
$$
Let $-\infty<\alpha<(m-1)/2$. Set 
\begin{equation}
\label{eq:6.18}
\phi(x):=
\left\{
\begin{array}{ll}
\langle e_1,x\rangle^{\frac{1}{\alpha}} & \mbox{if}\quad \alpha\not=0,\vspace{5pt}\\
e^{\langle e_1,x\rangle} & \mbox{if}\quad\alpha=0,
\end{array}
\right.
\end{equation}
for $x\in\Omega$. Then $\phi$ is $\alpha$-concave in $\Omega$. 

Let $u$ be the solution of (PM) with the initial datum $\phi$. 
It follows from the comparison principle that $u>0$ in $\Omega\times[0,\infty)$. 
Then parabolic regularity theorems (see \cite{LSU}*{Chapter~V}) imply that $u\in C^{2;1}(\Omega\times[0,\infty))$. 
Assume that $u(\cdot,t)$ is $\alpha$-concave in $\Omega$ for all $t>0$, 
that is,  
$$
\Phi_\alpha(u((1-\lambda)x+\lambda y,t))
-(1-\lambda)\Phi_\alpha(u(x,t))-\lambda\Phi_\alpha(u(y,t))\ge 0
$$
for $x$, $y\in\Omega$, $t>0$, and $\lambda\in(0,1)$. 
Furthermore, by \eqref{eq:6.18} we have
$$
\Phi_\alpha(u((1-\lambda)x+\lambda y,0))
-(1-\lambda)\Phi_\alpha(u(x,0))-\lambda\Phi_\alpha(u(y,0))=0
$$
for $x$, $y\in\Omega$. 
Then we observe that
\begin{equation}
\label{eq:6.19}
\left(\partial_t\Phi_\alpha(u)\right)((1-\lambda)x+\lambda y,0))
-(1-\lambda)\left(\partial_t\Phi_\alpha(u)\right)(x,0)-\lambda\left(\partial_t\Phi_\alpha(u)\right)(y,0)\ge 0
\end{equation}
for $x$, $y\in\Omega$. 

On the other hand, 
it follows from (PM) and \eqref{eq:6.18} that 
$$
\left(\partial_t\Phi_\alpha(u)\right)(z,0)=\phi(z)^{\alpha-1}(\Delta\phi^m)(z)\\
=
\left\{
\begin{array}{ll}
\displaystyle{\frac{m}{\alpha}\left(\frac{m}{\alpha}-1\right)
\langle e_1,z\rangle^{-1+\frac{m-1}{\alpha}}} & \mbox{if}\quad \alpha\not=0,\vspace{7pt}\\
m^2e^{(m-1)\langle e_1,z\rangle}\quad & \mbox{if}\quad \alpha=0,
\end{array}
\right.
$$
for $z\in\Omega$. 
Then relation~\eqref{eq:6.19} leads 
the concavity of $\left(\partial_t\Phi_\alpha(u)\right)(\cdot,0)$, and we have
$$
0\le -1+\frac{m-1}{\alpha}\le 1\quad\mbox{if}\quad\alpha\not=0,
\qquad
m^2(m-1)^2\le 0 \quad\mbox{if}\quad \alpha=0.
$$
These are both contradictions. 
Thus $\alpha$-concavity is not preserved by the solution~$u$. 
Finally, 
similarly to Lemma~\ref{Theorem:3.4}, 
we approximate $\phi$ by positive functions belonging to ${\mathcal C}_\Omega[\Phi_\alpha]\cap BC_0(\overline{\Omega})$ 
to obtain assertion~(1). 
Applying a similar argument with $\Omega$ replaced by $D$,
we obtain assertion~(2). Thus Theorem~\ref{Theorem:6.12} follows. 
$\Box$\vspace{5pt}

Finally, we consider the Cauchy--Dirichlet problem for the parabolic $p$-Laplace equation
\begin{equation}
\tag{PP}
\left\{
\begin{array}{ll}
\partial_t u=\mbox{div}\,(|\nabla u|^{p-2}\nabla u) & \quad\mbox{in}\quad\Omega\times(0,\infty),\vspace{3pt}\\
u=0 & \quad\mbox{on}\quad\partial\Omega\times(0,\infty)\mbox{ if }\partial\Omega\not=\emptyset,\vspace{3pt}\\
u(\cdot,0)=\phi\ge 0 & \quad\mbox{in}\quad\Omega,
\end{array}
\right.
\end{equation} 
where $p>2$ and $\phi\in L^\infty(\Omega)$. 
See e.g. \cite{DiB} for the existence, the uniqueness, and the regularity of solutions of problem~(PP). 
Little is known concerning the preservation of concavity properties by solution of the parabolic $p$-Laplace equation, 
and the only available result is in \cite{L}, 
where Lee proved that 
$(p-2)/p$\,-concavity is preserved by solutions of problem~(PP) with $\Omega={\mathbb R}^n$ 
under the following non-degeneracy condition on $\varphi_p:=\phi^{(p-2)/(p-1)}$:
\begin{equation}
\label{eq:6.20}
\left\{
\begin{array}{l}
  \mbox{$\varphi_p$ has a compact support $D$};\vspace{3pt}\\
  \mbox{$\varphi_p$ is smooth on $D$ and $\displaystyle{\min_{D}\,(\varphi_p+|\nabla\varphi_p|)>0}$}. 
\end{array}
\right.
\end{equation}
On the other hand, 
by the same argument as in the proof of Theorem~\ref{Theorem:6.12}, 
we have the following result for problem~(PP).
\begin{theorem}
\label{Theorem:6.13}
Let $n\ge 1$, $p>2$, and $-\infty<\alpha<(p-2)/p$. 
\begin{itemize}
  \item[{\rm (1)}] 
  Let $\Omega$ be a bounded convex domain in ${\mathbb R}^n$. 
  Then there exists $\phi\in {\mathcal C}_\Omega[\Phi_\alpha]\cap BC_0(\overline{\Omega})$ with $\phi>0$ in $\Omega$ such that 
  the corresponding solution of problem~{\rm (PP)} is not $\alpha$-concave in $\Omega$ for some $t>0$.  
  \item[{\rm (2)}]  
  Let $D$ be convex and compact in ${\mathbb R}^n$. 
  There exists $\phi\in {\mathcal C}_{{\mathbb R}^n}[\Phi_\alpha]$ such that 
  $\varphi_p=\phi^{(p-2)/(p-1)}$ satisfies \eqref{eq:6.20} and 
  the corresponding solution of problem~{\rm (PP)} with $\Omega={\mathbb R}^n$ is not $\alpha$-concave in ${\mathbb R}^n$ for some $t>0$
\end{itemize}
\end{theorem}
{\bf Proof.}
We prove assertion~(1). 
Thanks to the boundedness of $\Omega$ and \eqref{eq:1.2}, 
we can assume, without loss of generality, that 
$$
\langle e_1,x\rangle>0\quad\mbox{for}\quad x\in\Omega. 
$$
Let $-\infty<\alpha<(p-2)/p$ and let $\phi$ be as in the proof of Theorem~\ref{Theorem:6.12}. 
Since $\nabla\phi\not=0$ in $\Omega$, 
for any $x\in\Omega$, 
the corresponding solution~$u$ of problem~(PP) is $C^{2;1}$-smooth in a neighborhood of $(x,0)$ in $\Omega\times[0,\infty)$. 
(See e.g. \cite{DiB} and \cite{LSU}*{Chapters~IV and V}.)
Furthermore, 
\begin{equation}
\label{eq:6.21}
\phi(z)^{\alpha-1}\mbox{div}\,(|\nabla\phi|^{p-2}\nabla\phi)(z)\\
 =
\left\{
\begin{array}{ll}
\displaystyle{\frac{p-1}{|\alpha|^{p-2}\alpha}\frac{1-\alpha}{\alpha}
\langle e_1,z\rangle^{\frac{p-2}{\alpha}-(p-1)}}\quad & \mbox{if}\quad \alpha\not=0,\vspace{5pt}\\
(p-1)e^{(p-2)\langle e_1,z\rangle} & \mbox{if}\quad \alpha=0,
\end{array}
\right.
\end{equation}
for $z\in\Omega$. 

Assume that $u(\cdot,t)$ is $\alpha$-concave in $\Omega$ for all $t>0$. 
Applying the same argument as in the proof of  Theorem~\ref{Theorem:6.12}, 
by \eqref{eq:6.21} we see $\alpha\not=0$ and 
$$
0\le \frac{p-2}{\alpha}-(p-1)\le 1.
$$
This implies that $\alpha\ge (p-2)/p$, which yields a contradiction. 
Thus $\alpha$-concavity is not preserved by the solution~$u$. 
Finally, 
we approximate $\phi$ by functions belonging to ${\mathcal C}_\Omega[\Phi_\alpha]\cap BC_0(\overline{\Omega})$ 
to obtain assertion~(1). 
Applying a similar argument with $\Omega$ replaced by $D$,
we obtain assertion~(2) is similar. Thus Theorem~\ref{Theorem:6.13} follows. 
$\Box$
\bigskip

{\bf Statements and Declarations}:
\begin{itemize}
  \item 
  'Declarations of interest: none'
  \item 
  'Funding Source Declarations: The second author was partially financed by INdAM through a GNAMPA Project'.
  \item 
  'Competing Interests Declarations: The author has no relevant financial or non-financial interests to disclose and has no competing interests to declare that are relevant to the content of this article'.
  \item
  'Conflict of interest statement: the authors state that there is no conflict of interest'.
\end{itemize}

 
\end{document}